\documentclass[oneside,reqno]{amsart}
\usepackage{amsmath}
\usepackage{latexsym}
\usepackage{amsfonts,amssymb,amsxtra}

\usepackage{amscd}

\usepackage{latexsym}

\usepackage{mathptmx}

\newtheorem{theorem}{Theorem}[section]

\newtheorem{claim}[theorem]{Claim}

\newtheorem{corollary}[theorem]{Corollary}

\newtheorem{lemma}[theorem]{Lemma}

\newtheorem{proposition}[theorem]{Proposition}

\theoremstyle{definition}

\newtheorem{definition}[theorem]{Definition}
\newtheorem{example}[theorem]{Example}

\newtheorem{remark}[theorem]{Remark}

\newcommand\End{\operatorname{End}}

\newcommand\supp{\operatorname{supp}}

\newcommand\red{\operatorname{red}}

\newcommand\tr{\operatorname{tr}}

\newcommand\spa{\operatorname{span}}

\newcommand\ad{\operatorname{ad}}

\newcommand\rad{\operatorname{rad}}

\newcommand\rank{\operatorname{rank}}

\newcommand\der{\operatorname{Der}}

\newcommand\oder{\operatorname{Oder}}

\newcommand\sll{\operatorname{sl}}

\newcommand\gl{\operatorname{gl}}

\newcommand\dia{\operatorname{diag}}

\newcommand\fin{\operatorname{fin}}

\newcommand\sh{\operatorname{sh}}

\newcommand\lng{\operatorname{lg}}

\newcommand\ex{\operatorname{ex}}

\newcommand\vvp{\varphi}

\usepackage{color}
\RequirePackage[normalem]{ulem} 
\RequirePackage{color}\definecolor{RED}{rgb}{1,0,0}\definecolor{BLUE}{rgb}{0,0,1} 
\providecommand{\DIFaddbegin}{} 
\providecommand{\DIFaddend}{} 
\providecommand{\DIFdelbegin}{} 
\providecommand{\DIFdelend}{} 

\begin{document}

\today

\enskip

\title{Locally loop algebras and locally affine Lie algebras}

\author{Jun Morita}
\address[Jun Morita]{Institute of Mathematics\\
University of Tsukuba\\Tsukuba, Ibaraki, 305-8571 Japan}
\email{morita@math.tsukuba.ac.jp}
\author{Yoji Yoshii}
\address[Yoji Yoshii]
{Iwate University\\
Ueda 3-18-33 Morioka, Iwate, Japan 020-8550}
\email{yoshii@iwate-u.ac.jp}
\DIFdelbegin 
\DIFdelend 
\DIFaddbegin \thanks{Jun Morita gratefully acknowledges partial financial support from a
Grant-in-aid for Science Research (MONKASHO KAKENHI) in Japan.}
\DIFaddend 
\date{}

\begin{abstract}
In this study, we  investigate a new class of Lie algebras, i.e.,
tame locally extended affine Lie algebras of nullity $1$, which
are an infinite-rank analog of affine Lie algebras. This type of
algebra is called a locally affine Lie algebra. A certain ideal of
a locally affine Lie algebra, called a core, is a universal
extension of a local version of a loop algebra, which is  called a
locally loop algebra. We classify locally loop algebras and
locally affine Lie algebras.

\end{abstract}

\maketitle

Throughout this study, $F$ is a field of characteristic 0. All of
the algebras are assumed to be unital, except the Lie algebras.
The tensor products are over $F$.

\enskip

\section{Introduction}
\label{sec:intro}

Historically, root systems have played very important roles in Lie
theory and many other areas. To obtain a root system, we usually
need a certain ad-diagonalizable subalgebra $\mathcal H$ of a Lie
algebra $\mathcal L$ over $F$. Then, we have a decomposition:
$$\mathcal L = \bigoplus_{\xi \in \mathcal H^*} \mathcal L_\xi,$$
where $\mathcal H^*$ is the dual space of $\mathcal H$ and
$\mathcal L_\xi = \{ x \in \mathcal L \mid [ h , x ] = \xi(h)x\
\mbox{for all}\ h \in \mathcal H \}$. An element $\xi \in \mathcal
H^*$ is called a {\bf root} if $\mathcal L_\xi \not= 0$, and the
{\bf set of roots} is defined by
$$R = \{ \xi \in \mathcal H^* \mid \mathcal L_\xi \not= 0 \}.$$
The subspace $\mathcal L_\xi$ is called the {\bf root space} of
$\xi$ and the direct sum above is called the {\bf root space
decomposition} of $\mathcal L$ associated with $\mathcal H$. In
many cases, a root has its own length, which may come from a
symmetric invariant bilinear form $\mathcal B$ on $\mathcal L$.
Therefore, it is natural to consider a triplet $(\mathcal
L,\mathcal H,\mathcal B)$ in general.

\bigskip

Such a triplet $(\mathcal L,\mathcal H, \mathcal B)$ is called a
{\bf locally extended affine Lie algebra} if 
conditions (A1) -- (A4)
are satisfied, as follows:

\begin{itemize}
\item[(A1)] $\mathcal H$ is ad-diagonalizable and
self-centralizing, i.e.,
$$\mathcal L =\bigoplus_{\xi\in \mathcal H^*} \mathcal L _\xi
\quad\text{and}\quad \mathcal H = \mathcal L _0,$$

\item[(A2)] ${\mathcal B}$ is nondegenerate,

\item[(A3)] $\ad x\in\End_F \mathcal L $ is locally nilpotent for
all $\xi\in R^\times$ and all $x\in \mathcal L _\xi$, where
$R^\times = \{ \xi \in R \mid \mathcal{B}(t_\xi,t_\xi) \not= 0 \}$
and $t_\xi$ is an element of $\mathcal{H}$ such that $\xi(h) =
\mathcal{B}(t_{\xi},h)$ for all $h \in \mathcal{H}$.

\item[(A4)] $R^\times$ is irreducible, i.e., there is no
nontrivial partition, $R^\times = R_1 \cup R_2$, of $R^\times$
such that $\mathcal B(t_{\xi_1},t_{\xi_2}) = 0$ for all $\xi_1 \in
R_1,\ \xi_2 \in R_2$.

\end{itemize}

\smallskip

A locally extended affine Lie algebra is abbreviated as a {\bf
LEALA}. We note that the class of LEALAs contains:

(1) all finite dimensional split simple Lie algebras,

(2) all affine Lie algebras,

(3) all locally finite split simple Lie algebras,

(4) all extended affine Lie algebras, and

(5) all Heisenberg Lie algebras (as null systems).

\noindent We can refer to them all as LEALAs and study them
uniformly.
The class of LEALAs is assumed to be the best (or the widest)
class of $(\mathcal{L},\mathcal{H},\mathcal{B})$ in the following
sense.

If we fix $\xi \in R^\times$ and select  $a \in F^\times$ such that
$a \mathcal{B}(t_\xi,t_\xi) \in \mathbb{Q}_{> 0}$, as well as
defining $( \cdot , \cdot  )$ as a symmetric bilinear form on $V =
\sum_{\xi \in R}\ \mathbb{Q} \xi \subset \mathfrak{h}^*$ by
$$\left( \sum_{\xi \in R} p_\xi \xi \ ,\ \sum_{\eta \in R} q_\eta \eta \right) =
a\, \mathcal{B} \left( \sum_{\xi \in R} p_\xi t_\xi \ ,\
\sum_{\eta \in R} q_\eta t_\eta \right)$$ for $p_\xi, q_\eta \in
\mathbb{Q}$, then we find that the properties of the form $( \cdot
, \cdot )$ are actually  $\Bbb Q$-valued and that
$$
\text{$( \cdot , \cdot )$ is positive semi-definite, while $( \xi
, v ) = 0$ for all $\xi \in R^0$ and all $v \in V$},
$$
where $R^0 = \{ \xi \in R \mid (\xi,\xi) = 0 \}$. This property is
usually called the {\bf Kac Conjecture}, which was proved
 for EALAs
(defined below)
 in \cite{AABGP}
and for LEALAs in \cite{MY}.

\smallskip

Note that $R^\times = \{ \xi \in R \mid (\xi,\xi) \not= 0 \}$. An
element of $R^\times$ is called an {\bf anisotropic root}, and an
element of $R^0$ is called an {\bf isotropic root}. We can also
refer to $R^\times$ as the {\bf root system}. Note that $R^\times$
is a finite root system and $R^0 = \{ 0 \}$ when $\mathcal L$ is a
finite-dimensional split simple Lie algebra (cf. \cite{Bo}), while
$R^\times$ is an affine root system and $R^0 = \mathbb{Z} \xi_0$
for some $\xi_0 \in \mathcal H^*$ when $\mathcal L$ is an affine
Lie algebra (cf. \cite{Ma}). A LEALA is called an {\bf extended
affine Lie algebra} ({\bf EALA}), when $\dim \mathcal H < \infty$.

\smallskip

The {\bf nullity} of a LEALA
 is defined as the rank of the additive group
 generated by $R^0$.
In particular, we only use the term `nullity' when the additive group is free
 (see Remark \ref{onlyforfree}).
The {\bf core} $\mathcal L_c$ of a LEALA $\mathcal L$
 is defined as
 the subalgebra of $\mathcal L$ generated by $\mathcal L_\xi$
 for all $\xi\in R^\times$.
In fact, $\mathcal L_c$ is an ideal of $\mathcal L$, which is
obtained by the Kac Conjecture. If the centralizer $C_{\mathcal
L}(\mathcal L_c)$ of $\mathcal L_c$ in $\mathcal L$ is contained
in $\mathcal L_c$, this LEALA $\mathcal L$ is referred to as {\bf
tame}. The core $\mathcal L_c$ modulo of its center $Z(\mathcal
L_c)$, i.e., the quotient Lie algebra $\mathcal L_c/Z(\mathcal
L_c)$, is called the {\bf centerless core} of $\mathcal L$.

 \smallskip

Previously, we classified LEALAs of nullity $0$ in \cite{MY}. The
second simplest class comprises LEALAs of nullity $1$. The main
aim of the present study is to classify the class of tame LEALAs
of nullity $1$, which we call a {\bf locally affine Lie algebra}
({\bf LALA}).

\smallskip

The centerless core
$$L := \mathcal{L}_c/Z(\mathcal{L}_c)$$ of a LALA $\mathcal{L}$
is a local version of a loop algebra, which we call a {\bf locally
loop algebra}. In fact, we show that a locally loop algebra is a
directed union of loop algebras. We also show that the core
$\mathcal{L}_c$ of a LALA $\mathcal{L}$ is a universal covering of
a locally loop algebra $L$.  (These classifications were also shown
by Neeb \cite[Cor. 3.13]{N2} in a different manner.) 
Thus, we may
say that a LALA is also a local analog of an affine Lie algebra.
However,
 a LALA $\mathcal{L}$ has a more complex structure in
a complement of the core $\mathcal{L}_c$, such as $\mathcal{L} =
\mathcal{L}_c \oplus D$. Note that for an affine Lie algebra, the
complement $D$ is simply a $1$-dimensional space spanned by the
degree derivation. However, for a LALA, the corresponding
complement $D$ is rather large in general. Due to tameness, $D$
can be embedded into $\der_F \mathcal{L}_c$. Then, $d \in \der_F
\mathcal{L}_c$ induces a derivation of $L =
\mathcal{L}_c/Z(\mathcal{L}_c)$ since $d(Z(\mathcal{L}_c)) \subset
Z(\mathcal{L}_c)$, and we see that $d(x) \in Z(\mathcal{L}_c)$ for
$x \in \mathcal{L}_c$ implies that $d = 0$ in $\der_F
\mathcal{L}_c$ since $\mathcal{L}_c = [ \mathcal{L}_c ,
\mathcal{L}_c ]$. Therefore, we can find $D \subset \der_F
\mathcal{L}_c \subset \der_F L$, but also $D \subset \oder_F
\mathcal{L}_c \subset \oder_F L$, where $\oder_F( \cdot ) =
\der_F( \cdot ) / \ad( \cdot )$, which comprises the outer
derivations. There is a unique maximal choice of $D$ in $\oder_F
L$, such as $D^{max}$ in this case, and there are many minimal
choices of $D$ in $\oder_F L$, such as $D(p)$ with a specific
diagonal matrix $p$. Thus, a homogeneous space $D$ such that $D(p)
\subset D \subset D^{max}$ leads to our classification.  Thus, the
classification of LALAs is obtained by saying that any homogeneous
subalgebra $\mathcal{L}$ of $\mathcal{L}^{max} := \mathcal{L}_c
\oplus D^{max}$ that satisfies
$$\mathcal{L}(p) := \mathcal{L}_c \oplus D(p) \subset \mathcal{L} \subset \mathcal{L}^{max}$$
is a LALA, and that any LALA can be obtained in this manner. We
roughly explain the LALAs of type $\text A^{(1)}_\frak I$ and
$\text {C}^{(2)}_\frak I$ to obtain a better understanding.

\smallskip

Let ${\frak I}$ be an index set. We suppose that $\mathfrak{I}$ is
any index set, i.e., $\mathfrak{I}$ can be finite or infinite. Let
$M_\frak{I}(F) = \left\{ (a_{ij})_{i,j \in \frak{I}} \mid a_{ij}
\in F \right\}$ be the vector space of matrices of size $\frak I$,
and let $T_\frak{I} = T_\frak{I}(F) = \left\{ (a_{ij}) \in
M_\frak{I}(F) \mid a_{ij} = 0\ \mbox{for}\ i \not= j \right\}$ be
the subspace of $M_{\frak I}(F)$, which comprises all diagonal
matrices.

\smallskip

[$\text{A}_\frak I^{(1)}$]\ \ First, we explain the untwisted type
$\text{A}_\frak I^{(1)}$. Let $\sll_\frak{I}(F)$ be the subspace
of $M_{\frak I}(F)$ that comprises trace $0$ matrices with only
finitely many nonzero entries. Let $F[t^{\pm 1}]$ be the algebra
of Laurent polynomials, and let $\sll_\frak I(F[t^{\pm 1}])$ be
the Lie algebra $\sll_\frak{I}(F)\otimes F[t^{\pm 1}]$. For
example, if $\frak{I} = \mathbb{N}$ (the natural numbers), then we
see that
$$\sll_\mathbb{N}(F[t^\pm]) = \bigcup_{n=2}^\infty \sll_n(F[t^\pm])
= \bigcup_{n} \left(
\begin{array}{c|c}
\sll_n(F[t^\pm]) & O\\ \hline O & O
\end{array}
\right) \ .$$ We refer to $\sll_\frak I(F[t^{\pm 1}])$ as a
locally loop algebra of type $\text A^{(1)}_\frak I$, which is
simply an infinite-rank analog of a loop algebra
$\sll_{\ell+1}(F[t^{\pm 1}])$ of type $\text A^{(1)}_\ell$. We use
the following conventions.
$$
\begin{cases}
\text{$\sll_\frak{I}$ has type $\text A_\frak{I}$}
&\text{if $\frak{I}$ is an infinite index set}, \\
\text{$\sll_\frak{I} = \sll_{\ell + 1}$ has type $\text A_\ell$}
&\text{if $\frak{I}$ is a finite index set that comprises $\ell +
1$ elements}.
\end{cases}
$$
As in the case of $\sll_{\ell+1}(F[t^{\pm 1}])$, a universal
covering $\sll_\frak I(F[t^{\pm 1}])\oplus Fc$ of $\sll_\frak
I(F[t^{\pm 1}])$ exists, where $Fc$ is the $1$-dimensional center.
Then we can construct the Lie algebra
\begin{equation}
\mathcal L^{ms} :=\sll_\frak I(F[t^{\pm 1}])\oplus Fc\oplus
Fd^{(0)}, \label{minimalst}
\end{equation}
where $\displaystyle{d^{(0)}=t\frac{d}{dt}}$ is the degree
derivation. This $\mathcal L^{ms}$ is the simplest example of a
LALA, which is called a {\bf minimal standard LALA} of type $\text
A^{(1)}_\frak I$. In contrast to the affine Lie algebra case,
there are more examples of type $\text A^{(1)}_\frak I$, which are
obtained by adding diagonal derivations of $\sll_\frak I(F[t^{\pm
1}])$, and we explain these as follows. First, note that
$$\sll_\frak I(F)+T_{\frak I}$$
is a Lie algebra with center $F\iota$, where $\iota=\iota_{\frak
I} = (\delta_{ij})_{i,j \in \mathfrak{I}}$ is in $T_{\frak I}$.
Let
\begin{equation}\label{newdefleala0}
\mathcal A_\frak I:=\big(\sll_\frak I(F)+T_{\frak I}\big)/F\iota
\end{equation}
be the quotient Lie algebra. We identify the subalgebra
$$\overline{\sll_\frak I(F)}=\big(\sll_\frak I(F)+F\iota\big)/F\iota$$
of $\mathcal A_\frak I$ with $\sll_\frak I(F)$. Consider the Lie
algebra $\mathcal A_\frak I\otimes F[t^{\pm 1}]$. We construct the
Lie algebra
\begin{equation}
\hat{\mathcal A}_\frak I:=\mathcal A_\frak I\otimes F[t^{\pm
1}]\oplus Fc \oplus Fd^{(0)}
 \label{maxreali}
\end{equation}
as described in \eqref{minimalst}, which contains $\mathcal
L^{ms}$.

\begin{theorem}
$\mathcal{L}^{max} :=\hat{\mathcal A}_\frak I$ is a {\bf maximal
LALA} of type $\text A^{(1)}_\frak I$, i.e., any LALA of type
$\text A^{(1)}_\frak I$ is a subalgebra of $\hat{\mathcal A}_\frak
I$. In addition, any LALA of type $\text A^{(1)}_\frak I$ contains
a LALA
\begin{equation}
\mathcal L(p): =\sll_\frak I(F[t^{\pm 1}])\oplus Fc\oplus
F(p+d^{(0)}) \label{minreali}
\end{equation}
for some $p\in T_\frak I$.
\end{theorem}

This $\mathcal L(p)$ is called a {\bf minimal LALA determined by
$p$}. In general, we note that $\mathcal L(p)$ may be isomorphic to
$\mathcal L^{ms}=\mathcal L(0)$, but not always isomorphic to
$\mathcal L^{ms}$ (see Example \ref{counterexofst}). Note that a
LALA $\mathcal L$ of type $\text{A}_\frak I^{(1)}$ has a
decomposition $\mathcal{L} = \mathcal{L}_c \oplus D$ for a
homogeneous complement $D$ that satisfies $D(p) = F(p + d^{(0)})
\subset D \subset D^{max}$ and $\mathcal{L}_c \oplus D^{max} =
\hat{\mathcal{A}}_\frak I$.

\medskip

[$\text{C}_\frak I^{(2)}$]\ \ Next, we explain the twisted type
$\text {C}^{(2)}_\frak I$. Let $s=
\begin{pmatrix}
   0   &-\iota       \\
  \iota   & 0      \\
\end{pmatrix}$
be the matrix of size $2\frak I$, where $\iota=\iota_{\frak I}$,
as described above. Define an automorphism $\sigma$ of period $2$
on $\sll_{2\frak I}(F)+T_{2\frak I}$ by
$$\sigma(x)=sx^Ts$$
for $x\in \sll_{2\frak I}(F)+T_{2\frak I}$, where $x^T$ is the
transpose of $x$. Let $\text{sp}_{2\frak I}(F)$ be the fixed
subalgebra of $\sll_{2\frak I}(F)$ by $\sigma$, which is of type
$\text C_\frak I$. Let $\frak s$ be the $(-1)$-eigenspace of
$\sigma$ such that
$$\sll_{2\frak I}(F)=\text{sp}_{2\frak I}(F)\oplus\frak s.$$
Moreover, let $T^{+}$ be the $1$-eigenspace and $T^{-}$ is the
$(-1)$-eigenspace of $T_{2\frak I}$ relative to $\sigma$. Note
that $T_{2\frak I} = T^+ \oplus T^-$, and thus
$$\sll_{2\frak I}(F)+T_{2\frak I}=\big(\text{sp}_{2\frak I}(F)+T^+\big)\oplus\big(\frak s+T^-\big).$$
In addition, note that $F\iota_{2\frak I}$ is $\sigma$-invariant
and $F\iota_{2\frak I}\subset T^-$. Let
$$\mathcal A_{2\frak I}:=\big(\sll_{2\frak I}(F)+T_{2\frak I}\big)/F\iota_{2\frak I},$$
as described in \eqref{newdefleala0}. We have the induced
automorphism on $\mathcal A_{2\frak I}$, which is also denoted by
$\sigma$ for simplicity. Thus, we obtain the fixed algebra
$$\mathcal A_{2\frak I}^\sigma
=\big(\text{sp}_{2\frak I}(F)+T^+\big),$$
where we again omit the bars. Let
$$\hat{\mathcal A}_{2\frak I}:=\mathcal A_{2\frak I}\otimes F[t^{\pm 1}]\oplus Fc \oplus Fd^{(0)},$$
as in \eqref{maxreali}. We extend $\sigma$ to $\hat{\mathcal
A}_{2\frak I}$ as
\begin{equation*}
\hat\sigma(x\otimes t^k):=(-1)^k\sigma(x)\otimes t^k,
\label{extendedauto}
\end{equation*}
and identically on $Fc\oplus Fd^{(0)}$. Then, we obtain the fixed
algebra
\begin{equation}\label{maxtwisted}
\hat{\mathcal A}_{2\frak I}^{\hat\sigma}= \big((\text{sp}_{2\frak
I}(F)+T^{+})\otimes F[t^{\pm 2}]\big) \oplus\big((\frak
s+T^-)\otimes tF[t^{\pm 2}]\big) \oplus Fc\oplus Fd^{(0)}.
\end{equation}
Note that $\hat{\mathcal A}_{2\frak I}^{\hat\sigma}$ contains the
subalgebra
$$\mathcal L^{ms}:=
\big(\text{sp}_{2\frak I}(F)\otimes F[t^{\pm 2}]\big)
\oplus\big(\frak s\otimes tF[t^{\pm 2}]\big) \oplus Fc\oplus
Fd^{(0)},
$$
which is called a {\bf minimal standard twisted LALA} of type
$\text {C}_{\frak I}^{(2)}$. As in the case of type $\text
{A}_{\frak I}^{(1)}$, we have
\begin{theorem}
$\mathcal{L}^{max} := \hat{\mathcal A}_{2\frak I}^{\hat\sigma}$ is
a {\bf maximal twisted LALA} of type $\text {C}_{\frak I}^{(2)}$,
i.e., any LALA of type $\text {C}_{\frak I}^{(2)}$ is a subalgebra
of $\hat{\mathcal A}_{2\frak I}^{\hat\sigma}$. Moreover, any LALA
of type $\text {C}_{\frak I}^{(2)}$ contains a {\bf minimal
twisted LALA}
\begin{equation*}
\mathcal L(p): =\big(\text{sp}_{2\frak I}(F)\otimes F[t^{\pm
2}]\big) \oplus\big(\frak s\otimes tF[t^{\pm 2}]\big) \oplus
Fc\oplus F(p+d^{(0)}) \label{minreali2}
\end{equation*}
for some $p\in T^+$.
\end{theorem}

We must emphasize that the usual twisting process works for the
locally loop algebra $\sll_{2\frak I}(F)\otimes F[t^{\pm 1}]$ but
also for the bigger algebra $\big(\sll_{2\frak I}(F)+T_{2\frak
I}\big)\otimes F[t^{\pm 1}]$. Note that a LALA $\mathcal L$ of
type $\text{C}_\frak I^{(2)}$ has a decomposition $\mathcal{L} =
\mathcal{L}_c \oplus D$ for a homogeneous complement $D$ that
satisfies $D(p) = F(p + d^{(0)}) \subset D \subset D^{max}$ and
$\mathcal{L}_c \oplus D^{max} =
\hat{\mathcal{A}}_{2\frak{I}}^{\hat{\sigma}}$.

\medskip

Next, we explain how the classification of LALAs is conducted.
First,  we classify the cores of the LALAs. We show that the core
of a LALA is a locally Lie $1$-torus and that the core is a
universal covering of a locally loop algebra. We also show that
there is a one to one correspondence between reduced root systems
extended by $\Bbb Z$ (which are classified in \cite[Cor.15]{Y3},
as the class of reduced locally affine root systems) and the cores
of LALAs.

\smallskip

The second step of the classification process involves determining
a complement $D$ of the core $\mathcal{L}_c$ of a LALA
$\mathcal{L} = \mathcal{L}_c \oplus D$. As explained above, we can
obtain $D \subset \der_F\mathcal L_c \subset \der_F L$ or $D
\subset \oder_F \mathcal{L}_c \subset \oder_F L$, where $
L=\mathcal L_c/Z(\mathcal L_c)$ is the centerless core (which is a
locally loop algebra).

\smallskip

Now,  we need some information about $\der_F L$. Derivations of
this type of algebra were studied in \cite{BM}, \cite{B}, and
\cite{NY}. However, the derivations of a locally loop algebra are
new. We can use some results from 
\cite{A1}
for the untwisted case
since $L$ is a tensor product algebra (see Remark
\ref{azamresult}). 
However, we need to determine the twisted case.
Thus, we propose a new method. Clearly, we need to use the
classification of $\der_F \frak g$ for a locally finite split
simple Lie algebra $\frak g$, as described by Neeb in \cite{N1}.
Fortunately, we do not need all of the information about $\der_F
L$ to classify $D$. In fact, we only need to know the diagonal
derivations of degree $m$. To explain this, we note that $L$ has
double grading, i.e.,
$$ L
=\bigoplus_{\alpha\in\Delta\cup\{0\}}\ \bigoplus_{k\in\Bbb Z}\
L_\alpha^k,$$ where $\Delta$ is a locally finite irreducible root
system. The diagonal derivations of degree $m$ denote the space
$$(\der_F L)^m_0
:=\{d\in\der_F L\mid d(L_\alpha^k)\subset L_\alpha^{k+m} \
\text{for all $\alpha\in\Delta$ and $k\in\Bbb Z$}\}.$$ It is
crucial to determine the case where $m=0$, i.e., $(\der_F L)^0_0$.
Next, $(\der_F L)^m_0$ can be determined easily for the untwisted
case. However, for the twisted case, $(\der_F L)^m_0$ is still
difficult when $m$ is odd. Finally, using some new techniques (see
Lemma \ref{oddcomshift} and Lemma \ref{extder}), the
classification of $(\der_F L)^m_0$ is completed in Theorem
\ref{odddegreeth}.

\smallskip

If we take $D$ as a homogeneous complement of the graded algebra
$\mathcal L_c$, then $D$ has $\Bbb Z$-grading, e.g.,
$$D=\bigoplus_{m\in\Bbb Z}\ D^m,$$
and each $D^m$ can be identified with a subspace of the known
space $(\der_FL)_0^m$. Finally, we classify the Lie brackets on
$D$ and the concrete brackets are described in Example
\ref{exampleLALA}.

\smallskip

The remainder of this paper is organized as follows. In Section 2,
we define a locally Lie $G$-torus and we consider a locally Lie
$1$-torus as a special case. In Section 3, we introduce a locally
loop algebra, which is a centerless locally Lie $1$-torus. We
classify locally Lie $1$-tori in general. We prove that a
centerless locally Lie $1$-torus is uniquely determined by a root
system extended by $\Bbb Z$, and that a locally Lie $1$-torus is a
locally loop algebra or a universal covering of a locally loop
algebra. In Section 4, we recall the definition of a LEALA and we
prove some general properties of a LEALA. In Section 5, we
summarize and prove several properties related to LEALAs of
nullity $0$. In Section 6, we define a LALA. We show that the core
of a LALA is a universal covering of a locally loop algebra and we
then construct many examples of LALAs. In Sections 7 and 8, we
classify untwisted LALAs and twisted LALAs. Finally, we provide
our main theorem.
\begin{theorem}
The examples in Example \ref{exampleLALA} comprise all LALAs.
\end{theorem}
In Section 9, we discuss standard LALAs.

\enskip

The authors thank Karl-Hermann Neeb and Erhard Neher for helpful
discussions and suggestions regarding this study.

\enskip

\section{Locally Lie $G$-tori}
\label{sec:locallylieg}

\noindent To classify $\mathcal{L}_c$ and $L =
\mathcal{L}_c/Z(\mathcal{L}_c)$, we need to study localy Lie
$G$-tori, which are very useful. Let $\Delta$ be a locally finite
irreducible root system (see \cite{LN1}), and we denote the Cartan
integer
$$\frac{2(\mu,\nu)}{(\nu,\nu)}$$ by
$\langle\mu,\nu\rangle$ for $\mu,\nu\in\Delta$, while we also let
$\langle 0,\nu\rangle:=0$ for all $\nu\in\Delta$. Recall that
$\Delta$ is called {\bf reduced} if $2\alpha\notin\Delta$ for all
$\alpha\in\Delta$. We define the subset
$$\Delta^{\red}:=
\{\alpha\in\Delta\mid\ \frac{1}{2}\alpha\notin\Delta\}
$$
 of $\Delta$,
 which is a reduced locally finite irreducible root system.
 Note that $\Delta=\Delta^{\red}$ if $\Delta$ is reduced.
To simplify the description later, we partition the locally finite
irreducible root system $\Delta$ according to length. The roots of
$\Delta$ of minimal length are called {\bf short}.  The roots of
$\Delta$, which are two times a short root of $\Delta$, are called
{\bf extra long}. Finally, the roots of $\Delta$, which are
neither short nor extra long,  are called {\bf long}. We denote the
subsets of the short, long,  and extra long roots of $\Delta$ by
$\Delta_{\sh}$, $\Delta_{\lng}$, and $\Delta_{\ex}$, respectively.
Thus,
$$
\Delta=\Delta_{\sh}\sqcup\Delta_{\lng}\sqcup\Delta_{\ex}.
$$
Clearly, the last two terms in this union may be empty. Indeed,
$$\Delta_{\lng}=\emptyset \quad \Longleftrightarrow\quad
\text{$\Delta$ is a simply laced type or type $\text{BC}_1$,}$$
and
$$\Delta_{\ex}=\emptyset \quad \Longleftrightarrow\quad
\Delta=\Delta^{\red}.
$$

Let $G= (G,+,0)$ be an arbitrary abelian group. In general, for a
subset $S$ of $G$, the subgroup generated by $S$ is denoted by
$\langle S\rangle$.

\begin{definition}
A Lie algebra ${\mathcal L}$ is called a {\bf locally Lie
$G$-torus of type $\Delta$} if:
\begin{itemize}
\item[(LT1)] ${\mathcal L}$ has a decomposition into subspaces
$${\mathcal L} =  \bigoplus_{\mu \in \Delta \cup \{0\}, \ { g \in G}} \  {\mathcal L}^g_\mu$$
such that $[{\mathcal L}^g_\mu, {\mathcal L}^h_\nu]  \subset
{\mathcal L}^{g+h}_{\mu+\nu}$ for $\mu,\nu, \mu+\nu \in \Delta
\cup \{0\}$ and $g,h\in G$; \item[(LT2)]
 For every $g \in G$, \ ${\mathcal L}_0^g =  \sum_{\mu \in \Delta,\
h \in G} \  [{\mathcal L}_\mu^h, {\mathcal L}_{-\mu}^{g-h}]$;
\item[(LT3)] For each nonzero $x \in {\mathcal L}_\mu^g$ ($\mu \in
\Delta, g \in G)$, an element $y \in {\mathcal L}_{-\mu}^{-g}$
exists such that $t : = [x, y] \in {\mathcal L}_0^0$ satisfies
$[t, z] = \langle \nu,\mu\rangle z$ for all $z \in {\mathcal
L}_\nu^h$ $(\nu \in \Delta \cup \{0\}, h \in G)$; \item[(LT4)]
$\dim {\mathcal L}_\mu^g \leq 1$ for $\mu \in \Delta$ and $g \in
G$, and $\dim {\mathcal L}_\mu^0 = 1$ if $\mu \in \Delta^{\red}$;
\item[(LT5)] $\langle\supp{\mathcal L}\rangle=G$, where
$\supp{\mathcal L} = \{g \in G \mid {\mathcal L}_\mu^g \neq 0 \
\text{for some} \ \mu \in \Delta \cup \{0\}\}$.
\end{itemize}

If $\Delta$ is finite, we omit the term `locally' and simply call
it a {\bf Lie $G$-torus}. Furthermore, if $G\cong\Bbb Z^n$, then
${\mathcal L}$ is called a {\bf locally Lie $n$-torus}, or simply
a {\bf locally Lie torus}. We refer to the rank of $\Delta$ as the
{\bf rank} of ${\mathcal L}$.
\end{definition}

\remark (i)  Condition (LT5) is simply for convenience but if it
fails to hold, we may replace $G$ by the subgroup generated by $
\supp\mathcal L $.

(ii)  It follows from (LT1) and (LT3) that $\mathcal L $ admits a
grading by the root lattice $\langle\Delta\rangle$.

\noindent Let
\begin{equation}\label{rootspaces}
\mathcal L _\lambda := \bigoplus_{g \in G}\  \mathcal L _\lambda^g
\end{equation} for
$\lambda \in \langle\Delta\rangle$, where $\mathcal L _\lambda^g =
0$ if $\lambda \not \in \Delta \cup \{0\}$. Then, $\mathcal L  =
\oplus_{\lambda \in \langle\Delta\rangle}\ \mathcal L _\lambda$
and $[\mathcal L _\lambda, \mathcal L _\mu] \subset \mathcal L
_{\lambda+\mu}$.

(iii) $\mathcal L $ is also graded by the group $G$, i.e., if
\begin{equation}
\mathcal L ^g := \bigoplus_{\mu \in \Delta \cup \{0\}}\ \mathcal L
_\mu^g,
\end{equation} then
$\mathcal L  = \oplus_{g \in G}\ \mathcal L ^g$ and $[\mathcal L
^g,\mathcal L ^h] \subset \mathcal L ^{g+h}$. In addition,
$\supp\mathcal L = \{g \in G \mid \mathcal L^g \neq 0\}$.

(iv)
 From (LT3), for $\mu \in \Delta^{\red}$, we can see that the
elements $e_\mu \in \mathcal L _\mu^0$, $f_\mu \in \mathcal L _{-\mu}^0$,
 and $\mu^\vee=\mu^\vee_\mathcal L: = [e_\mu,f_\mu]$ exist such that $[\mu^\vee, z] = \langle \nu,\mu \rangle z$
 for all $z \in \mathcal L _\nu^h$,
$\nu \in \Delta$ and $h \in G$. Thus, the elements
$e_\mu,f_\mu,\mu^\vee$ determine a canonical basis for a copy of
the Lie algebra $\sll_2(F)$. (Note that $\mu^\vee$ is a unique
element in $ [\mathcal L_\mu^0,\mathcal L_{-\mu}^0]$ that
satisfies the property.)
The subalgebra $\mathfrak g$ of $\mathcal L $ generated by the
subspaces $\mathcal L _\mu^0$ for $\mu \in \Delta^{\red}$ is a
locally finite split simple Lie algebra with the split Cartan
subalgebra
$$\mathfrak h: = \sum_{\mu \in \Delta^{\red}}\
[\mathcal L _\mu^0, \mathcal L _{-\mu}^0],$$
and $\mu^\vee$ are
the coroots in $\mathfrak h$. (We can show this in the same manner
as the proof of \cite[Prop.8.3]{MY}, or see \cite[Sec.III]{St}).
Note that if $\Delta$ is finite, then $\mathfrak g$ is a
finite-dimensional split simple Lie algebra. Furthermore,
$\Delta^{\red}$ may be replaced by $\Delta$ in the definition of
$\mathfrak g$ and $\mathfrak h$ since it can be shown in the same
manner described by \cite[Thm.5.1]{Y1} that $\mathcal L _{2\nu}^0
= 0$ for all $\nu \in   \Delta^{\red}$. We say that the pair
$(\mathfrak g,\mathfrak h)=(\mathfrak g,\mathfrak h)_\mathcal L$
is the {\bf grading pair} of $\mathcal L$.

(v) A locally Lie $G$-torus is perfect, and thus it has a
universal covering.

(vi) Let $\mathcal L$ be a locally Lie $G$-torus and $Z$ is its
center. Then, we can see that $Z\subset  \mathcal L_0$. In
addition, $\mathcal L/Z$ is a locally Lie $G$-torus with a trivial
center. In general, a Lie algebra with a trivial center is called
{\bf centerless}.

 \endremark

We define the root systems of locally Lie $G$-tori. Let $\mathcal
L=\oplus_{\mu\in\Delta\cup\{0\}}\ \oplus_{g\in G}\ \mathcal
L_{\mu}^g$ be a locally Lie $G$-torus. For each $\mu\in\Delta$,
let
$$
S_\mu:=\{g\in G\mid \mathcal L_\mu^g\neq 0\},
$$
and we refer to
$$
\tilde\Delta:=\{S_\mu\}_{\mu\in\Delta}
$$
as the {\bf root system} of $\mathcal L$ (which is called an
extension datum in \cite{LN2}). This system fits into the system
introduced in \cite{Y1}. Let us state the precise definition. A
family of subsets $S_\mu$ of $G$ indexed by $\Delta$, such as
$\{S_\mu\}_{\mu\in\Delta}$, is called a {\bf root system extended
by $G$} if
\begin{align}
&\langle\cup_{\mu\in\Delta}\ S_{\mu}\rangle=G,
\label{S0} \\
&S_\nu-\langle\nu,\mu\rangle S_\mu\subset
S_{\nu-\langle\nu,\mu\rangle\mu} \quad\text{for all
$\mu,\nu\in\Delta$, and}
\label{S1} \\
&0\in S_{\mu} \quad\text{ for all $\mu\in\Delta^{\red}$.}
\label{S2}
\end{align}
Moreover, $\{S_\mu\}_{\mu\in\Delta}$ is called {\bf reduced} if
\begin{equation}
S_{2\mu}\cap 2S_{\mu}=\emptyset \quad\text{for all
$2\mu,\mu\in\Delta$.} \label{S3}
\end{equation}
In the same manner described in \cite[Thm 5.1]{Y1}, we can show
that the root system $\tilde\Delta$ of $\mathcal L$ is a reduced
root system extended by $G$, i.e., $\tilde\Delta$ satisfies
\eqref{S0}, \eqref{S1}, \eqref{S2}, and \eqref{S3}. Moreover,
$$S_\mu=S_\nu
\quad \text{if $\mu$ and $\nu$ are the same length, and}
$$
\begin{equation}\label{shortbig}
S_\nu\subset S_\mu \quad\text{ for all $\nu\in\Delta$ if $\mu$ is
a short root}.
\end{equation}
Finally, if we let
\begin{equation}
S_0:=\{g\in G\mid\mathcal L_0^g\neq 0\},
\end{equation}
then we obtain
\begin{equation}
S_0=S_\mu+S_\mu \label{S4}
\end{equation}
for a short root $\mu$.

\begin{lemma}\label{locallystatements}
A locally Lie $G$-torus $\mathcal L$ of type $\Delta$ is a
directed union of Lie $G$-tori. In particular, $\mathcal
L=\bigcup_{\Delta'} \ \mathcal L_{\Delta'}$, where $\Delta'$ is a
finite irreducible full subsystem of $\Delta$ that contains a
short root and $\mathcal L_{\Delta'}$ is the subalgebra of
$\mathcal L$ generated by $\mathcal L_\alpha$ for all
$\alpha\in\Delta'$.

Furthermore, if $G$ is torsion-free, then a locally Lie $G$-torus
$\mathcal L$ of type $\Delta$ is a directed union of Lie $n$-tori,
where $n$ runs over a certain subset of $\Bbb N$. In particular,
$\mathcal L=\bigcup_{\Delta', G'} \ \mathcal L_{\Delta'}^{G'}$,
where $G'$ is a finitely generated subgroup of $G$ and
 $\mathcal L_{\Delta'}^{G'}$ is the subalgebra of $\mathcal L$ generated by
$\mathcal L_\alpha^g$ for all $\alpha\in\Delta'$ and $g\in G'$.

\end{lemma}

\proof Since $S=S_\mu$ generates $G$ for a short root $\mu$ by
\eqref{shortbig}, then it is easy to check that $\mathcal
L_{\Delta'}$ is a Lie $G$-torus. Hence, the statement is true
since $\Delta$ is a directed union of finite irreducible full
subsystems that contain a short root (see \cite[3.15 (b) and the
proof]{LN2}). The second statement follows from the fact that $G$
is a directed union of finitely generated subgroups, and the fact
that a finitely generated torsion-free abelian group is free. \qed

\remark\label{subtori} Let $\Delta$ and $G$ be as given in Lemma
\ref{locallystatements}. For a locally finite irreducible full
subsystem $\Delta'$ of $\Delta$, and for a subgroup $G'$ of $G$,
we put $\mathcal{M} = \mathcal{L}_{\Delta'}^{G'}$, which can be
defined as given in Lemma \ref{locallystatements}. Then,
$$\mathcal{M} = \bigoplus_{\mu' \in \Delta' \cup \{ 0 \},\ g' \in G'} \mathcal{M}_{\mu'}^{g'},$$
where
$$\mathcal{M}_{\mu'}^{g'} = \mathcal{M} \cap \mathcal{L}_{\mu'}^{g'}\quad (\mu' \in \Delta' \cup \{ 0 \},\ g' \in G').$$
In fact, we obtain
$$\mathcal{M}_{\mu'}^{g'} = \mathcal{L}_{\mu'}^{g'}\quad (\mu' \in \Delta',\ g' \in G'),$$
and since
 $\mathcal M$ is generated by
$\mathcal L_\mu^g$ for all $\mu\in\Delta'$ and $g\in G'$, then for
$g' \in G'$, we have
$$\mathcal{M}_0^{g'} = \sum_{\mu' \in \Delta',\ h' \in G'} [ \mathcal{L}_{\mu'}^{h'} , \mathcal{L}_{-\mu'}^{g'-h'} ]
= \sum_{\mu' \in \Delta',\ h' \in G'} [ \mathcal{M}_{\mu'}^{h'} ,
\mathcal{M}_{-\mu'}^{g'-h'} ].
$$
Thus, we can check conditions (LT1) -- (LT5) for $\mathcal{M}$,
which implies that $\mathcal{M}$ is a locally Lie $G'$-torus.

\endremark

\enskip

\section{Locally loop algebras}
\label{sec:locallyloop}

\noindent For any index set ${\mathfrak I}$, in the introduction,
we defined
$$M_\frak I(F) = \left\{ \left. (a_{ij})_{i,j \in \frak I} \right| a_{ij} \in F \right\}
\approx \text{Map}(\frak I \times \frak I,F),$$
as the set of all
matrices of size $\frak I$, which is naturally a vector space over
$F$. Let $\gl_{\mathfrak I}(F)$ be the subspace of $M_\frak I(F)$
that comprises matrices with only a finite number of nonzero
entries. Then, $\gl_{\mathfrak I}(F)$ is an associative algebra
and a Lie algebra with the usual commutator bracket. 
Furthermore,
we can define the trace of a matrix in $\gl_{\mathfrak I}(F)$, and
the subalgebra of $\gl_{\mathfrak I}(F)$ that comprises trace $0$
matrices is denoted by $\sll_{\mathfrak I}(F)$, as follows.
$$
\sll_{\mathfrak I}(F)=\{x\in \gl_{\mathfrak I}(F)\mid \tr (x)=0\}
$$
We note that $M_\frak I(F)$ is not an algebra if $\frak I$ is
infinite, but
\begin{equation*}\label{}
M_\frak I^{\fin}(F):= \{x\in M_\frak I(F)\mid \text{{\small each
row and column of} $x$ {\small have only finitely many nonzero
entries}}\}
\end{equation*}
is an associative algebra with the identity matrix
$\iota=\iota_{\mathfrak I}$, and a Lie algebra with the commutator
bracket. In fact, this gives the Lie algebra of derivations of
$\sll_\frak I(F)$, as described by Neeb \cite{N1}. In particular,
we have
$$
[M_\frak I^{\fin}(F),\sll_\frak I(F)]\subset\sll_\frak I(F)\quad
\mbox{and}\quad \der_F(\sll_\frak I(F)) \simeq \ad(M_\frak
I^{\fin}(F)).
$$

\medskip

As a result, we note that there are 14 types of locally loop
algebras, i.e., we obtain:
$${\rm A}_\mathfrak{I}^{(1)},\ {\rm B}_\mathfrak{I}^{(1)},\ {\rm C}_\mathfrak{I}^{(1)},\ {\rm D}_\mathfrak{I}^{(1)},\
{\rm B}_\mathfrak{I}^{(2)},\ {\rm C}_\mathfrak{I}^{(2)},\ {\rm
BC}_\mathfrak{I}^{(2)},\ {\rm E}_6^{(1)},\ {\rm E}_7^{(1)},\ {\rm
E}_8^{(1)},\ {\rm F}_4^{(1)},\ {\rm G}_2^{(1)},\ {\rm F}_4^{(2)},\
{\rm G}_2^{(3)},
$$
where we mainly assume that $\mathfrak{I}$ is
infinite since we already know the affine Lie algebras.

\medskip

The locally finite split simple Lie algebra of type $\text
X_{\mathfrak I}$ is defined as a subalgebra of $\sll_\frak I(F)$,
$\sll_{2\frak I+1}(F)$ or $\sll_{2\frak I}(F)$ as follows:

Type $\text A_{\mathfrak I}$: $\sll_{\mathfrak I}(F)$;

Type $\text B_{\mathfrak I}$: $\text o_{2{\mathfrak I}+1}(F)
=\{x\in \sll_{2{\mathfrak I}+1}(F)\mid sx=-x^Ts\}$;

Type $\text C_{\mathfrak I}$: $\text {sp}_{2{\mathfrak I}}(F)
=\{x\in \sll_{2{\mathfrak I}}(F)\mid sx=-x^Ts\}$;

Type $\text D_{\mathfrak I}$: $\text o_{2{\mathfrak I}}(F) =\{x\in
\sll_{2{\mathfrak I}}(F)\mid sx=-x^Ts\}$,

\noindent where $\frak I$ is assumed to be infinite,
 $x^T$ is the transpose of $x$,
 and
\begin{equation}\label{definings}
s=
\begin{pmatrix}
   0   &\iota   & 0   \\
  \iota   & 0    & 0  \\
   0   & 0   &   1
\end{pmatrix}
\ \text{for $\text B_{\mathfrak I}$},\ s=
\begin{pmatrix}
   0   &-\iota       \\
  \iota   & 0   \\
\end{pmatrix}
\ \text{for $\text C_{\mathfrak I}$},\ \ \text{or}\ \ s=
\begin{pmatrix}
   0   &\iota       \\
  \iota   & 0   \\
\end{pmatrix}
\ \text{for $\text D_{\mathfrak I}$}.
\end{equation}
Note that $s\in M_{2\frak I + 1}^{\fin}(F)$ for ${\rm B}_\frak I$
and $s \in M_{2\frak I}^{\fin}(F)$ for ${\rm C}_\frak I$ or ${\rm
D}_\frak I$, and that $s^2=\iota_{2\frak I+1}$ for $\text
B_{\mathfrak I}$, $s^2=-\iota_{2\frak I}$ for $\text C_{\mathfrak
I}$ and $s^2=\iota_{2\frak I}$ for $\text D_{\mathfrak I}$. In
addition, $\text B_{\mathfrak I}$, $\text C_{\mathfrak I}$, or
$\text D_{\mathfrak I}$ is the fixed algebra of $\sll_{2\frak
I+1}(F)$ or $\sll_{2\frak I}(F)$ by an automorphism $\sigma$,
which are defined as
\begin{equation}\label{originalinv}
\text{ $\sigma(x)=-sx^Ts\ $\ for $\text B_{\mathfrak I}$ or $\text
D_{\mathfrak I}$, and $\sigma(x)=sx^Ts\ $\ for $\text C_{\mathfrak
I}$.}
\end{equation}

In
\cite{NS}, Neeb and Stumme showed that these algebras comprise
all of the infinite-dimensional locally finite split simple Lie
algebras. In addition, they are considered to be locally Lie
$0$-tori (in the case where $G=\{0\}$). Moreover, since locally
finite split simple Lie algebras are centrally closed (see
\cite{NS}), we have the equality \{infinite-dimensional locally
Lie $0$-tori\} =\{infinite-dimensional locally finite split simple
Lie algebras\}. We note that Lie $0$-tori are exact
finite-dimensional split simple Lie algebras. In the present
study, we are interested in the class of locally Lie $1$-tori.

\medskip

Let $F[t^{\pm 1}]$ be the algebra of Laurent polynomials over $F$.
We call one of the following four Lie algebras an {\bf untwisted
locally loop algebra}:

(1) Type $\text A_{\mathfrak I}^{(1)}$: $\sll_{\mathfrak
I}(F)\otimes F[t^{\pm 1}]$;

(2) Type $\text B_{\mathfrak I}^{(1)}$: $\text o_{2{\mathfrak
I}+1}(F)\otimes F[t^{\pm 1}]$;

(3) Type $\text C_{\mathfrak I}^{(1)}$: $\text {sp}_{2{\mathfrak
I}}(F)\otimes F[t^{\pm 1}]$;

(4) Type $\text D_{\mathfrak I}^{(1)}$: $\text o_{2{\mathfrak
I}}(F)\otimes F[t^{\pm 1}]$.

\noindent (In addition, it is called an untwisted loop algebras if
$\mathfrak I$ is finite.) Each of the following three Lie algebras
is called a {\bf twisted locally loop algebra}:

(5) Type $\text B_{\mathfrak I}^{(2)}$: $(\text o_{2{\mathfrak
I}+1}(F)\otimes F[t^{\pm 2}]\oplus (\frak s\otimes  tF[t^{\pm
2}])$,

\noindent where $\frak s=F^{(2{\mathfrak I}+1)}$ is the natural
$\text o_{2{\mathfrak I}+1}(F)$-module;

(6) Type $\text C_{\mathfrak I}^{(2)}$: $(\text {sp}_{2{\mathfrak
I}}(F)\otimes F[t^{\pm 2}]) \oplus (\frak s\otimes  tF[t^{\pm
2}])$,

\noindent where $\frak s=\{x\in \sll_{2{\mathfrak I}}(F)\mid
sx=x^Ts\}$;

(7) Type $\text {BC}_{\mathfrak I}^{(2)}$: $(\text o_{2{\mathfrak
I}+1}(F)\otimes F[t^{\pm 2}]) \oplus (\frak s\otimes  tF[t^{\pm
2}])$,

\noindent where $\frak s=\{x\in \sll_{2{\mathfrak I}+1}(F)\mid
sx=x^Ts\}$. (In addition, it is called a twisted loop algebra if
$\mathfrak I$ is finite.) Note that $\sll_{2{\mathfrak
I}}(F)=\text {sp}_{2{\mathfrak I}}(F)\oplus \frak s$ for ${\rm
C}_\frak I^{(2)}$ and $\sll_{2{\mathfrak I}+1}(F)=\text
o_{2{\mathfrak I}+1}(F)\oplus \frak s$ for ${\rm BC}_\frak
I^{(2)}$.

\medskip

The Lie bracket of each untwisted type is natural, i.e.,
$[x\otimes t^m, y\otimes t^n]=[x,y]\otimes t^{m+n}$. The Lie
bracket of type $\text C_{\mathfrak I}^{(2)}$ or $\text
{BC}_{\mathfrak I}^{(2)}$ is also natural, and we have
\begin{align*}
&[\text {sp}_{2{\mathfrak I}}(F), \frak s]\subset  \frak s \quad
\text{and} \quad
[\frak s,\frak s]\subset \text {sp}_{2{\mathfrak I}}(F)\quad \text{for}\ {\rm C}_\frak I^{(2)},\\
&[\text o_{2{\mathfrak I}+1}(F), \frak s]\subset  \frak s
\quad\text{and}\quad [\frak s,\frak s]\subset \text o_{2{\mathfrak
I}+1}(F)\quad \text{for}\ {\rm BC}_\frak I^{(2)}.
\end{align*}
Note that $\text C_{\mathfrak I}^{(2)}$ or $\text {BC}_{\mathfrak
I}^{(2)}$ is the fixed subalgebra of $\sll_{2{\mathfrak
I}}(F)\otimes F[t^{\pm 1}]$ or $\sll_{2{\mathfrak I}+1}(F)\otimes
F[t^{\pm 1}]$ by the automorphism $\hat\sigma$, which is defined
as
\begin{equation}\label{tildeinv}
\hat\sigma(x\otimes t^m):=(-1)^m\sigma(x)\otimes t^m
\end{equation}
(see \eqref{originalinv}). This construction is called a {\bf
twisting construction} by an automorphism $\sigma$.

\smallskip

For $\text B_{\mathfrak I}^{(2)}$, we have $\text o_{2{\mathfrak
I}+1}(F)\frak s\subset \frak s$, and thus we define the bracket of
$\text o_{2{\mathfrak I}+1}(F)$ and $\frak s$ by the natural
action, i.e., $[x,v]=xv=-[v,x]$ for $x\in \text o_{2{\mathfrak
I}+1}(F)$ and $v\in \frak s$. We define a bracket on $\frak s$
such that $[\frak s,\frak s]\subset \text o_{2{\mathfrak I}+1}(F)$
as follows. First, let $(\cdot,\cdot)$ be the bilinear form on
$\frak s$ determined by $s$. Then, there is a natural
identification
$$
\text o_{2{\mathfrak I}+1}(F)=D_{\frak s,\frak
s}:=\spa_F\{D_{v,v'}\mid v,v'\in \frak s\},
$$
where $D_{v,v'}\in\End (\frak s)$ is defined by
$D_{v,v'}(v'')=(v',v'')v-(v,v'')v'$ for $v''\in \frak s$. Thus, we
define $[v,v']:=D_{v,v'}$. Note that $[v',v]=-[v,v']$. It is easy
to check that the bracket
\begin{align*}
&[x\otimes t^{2m}+v\otimes  t^{2m'+1},x'\otimes t^{2n}+v'\otimes  t^{2n'+1}]\\
=&[x,x']\otimes t^{2(m+n)}+D_{v,v'}\otimes t^{2(m'+n'+1)}
+xv'\otimes t^{2(m+n')+1}-x'v\otimes t^{2(m'+n)+1}
\end{align*}
defines a Lie bracket for $m,m',n,n'\in\Bbb Z$.

\smallskip

There is a twisting construction for $\text B_{\mathfrak I}^{(2)}$
(see \cite{N2}), which we discuss in Section 7, but we also
consider that the simple description of $\text B_{\mathfrak
I}^{(2)}$ is important for developing the theory of locally Lie
tori.

\remark We often omit the term `untwisted' or `twisted' and we
simply refer to a locally loop algebra. In addition, a locally loop
algebra can be simply called a loop algebra in more general
theory. For example, $A\otimes F[t^{\pm 1}]$ for any algebra $A$
is called a loop algebra of $A$. However, we use the term
`locally' in this study to distinguish the familiar loop algebras
in Kac-Moody theory.
\endremark

\medskip

We can easily check that
$$
\text{\bf all locally loop algebras are centerless locally Lie
1-tori}.
$$
For example, let $\Delta$ be the root system of type
$\text{BC}_{\mathfrak I}$, and we put $\mathfrak g=\text
o_{2\mathfrak I+1}(F)$ and $\frak s\subset \sll_{2{\mathfrak
I}+1}(F)$, as defined above. Let $\mathfrak h$ be the Cartan
subalgebra of $\mathfrak g$ that comprises diagonal matrices.
Then, $\mathfrak h$ decomposes $\mathfrak g$ into the root spaces,
such as $\mathfrak g=\mathfrak
h\oplus\bigoplus_{\mu\in\Delta^{\red}}\ \mathfrak g_\mu$, and
$\frak s$ into the weight spaces, such as $\frak
s=\bigoplus_{\mu\in\Delta \cup \{ 0 \}}\ \frak s_\mu$, where
$\Delta^{\rm red}$ is of type $\text B_\frak I$. Therefore, the
twisted locally loop algebra $\mathcal L:= (\mathfrak g\otimes
F[t^{\pm 2}])\oplus (\frak s\otimes  tF[t^{\pm 2}])$ of type
$\text{BC}_{\mathfrak I}^{(2)}$ is decomposed into
$$\bigoplus_{m\in \Bbb Z}\ \bigg((\mathfrak h\otimes Ft^{2m})
\oplus\bigoplus_{\mu\in\Delta^{\red}}\ (\mathfrak g_\mu\otimes
Ft^{2m}) \oplus \bigoplus_{\mu\in\Delta\cup\{0\}}\ (\frak
s_\mu\otimes Ft^{2m+1})\bigg).$$ This gives a natural double
grading by the groups $\langle\Delta\rangle$ and $\Bbb Z$, and we
can check the axioms of a locally Lie torus. In addition, the
center is contained in $\mathcal L_0=\mathfrak h\otimes F[t^{\pm
2}]$, and thus $\mathcal L$ is a centerless locally Lie 1-torus.
The grading subalgebra is equal to $\frak g=\text o_{2\mathfrak
I+1}(F)$. We refer to the $\frak g$-module $\frak s$ as the {\bf
grading module}.

\medskip

 The following lemma was proved for the base field $\Bbb C$ in [ABGP],
 but it also works for our base field $F$.
 We use the notation
$$
\tilde\Delta:=\{S_\mu\}_{\mu\in\Delta}
$$
(defined in Section 2) for the case where
$\langle\cup_{\mu\in\Delta}\ S_{\mu}\rangle=\Bbb Z$ (the root
system $\Delta$ extended by $\Bbb Z$).

\lemma Let $\Delta$ be a {\bf finite} irreducible root system. Let
$\mathcal L =  \oplus_{\mu \in \Delta \cup \{0\},\  { m  \in \Bbb
Z}} \  \mathcal L^m _\mu$ and $\mathcal M =  \oplus_{\mu \in
\Delta \cup \{0\},\  { m  \in \Bbb Z}} \  \mathcal M^m _\mu$ be
centerless Lie 1-tori, which have the same root system
$\tilde\Delta$ extended by $\Bbb Z$. Then, an isomorphism
$\vvp:\mathcal L\longrightarrow\mathcal M$ exists such that
\begin{equation}\label{dagger}
\vvp(\mu^\vee_\mathcal L)=\mu^\vee_\mathcal M \quad\text{and}\quad
\vvp(\mathcal L_\mu^m )=\mathcal M_\mu^m \quad\text{ for all
$\mu\in\Delta$ and $m \in \Bbb Z$.}
\end{equation}
\endlemma
\remark If $\mathcal L$ is a loop algebra, then $\tilde\Delta$
determines $\mathcal L$, i.e., there is a one to one
correspondence between loop algebras and root systems extended by
$\Bbb Z$ (see \cite{Y3}). In particular, $\tilde\Delta$ determines
whether the loop algebra is untwisted or twisted.
\endremark
\proof Let $0\neq e_{\pm\mu}\in\mathcal L_{\pm\mu}^0$ and
$\mu^\vee$ be an $\sll_2$-triple for $\mu\in\Pi$, where $\Pi$ is a
root base of $\Delta$ and let $0\neq x_{\pm\nu}\in\mathcal
L_{\pm\nu}^{\mp 1}$ and $\nu^\vee$ be an $\sll_2$-triple, where
$\nu\in\Delta$ is the highest long (or short) root relative to
$\Pi$ (depending on the type $\tilde\Delta$). Then, the set
$$\{e_{\pm\mu},\ \mu^\vee,\ x_{\pm\nu},\ \nu^\vee
\mid \mu\in\Pi\}$$ satisfies the Serre relations. Hence, using the
Gabber-Kac Theorem (e.g., see \cite[Thm 4, p.381]{MP}),
a homomorphism $\psi$ exists
 from the derived affine Lie algebra $A$
 (which is a 1-dimensional central extension of a loop algebra), which is
determined by
 $\Delta$ and $\nu$
 (or $\tilde\Delta$) into $\mathcal L$.
 Let
 $$A=\bigoplus_{\mu \in \Delta \cup \{0\},\  { m  \in \Bbb Z}} \ A_\mu^m$$
 be the loop realization of $A$
 (which could be twisted)
 viewed as a Lie $1$-torus
 such that $\psi(A_{\pm\mu}^0)=Fe_{\pm\mu}$
 and
 $\psi(A_{\pm\nu}^{\mp 1})=Fx_{\pm\nu}$.
 Then, $\psi$ is graded
 relative to
the $\Bbb Z$-grading but also to
the double grading $\langle\Delta\rangle\times \Bbb Z$.
 Note that a centerless Lie torus is $\Bbb Z$-graded simple
 (see \cite[Lem.4.4]{Y1}).
Thus, the nontrivial $\Bbb Z$-graded ideal of $A$
 is exactly the $1$-dimensional center $Fc$,
 and
 the image of $\psi$ contains
 $\cup_{\mu\in\Delta}\mathcal L_\mu=\cup_{\mu\in\Delta}\oplus_{m\in\Bbb Z}\mathcal L_\mu^m$.
 Therefore, $\psi$ is onto
  since $\mathcal L$ is generated by $\cup_{\mu\in\Delta}\mathcal L_\mu$.
Thus, the induced graded isomorphism from
 the loop algebra $A/Fc$ onto $\mathcal L$ exists.
 Similarly, we obtain a graded isomorphism from the loop algebra $A/Fc$
 onto $\mathcal M$.
 Based on these isomorphisms, we obtain the graded isomorphism
 $\varphi$ described above.
 \qed

 \enskip

Thus, a centerless Lie 1-torus is isomorphic to a loop algebra,
and a Lie 1-torus with nontrivial center is isomorphic to a
derived affine Lie algebra, which has a 1-dimensional center.

For a Lie 1-torus $\mathcal L =  \oplus_{\mu \in \Delta \cup
\{0\},\  { m  \in \Bbb Z}} \  \mathcal L^m _\mu$, we have
\begin{equation}\label{1dim}
\text{ $\dim\mathcal L^m _\mu\neq 0$ (so $\dim\mathcal L^m
_\mu=1$) for all $\mu \in \Delta_{sh}$ and $m \in \Bbb Z$ (so
$S_\mu=\Bbb Z$), and}
\end{equation}
\begin{equation}\label{dagger2}
\text{ the center of $\mathcal L$ is equal to $[\mathcal L^m
_0,\mathcal L^{-m }_0]$ for any $0\neq m  \in \Bbb Z$.}
\end{equation}
This can be seen easily from the loop realization. Furthermore,
we have
\begin{equation}\label{dagger3}
\dim\sum_{m\in\Bbb Z}[\mathcal L_\mu^m,\mathcal L_{-\mu}^{-m}]=
\begin{cases}
1&\text{if $\mathcal L$ is loop}\\
2&\text{if $\mathcal L$ is derived affine}
\end{cases}
\end{equation}
since
$$
\sum_{m\in\Bbb Z}[\mathcal L_\mu^m,\mathcal L_{-\mu}^{-m}]=
\begin{cases}
F\mu^\vee&\text{if $\mathcal L$ is loop}\\
F\mu^\vee+Fc&\text{if $\mathcal L$ is derived affine}
\end{cases}
$$
for $\mu\in\Delta$ and a nontrivial central element $c$.

\begin{lemma}\label{lemma2.1}
 The center of a locally Lie 1-torus is at most 1-dimensional.
In particular, for
 a locally Lie 1-torus
 $\mathcal L=  \oplus_{\mu \in \Delta \cup \{0\},\  { m  \in \Bbb Z}} \  \mathcal L^m _\mu$,
$$\text{$\mathcal L$ has a 1-dimensional center
$\Longleftrightarrow$ $\mathcal L$ is a directed union of derived
affine Lie algebras},
$$
and
$$
\text{$\mathcal L$ is centerless $\Longleftrightarrow$ $\mathcal
L$ is a directed union of loop algebras}
$$
in the following sense:
$$\mathcal L=\bigcup_{\Delta'\subset\Delta}\
\mathcal L_{\Delta'},$$ where $\Delta'$ is a finite irreducible
full subsystem of $\Delta$ and $\mathcal L_{\Delta'}$ is the
homogeneous subalgebra of $\mathcal L$ generated by $\mathcal
L_\mu$ for $\mu \in\Delta'$, and $\mathcal L_{\Delta'}$ is a
derived affine Lie algebra if the center of $\mathcal L$ is
1-dimensional and a loop algebra if $\mathcal L$ is centerless.

In particular, the properties \eqref{1dim}, \eqref{dagger2}, and
\eqref{dagger3} given above hold in a locally Lie $1$-torus.

\end{lemma}

\proof Most of the statements follow from Lemma
\ref{locallystatements}. In fact, Lie $1$-tori are either derived
affine Lie algebras or loop algebras, and thus $\mathcal L$ is a
directed union of
 derived affine Lie algebras or loop algebras.
Considering the loop realization of a derived affine Lie algebra,
we find \eqref{1dim}.

Suppose that $C$ is a $2$-dimensional subalgebra contained in the
center. Then, a derived affine Lie algebra or a loop algebra
exists that contains $C$. However, this is impossible because
their centers have to be $1$-dimensional or zero.

Now, we need to show that derived affine Lie algebras and loop
algebras cannot appear simultaneously. If this is case, e.g.,
$\mathcal L'$ is a derived affine subalgebra and $\mathcal L''$ is
a loop subalgebra, then a derived affine or a loop algebra exists
that contains both $\mathcal L'$ and $\mathcal L''$ as graded
subalgebras. Suppose that $\mathcal L'$ and $\mathcal L''$ are
contained in $\mathcal L'''$ for a loop algebra $\mathcal L'''$.
However, this is impossible because of property \eqref{dagger2}
above. Thus, suppose that $\mathcal L'$ and $\mathcal L''$ are
contained in $\mathcal L'''$ for a derived affine Lie algebra
$\mathcal L'''$. Then, this is also impossible because of property
\eqref{dagger3} above. Thus, a locally Lie $1$-torus is either a
directed union of derived affine Lie algebras, such as $\mathcal
L_{da}$, or a directed union of loop algebras, such as $\mathcal
L_{lo}$. It is now clear that the center of $\mathcal L_{lo}$ is
zero. To show the $1$-dimensionality of the center of $\mathcal
L_{da}$, let $C'\ (\not= 0)$ be a finite dimensional central
subspace of a derived affine subalgebra of $\mathcal L_{da}$. For
any $\mu \in \Delta$ and $m  \in \Bbb Z$, a derived affine
subalgebra $M$ exists that contains $\mathcal L^m _\mu$ and $C'$.
Considering the loop realization of $M$, we find that $C'$ is the
$1$-dimensional center of $M$ and, in particular, $C'$ is the
$1$-dimensional center of $\mathcal L_{da}$.

Finally, let $\mathcal L$ be a locally Lie $1$-torus. Then,
\eqref{dagger3} is clear. To show \eqref{dagger2}, let
$Z:=[\mathcal L^k_0,\mathcal L^{-k}_0]$ for $0\neq k\in\Bbb Z$.
For any $z\in Z$, $\mu \in \Delta$, and $m  \in \Bbb Z$, a derived
affine subalgebra or a loop subalgebra exists that contains $z$
and $\mathcal L^m _\mu$, and $z$ is in the center of the
subalgebra (by \eqref{dagger2} for a Lie $1$-torus as given
above). Hence $[z,\mathcal L^m _\mu]=0$ for all $\mu \in \Delta$
and $m  \in \Bbb Z$. Therefore, $Z$ is contained in the center of
$\mathcal L$. Thus, $Z=0$ or $\dim Z=1$. If $Z=0$, then a loop
subalgebra exists, and thus $\mathcal L=\mathcal L_{lo}$. Hence,
$Z=0$ is the center of $\mathcal L$. If $\dim Z=1$, then $Z$ is
the center of $\mathcal L$ since the center of $\mathcal L$ is at
most $1$-dimensional. \qed

\enskip

For any two elements $x\otimes t^m$ and $y\otimes t^n$, in each
locally loop algebra $\mathcal L$, we define the new bracket on a
1-dimensional central extension
$$\tilde{\mathcal L}:=\mathcal L\oplus Fc$$ by
\begin{equation}\label{star}
[x\otimes t^m,y\otimes t^n]:=[x,y]\otimes
t^{m+n}+m(x,y)\delta_{m+n,0}c,
\end{equation}
where $(x,y)$ is the trace form $\tr(xy)$, or for type
$B_{\mathfrak I}^{(2)}$, the direct sum of the trace form and the
bilinear form on $\frak s$ is determined by the symmetric matrix
$s$ given above. Indeed, this gives a central extension since
$\mathcal L$ is a directed union of loop algebras and
$\tilde{\mathcal L}$ is a derived LALA, i.e., a 1-dimensional
central extension of a loop algebra.

\begin{lemma}\label{uclla}
A universal covering of a locally loop algebra is given by
\eqref{star}.
\end{lemma}
\proof Suppose that $\hat{\mathcal L}$ is a universal covering of
a locally loop algebra $\mathcal L$. We know that $\dim_F
Z(\hat{\mathcal L})\geq 1$ since $\tilde{\mathcal L}$ is a
covering. Therefore, if $\dim Z(\hat{\mathcal L})>1$, then a
covering $\mathcal L\oplus Fc_1\oplus Fc_2$ of $\mathcal L$
exists. Let $x_1, y_1, \ldots, x_m, y_m, u_1, v_1, \ldots, u_n,
v_n\in\mathcal L$ be such that $\sum_{i=1}^m[x_i,y_i]=c_1$ and
$\sum_{i=1}^n[u_i,v_i]=c_2$. Let $\mathcal L'$ be
 a loop subalgebra of $\mathcal L$ that
 contains $x_i, y_i, u_j, v_j$ for $1\leq i\leq m$ and $1\leq j\leq n$.
 Then,
 $\mathcal L'\oplus Fc_1\oplus Fc_2$ is
 perfect, and thus this is
 a covering of $\mathcal L'$.
However, a universal covering of a loop algebra has a
$1$-dimensional center, which is a contradiction. Hence, $\dim
Z(\hat{\mathcal L})=1$. However, it is then clear that
 $\hat{\mathcal L}\cong\tilde{\mathcal L}$
 since the unique morphism from $\hat{\mathcal L}$ onto $\tilde{\mathcal L}$
has to be one to one. \qed

\remark By Lemma \ref{lemma2.1}, a locally Lie 1-torus has at most
a 1-dimensional center. Thus, if we show that $\hat{\mathcal L}$
is a locally Lie $1$-torus, then we also obtain a proof of Lemma
\ref{uclla}. In fact, Neher showed that a universal covering of a
locally Lie torus is a locally Lie torus in general (see
\cite{Ne3} and \cite{NeS}).

\endremark

\medskip

Now, we classify locally Lie 1-tori. The method we use is derived
from [NS]. In particular, we show that there is only one locally
Lie 1-torus for each reduced root system extended by $\Bbb Z$. The
root systems extended by $\Bbb Z$ were classified in
\cite[Cor.15]{Y3} as the class of locally affine root
systems (more general results are given in 
\cite{LN2}). The
following is a list of all the reduced root systems extended by
$\Bbb Z$ of infinite rank:
\begin{align*}
&\text A_{\mathfrak I}\times\Bbb Z,\quad \text B_{\mathfrak
I}\times\Bbb Z,\quad \text C_{\mathfrak I}\times\Bbb Z,\quad \text
D_{\mathfrak I}\times\Bbb Z,
\\
& \big((\text B_{\mathfrak I})_{\sh}\times\Bbb
Z\big)\sqcup\big((\text B_{\mathfrak I})_{\lng}\times 2\Bbb
Z\big),\quad
\big((\text C_{\mathfrak I})_{\sh}\times\Bbb Z\big)\sqcup\big((\text C_{\mathfrak I})_{\lng}\times 2\Bbb Z\big),\\
&\bigg(\big((\text{BC}_{\mathfrak I})_{\sh}
\sqcup(\text{BC}_{\mathfrak I})_{\lng}\big)\times\Bbb Z\bigg)
\sqcup\big((\text {BC}_{\mathfrak I})_{\ex}\times (2\Bbb
Z+1)\big),
\end{align*}
where we write $\sqcup_{\mu\in\Delta}(\mu\times S_\mu)$ for
$\{S_\mu\}_{\mu\in\Delta}$, and for a subset $\Delta'$ of
$\Delta$, if all $S_\mu$'s for $\mu\in\Delta'$ are the same set
$S$, we write $\Delta'\times S$ instead of
$\sqcup_{\mu\in\Delta'}(\mu\times S_\mu)$. Furthermore, we simply
use a type instead of writing $\Delta$, e.g., $\text A_{\mathfrak
I}$ for $\Delta$ of type $\text A_{\mathfrak I}$.

We can see that these seven systems are the exact root systems of
the locally loop algebras introduced above, and thus we label each
system by
$$
\text A_{\mathfrak I}^{(1)},\quad \text B_{\mathfrak
I}^{(1)},\quad \text C_{\mathfrak I}^{(1)},\quad \text
D_{\mathfrak I}^{(1)},\quad \text B_{\mathfrak I}^{(2)},\quad
\text C_{\mathfrak I}^{(2)},\quad \text {BC}_{\mathfrak I}^{(2)}.
$$
We also use the label for the root system as the {\bf type} of  a
locally Lie 1-torus. We efer to the first four types as {\bf
untwisted} and the last three types as {\bf twisted}. Note that
\begin{equation*}
\text{all $S_\mu$'s are $\Bbb Z$ for the untwisted type, and in
general,}
\end{equation*}
\begin{equation}\label{shortZ}
\text{$S_\mu=\Bbb Z$ for a short root $\mu$.}
\end{equation}

\medskip

First, we provide the following lemma when {\bf $\Delta$ is
finite}. Suppose that $\Pi\subset\Delta$ is an integral base,
i.e., $\Delta\subset\langle\Pi\rangle$, and $\Pi$ is linearly
independent in the vector space that defines $\Delta$, where
$\langle \Pi \rangle$ is the additive subgroup generated by $\Pi$,
i.e., $\langle \Pi \rangle$ is the $\mathbb{Z}$-span of $\Pi$.
Note that
\begin{equation}\label{intbasered}
\Pi\subset\Delta^{\red}.
\end{equation}

\begin{lemma}\label{Lemma 2.3}
Let $\mathcal L =  \oplus_{\mu \in \Delta \cup \{0\}, { m \in \Bbb
Z}} \  \mathcal L^m_\mu$ and $\mathcal M =  \oplus_{\mu \in \Delta
\cup \{0\}, { m \in \Bbb Z}} \  \mathcal M^m_\mu$ be centerless
Lie 1-tori of the same type $\tilde\Delta$, where {\bf $\Delta$ is
finite}. Let $\Pi$ be an integral base of $\Delta$ that contains a
fixed short root $\nu\in\Delta$ if $\tilde\Delta$ is of the
untwisted type or a fixed short root $\nu\in\Delta$ if
$\tilde\Delta$ is of the twisted type. Let $0\neq x_\mu\in\mathcal
L_\mu^0$ and $0\neq y_\mu\in\mathcal M_\mu^0$ for each $\mu\in\Pi$
(see \eqref{intbasered}). Furthermore, let $0\neq x \in \mathcal
L_{\nu}^1$ and $0\neq y \in \mathcal M_{\nu}^1$ (see
\eqref{shortZ}).

Then, a unique isomorphism $\psi$ from $\mathcal L$ onto $\mathcal
M$ exists such that $\psi(x )=y $, $\psi(\mu^\vee_\mathcal L)=
\mu^\vee_\mathcal M$ and $\psi(x_\mu)=y_\mu$ for all $\mu\in\Pi$.

\end{lemma}

\proof By \eqref{dagger}, an isomorphism $\vvp:\mathcal
L\longrightarrow\mathcal M$ exists such that
$\vvp(\mu^\vee_\mathcal L)=\mu^\vee_\mathcal M$ and $\vvp(\mathcal
L_\mu^m)=\mathcal M_\mu^m$ for all $\mu\in\Delta$ and $m \in \Bbb
Z$. Hence,  we have $y =a \vvp(x )$ and $y_\mu=a_\mu\vvp(x_\mu)=$
for some $a$ and $a_\mu\in F^\times$. Let
$f:\langle\Pi\rangle_\Bbb Z\times\Bbb Z\longrightarrow F^\times$
be the group homomorphism of the abelian groups defined by
$f(\mu,0)=a_\mu$ and $f(0,1)=a $. Let $D_f$ be the diagonal linear
automorphism on $\mathcal M$ defined by $D_f(y)=f(\mu,m)y$ for
$y\in\mathcal M_\mu^m$. Then, $D_f$ is an automorhism of Lie
algebras. Indeed, $D_f([y,y'])=f(\mu+\mu',m+m')[y,y']
=f((\mu,m)+(\mu',m'))[y,y'] =f(\mu,m)f(\mu',m')[y,y']
=[f(\mu,m)y,f(\mu',m')y']=[D_f(y),D_f(y')]$ for $y\in\mathcal
M_\mu^m$ and $y'\in\mathcal M_{\mu'}^{m'}$. Hence,
$\psi:=D_f^{-1}\circ\vvp$ is the required isomorphism.

For the uniqueness, we first note that this isomorphism is unique
on $\mathcal L_{-\nu}^{-1}$ and $\mathcal L_{-\mu}^0$ for all
$\mu\in\Pi$ since $[\mathcal L_{\nu}^1,\mathcal
L_{-\nu}^{-1}]=F\nu^\vee$ (since $\mathcal L$ is centerless) and
$[\mathcal L_{\mu}^0,\mathcal L_{-\mu}^0]=F\mu^\vee$. Thus, it is
sufficient to show that $\mathcal L$ is generated by $\mathcal
L_{\nu}^{1}$, $\mathcal L_{-\nu}^{-1}$, and $\mathcal L_{\pm
\mu}^0$ for all $\mu\in\Pi$. However, by a standard argument (or
see \cite[Prop.9.9]{St}), $\mathcal L^0$ (= the finite-dimensional
split simple Lie algebra $\mathfrak g$) is generated by $\mathcal
L_{\pm \mu}^0$ for all $\mu\in\Pi$. Then, we can choose a root
base of $\Delta$ such that $\nu$ is the negative highest long root
if $\tilde\Delta$ is of the untwisted type or the negative highest
short root if $\tilde\Delta$ is of the twisted type. Using the
loop realization of $\mathcal L$, it is clear that $\mathcal L$ is
generated by $\mathcal L^0=\mathfrak g$ and $\mathcal
L_{\pm\nu}^{\pm 1}$. \qed

\enskip

Now, we can prove that there is a one to one correspondence
between the class of centerless locally Lie 1-tori and the class
of reduced root systems extended by $\Bbb Z$, and that locally
loop algebras exhaust all of the centerless locally Lie 1-tori.
Note that this method works for any cardinality of $\Delta$.

\begin{theorem}\label{Theorem 2.4}
Let $\mathcal L=\oplus_{\mu\in\Delta\cup\{0\}}\oplus_{m\in \Bbb
Z}\ \mathcal L_{\mu}^m$ be a locally Lie 1-torus of type
$\tilde\Delta$. If $\mathcal L$ is centerless, then $\mathcal L$
is graded isomorphic to the locally loop algebra of type
$\tilde\Delta$, and if $\mathcal L$ has a nontrivial center, then
$\mathcal L$ is graded isomorphic to a universal covering of the
locally loop algebra of type $\tilde\Delta$ given by \eqref{star}.
\end{theorem}

\proof First, it should be noted that we already know this theorem
for Lie 1-tori, i.e., the case where $\Delta$ is finite. In
addition, it is sufficient to show the case where $\mathcal L$ is
centerless (see Lemma \ref{uclla}), and thus we assume that
$\mathcal L$ is centerless. Let $\mathcal
M=\oplus_{\mu\in\Delta\cup\{0\}}\oplus_{m\in \Bbb Z}\ \mathcal
M_{\mu}^m$ be a locally loop algebra of type $\tilde\Delta$.
Furthermore, let $\tilde\Delta=\{S_{\mu}\}_{\mu\in\Delta}$.

Fix a long root $\nu$ if $\tilde\Delta$ is of the untwisted type,
or a short root $\nu$ if $\tilde\Delta$ is of the twisted type,
and let $0\neq x\in \mathcal L_{\nu}^1$ and $0\neq e_\nu\otimes
t\in \mathcal M_{\nu}^1$
 (see \eqref{shortZ}).
Let $\Pi$ be an integral base of $\Delta$ that contains $\nu$. Let
$0\neq x_\mu\in\mathcal L_\mu^0$ and $0\neq e_\mu\otimes
1\in\mathcal M_\mu^0$ for each $\mu\in\Pi$
 (see \eqref{intbasered}).
Then, we claim that the map $\psi:\mu^\vee_\mathcal L\mapsto
\mu^\vee_\mathcal M$ and $x_\mu\mapsto e_\mu\otimes 1$ for all
$\mu\in\Pi$, and $x\mapsto e_\nu\otimes t$ extends to an
isomorphism from $\mathcal L$ onto $\mathcal M$. Indeed, if we let
$\Gamma\subset\Pi$ be a finite irreducible subset that contains
$\nu$, then $\Gamma$ is an integral base of the irreducible root
system $\Delta_\Gamma:=\Delta\cap\langle\Gamma\rangle$.

Let $\tilde\Delta_{\Gamma}=\{S_{\mu}\}_{\mu\in\Delta_\Gamma}$ be
the root system extended by $\Bbb Z$. Let $\mathcal L_{\Gamma}$ be
the subalgebra determined by $\Delta_\Gamma$, i.e., the subalgebra
of $\mathcal L$ generated by $\mathcal L_\mu^m$ for all
$\mu\in\Delta_\Gamma$ and $m\in\Bbb Z$, which is a centerless Lie
1-torus of type $\tilde\Delta_{\Gamma}$ (see Lemma
\ref{lemma2.1}). Similarly, let $\mathcal M_{\Gamma}$ be the
subalgebra of $\mathcal M$ determined by $\Delta_\Gamma$. Then, by
Lemma \ref{Lemma 2.3}, a unique graded isomorphism $\psi_{\Gamma}$
from $\mathcal L_{\Gamma}$ onto $\mathcal M_{\Gamma}$ exists such
that $\psi_{\Gamma}(x_\mu)=e_\mu\otimes 1$ for all $\mu\in\Gamma$
and $x\mapsto e_\nu\otimes t$.

Suppose that $\Gamma_1,\Gamma_2\subset\Pi$ are finite irreducible
subsets that contain $\nu$ such that $\mathcal L_{\Gamma_1}\subset
\mathcal L_{\Gamma_2}$. Then, the uniqueness of the isomorphisms
$\psi_{\Gamma_1}$ and $\psi_{\Gamma_2}$ implies that they agree on
$\mathcal L_{\Gamma_1}$. Since $\mathcal L$ is the directed union
of the subalgebras $\mathcal L_{\Gamma}$ (\ $\Gamma\subset \Pi$ is a
finite irreducible subset), we can define an isomorphism
$\psi:\mathcal L\longrightarrow \mathcal M$ by
$\psi(x)=\psi_{\Gamma}(x)$ for $x\in\mathcal L_{\Gamma}$, which
has the required properties. \qed

\medskip

Note that in \eqref{star}, we defined the Lie bracket of a
universal covering of a locally loop algebra using a symmetric
bilinear form $(\cdot, \cdot)$ on a locally loop algebra. In
particular, we can write $(\cdot, \cdot)=\tr(\cdot,
\cdot)\otimes\epsilon(\cdot, \cdot)$, where $\epsilon(t^m,
t^n)=\delta_{m+n,0}$. In fact, it is easy to check that this form
is invariant, graded (as a form of a Lie torus defined in
\cite{Y2}), and nondegenerate. We simply refer to a {\bf form} for
a symmetric invariant graded bilinear form on a Lie $G$-torus. We
use the following lemma later.

\lemma\label{uniquenessform} A nonzero form on a locally Lie
$1$-torus exists. In addition, this form is unique up to a nonzero
scalar. In particular, a form of a locally loop algebra is equal
to $c(\cdot,\cdot)$ for some $c\in F$, where $(\cdot, \cdot)$ is
used in \eqref{star}.

\endlemma

\proof Only the uniqueness part is not clear (since we already use
a form in \eqref{star}). However, this form is unique up to a
scalar for a Lie $1$-torus (e.g., see \cite{Y2}). Thus, the
uniqueness follows from a local argument since a locally Lie
$1$-torus is a directed union of Lie $1$-tori. \qed

\enskip

\section{LALAs}
\label{sec:locallyaffine}

Let us recall LEALAs in \cite{MY}.  A subalgebra $\mathcal H$ of
a Lie algebra $\mathcal L$ is called ad-diagonalizable if
$$\mathcal L =\bigoplus_{\xi\in \mathcal H^*}\ \mathcal L _\xi,$$ where
$\mathcal H^*$ is the dual space of $\mathcal H$ and
$$\mathcal L _\xi=\{x\in \mathcal L \mid [h,x]=\xi(h)x\
\text{for all}\ h\in \mathcal H\}.$$ This decomposition is called
the {\bf root space decomposition} (of $\mathcal L$ with respect
to an ad-diagonalizable subalgebra $\mathcal H$). Note that an
ad-diagonalizable subalgebra $\mathcal H$ is automatically
abelian. To confirm this, we need the well-known fact that every
submodule of a weight module is also a weight module. We can use a
common trick to obtain the proof, e.g., as given in
\cite[Prop.2.1]{MP}, but they assumed that $\mathcal H$ is
abelian. To ensure that this assumption is unnecessary, we prove
it here. First we show that:
\begin{claim}
$\mathcal H = \oplus_{\xi \in \mathcal H^*}\ \mathcal H_\xi$,
where $\mathcal H_\xi = \mathcal L_\xi \cap \mathcal H$.
\end{claim}
\proof Suppose that $\mathcal H \not= \oplus_{\xi \in \mathcal
H^*} \mathcal H_\xi$. Then, $x \in \mathcal H$ exists such that
$x$ can be written as $x = x_1 + \cdots + x_n$ with $n > 1$, which
satisfies $x_i \in \mathcal L_{\xi_i} \setminus \mathcal H$ for
all $i$. Take $x \in \mathcal H$ among all of these elements such
that $n$ is minimal, and choose $h \in \mathcal H$ such that
$\xi_1(h) \not= \xi_2(h)$. Then, $x' := \ad h(x) - \xi_1(h)x =
(\xi_2(h) - \xi_1(h))x_2 + \cdots +(\xi_n(h) - \xi_1(h))x_n \in
\mathcal H$. This contradicts the minimality of $n$. Hence, we
have $\mathcal H = \oplus_{\xi \in \mathcal H^*} \mathcal H_\xi$.
\qed

\medskip

Now, suppose that $h \in \mathcal H_\xi$ and $h' \in \mathcal
H_{\xi '}$. Then, $[ h , h' ] = \xi'(h)h' = - \xi(h')h$. Hence, if
$h$ and $h'$ are linearly independent, then $[ h , h' ] = 0$.
Furthermore, we can see that $[h,h'] = 0$ if they are linearly
dependent. Thus, $\mathcal H$ is always abelian.

In particular, we have
$$\mathcal H = \mathcal H_0 \subset \mathcal L_0 = C_{\mathcal L}(\mathcal H),$$
where $C_{\mathcal L}(\mathcal H)$ is the centralizer of $\mathcal
H$ in $\mathcal L$.

\enskip

An element of the set
$$R=\{\xi\in \mathcal H^*\mid \mathcal L _\xi\neq 0\}$$ is called a {\bf root}.
(We do not call this $R$ a root system and we simply call it the
{\bf set of roots}.)

\enskip

Let $\mathcal L$ be a Lie algebra, $\mathcal H$ is a subalgebra of
$\mathcal L$, and $\mathcal B$ is a symmetric invariant bilinear
form of $\mathcal L $. A triple $(\mathcal L ,\mathcal H,{\mathcal
B})$ (or simply $\mathcal L $) is called a {\bf LEALA} if it
satisfies the following four axioms (we explain $R^\times$
shortly):
\begin{itemize}
\item[(A1)] $\mathcal H$ is ad-diagonalizable and
self-centralizing, i.e.,
$$\mathcal L =\bigoplus_{\xi\in \mathcal H^*} \mathcal L _\xi
\quad\text{and}\quad \mathcal H = \mathcal L _0;$$

\item[(A2)] ${\mathcal B}$ is nondegenerate;

\item[(A3)] $\ad x\in\End_F \mathcal L $ is locally nilpotent for
all $\xi\in R^\times$ and all $x\in \mathcal L _\xi$,

\item[(A4)] $R^\times$ is irreducible.

\end{itemize}

Moreover,
\begin{itemize}
\item[(i)] If $\mathcal H$ is finite-dimensional, then $\mathcal L
$ is called an {\bf EALA}.

\item[(ii)] If $R^\times=\emptyset$, then $(\mathcal L ,\mathcal
H,{\mathcal B})$ is called a {\bf null LEALA} (or a {\bf null
EALA} if $\mathcal H$ is finite-dimensional) or simply a {\bf null
system}. Note that if $R^\times=\emptyset$, then the axioms (A3)
and (A4) are empty statements.
\end{itemize}
Now, using (A1) and (A2), we find that $\mathcal B_{\mathcal L_\xi
\times \mathcal L_{-\xi}}$ is nondegenerate for all $\xi\in R$. In
particular,
$$
\text{$\mathcal B_{\mathcal H \times \mathcal H}$ is
nondegenerate.}
$$

\begin{lemma}\label{excludeax}
For each $\xi\in R$, a unique $t_\xi \in\mathcal H$ exists such
that $\mathcal B(h,t_\xi)=\xi(h)$ for all $h\in\mathcal H$.
\end{lemma}
\proof By the nondegeneracy of $\mathcal B_{\mathcal L_\xi \times
\mathcal L_{-\xi}}$, $x\in\mathcal L_\xi$ and $y\in\mathcal
L_{-\xi}$ exist such that $\mathcal B(x,y)=1$. Let
$t_\xi:=[x,y]\in\mathcal H$. Then,
$$\mathcal B(h,t_\xi) =\mathcal B(h,[x,y]) = B([h,x],y)
= \xi(h)\mathcal B(x,y) = \xi(h)$$ for all $h\in\mathcal H$. The
uniqueness of $t_\xi$ follows from the nondegeneracy of $\mathcal
B_{\mathcal H \times \mathcal H}$. \qed

\medskip

Using these $t_\xi$s, we can define an {\bf induced form on the
vector space spanned by $R$ over $F$}, which is simply denoted as
$(\cdot,\cdot)$, by
$$(\xi,\eta):=\mathcal B(t_\xi,t_\eta)$$
for $\xi,\eta\in R$. Note that the form $( \cdot , \cdot )$ is
well defined, which is easily confirmed by:
$$\begin{array}{lllll}
\mathcal{B} ( \sum_\xi p_\xi t_\xi , \sum_\eta q_\eta t_\eta ) & =
& \sum_\xi p_\xi \mathcal{B} ( t_\xi , \sum_\eta q_\eta t_\eta ) &
= & \sum_\xi p_\xi \xi ( \sum_\eta q_\eta t_\eta )\\ [0.2cm] & = &
\sum_\xi p_\xi' \xi ( \sum_\eta q_\eta t_\eta ) & = & \sum_\xi
p_\xi' \mathcal{B} (  t_\xi , \sum_\eta q_\eta t_\eta )\\ [0.2cm]
& = & \mathcal{B} ( \sum_\xi p_\xi' t_\xi , \sum_\eta q_\eta
t_\eta ) & = & \mathcal{B} ( \sum_\eta q_\eta t_\eta , \sum_\xi
p_\xi' t_\xi )\\ [0.2cm] & = & \sum_\eta q_\eta \mathcal{B}(
t_\eta , \sum_\xi p_\xi' t_\xi ) & = & \sum_\eta q_\eta \eta (
\sum_\xi p_\xi' t_\xi ) \\ [0.2cm] & = & \sum_\eta q_\eta' \eta (
\sum_\xi p_\xi' t_\xi ) & = & \sum_\eta q_\eta' \mathcal{B}(
t_\eta , \sum_\xi p_\xi' t_\xi ) \\ [0.2cm] & = & \mathcal{B} (
\sum_\eta q_\eta' t_\eta , \sum_\xi p_\xi' t_\xi ) & = &
\mathcal{B} ( \sum_\xi p_\xi' t_\xi , \sum_\eta q_\eta' t_\eta )
\end{array}$$
for $\sum_\xi p_\xi \xi = \sum_\xi p_\xi' \xi$ and $\sum_\eta
q_\eta \eta = \sum_\eta q_\eta' \eta$.

Now we call an element of
$$R^\times:=\{\xi\in R\mid (\xi,\xi)\neq 0\}$$
an {\bf anisotropic root}.  Axiom (A4) means that $R^\times=R_1\cup
R_2$ and $(R_1,R_2)=0$, which imply that $R_1=\emptyset$ or
$R_2=\emptyset$.

\remark Null systems have not been studied widely. In
\cite{AABGP}, they assumed that $R^\times\neq\emptyset$ for an
EALA. We also assume that $R^\times\neq\emptyset$ throughout this
study.
\endremark

\remark We note that there was one more axiom for a LEALA in
\cite{MY}, but we showed that axiom is unnecessary by Lemma
\ref{excludeax} above.
\endremark

We say that a triple $(\mathcal L ,\mathcal H,{\mathcal B})$ is
{\bf admissible} if it satisfies (A1) and (A2). A fundamental
property of admissible triples is as follows.

\lemma\label{FPAT} For $\xi\in R$ and all $x\in \mathcal L_\xi$
and $y\in\mathcal L_{-\xi}$, we have
\begin{equation}\label{fpat}
[x,y]=\mathcal B(x,y)t_\xi,
\end{equation}
where $t_\xi$ is defined in Lemma \ref{excludeax}.

\endlemma

\proof Let $h:=[x,y]-\mathcal B(x,y)t_\xi\in\mathcal H$. Then, for
all $h'\in\mathcal H$, we have
$$\mathcal B(h,h')=\mathcal B(x,[y,h'])-\mathcal B(x,y)\mathcal B(t_\xi,h')
=\mathcal B(x,y)\xi(h')-\mathcal B(x,y)\xi(h')=0.$$ Hence, by the
nondegeneracy of $\mathcal B_{\mathcal H \times \mathcal H}$, we
obtain $h=0$. \qed

\medskip

We can scale the above form $(\cdot,\cdot)$ by a nonzero scalar
such that $(\xi,\eta)\in\Bbb Q$ for all $\xi,\eta\in R^\times$
(see \cite[p.16]{AABGP} or \cite[\S 3]{MY}). Let $V$ be the $\Bbb
Q$-span of $R$, such as
$$V:=\spa_\Bbb Q R.$$
We showed the Kac Conjecture in \cite[Thm 3.10]{MY}, which states
that
\begin{equation}
\label{kac} \text{the scaled form $(\cdot,\cdot)$ on $V$ is
positive semidefinite, and $(R^0,V) = 0$,}
\end{equation}
where
$$R^0:=\{\xi\in R\mid (\xi,\xi)= 0\},$$
the set of {\bf isotropic roots} or {\bf null roots}. As a
corollary, $(W, R^\times)$ becomes a reduced locally extended
affine root system (LEARS), where $W=\spa_\Bbb Q R^\times$ (see
\cite[\S 4]{MY} and \cite{Y3}). We simply refer to $R$ as the set
of roots, but we refer to $R^\times$ as a {\bf  LEARS}.
This  $R^\times$
satisfies the fundamental properties of classical finite
irreducible root systems, locally finite irreducible root systems,
and affine root systems in the sense of Macdonald \cite{Ma},  or
extended affine root systems in the sense of Saito \cite{S}. We do
not recall the definition of LEARS because is is not needed in
this study. The reader can find the precise definition in
\cite{Y3}.

The dimension of the radical of $V$ is called the {\bf null
dimension} for a LEALA. If the additive subgroup of $V$ generated
by $R^0$ is free, we call the rank the {\bf nullity} of a LEALA.
Thus, we only use the term {\it nullity} when $\langle R^0\rangle$
is a free abelian group. 
\remark\label{onlyforfree} Of course,
there is a notion of rank for non-free abelian groups, but to be
consistent with the original theory of EALAs, as given in
\cite{AABGP}
and  \cite{Ne2}, we assume that $\langle R^0\rangle$
is free for nullity. 
Thus, if we say that a LEALA $\mathcal L$ has
nullity, this means that $\langle R^0\rangle$ is a free abelian
group. (In \cite{MY}, we used the term {\it null rank} for
nullity, and {\it nullity} for null dimension, but we have changed
these terms to maintain consistency with the notion of nullity in
\cite{Ne2}.)
\endremark

The {\bf core} of a LEALA $\mathcal L$, denoted by $\mathcal L_c$,
is the subalgebra of $\mathcal L$ generated by the root spaces
$\mathcal L_\alpha$ for all $\alpha\in R^\times$. Then, by the Kac
Conjecture \eqref{kac}, $\mathcal L_c$ is an ideal of $\mathcal
L$. If the centralizer of $\mathcal L_c$ in $\mathcal L$ is
contained in $\mathcal L_c$, then $\mathcal L$ is called {\bf
tame}. Note that the core is zero for a null system (since it is
generated by an empty set), so a null system is not tame.

Now, as mentioned earlier, $(W, R^\times)$ is a reduced LEARS.
Thus, by \cite{Y3}, a locally finite irreducible root system
$\Delta$ and
 a {\bf reflectable section} $W'$ of $W$
exist such that $\Delta^{\red}$
is contained in $R^\times\cap W'$. 
In particular, $W'$ is a
complement of $\rad W$, such as $W=W'\oplus\rad W$, where $\rad W$
is the radical of $W$ relative to the defining positive
semidefinite form of the LEARS $(W, R^\times)$. Moreover, a family
of subsets $\{S_\mu\}_{\mu\in\Delta}$ of $\rad W$ indexed by
$\Delta$ exists such that
\begin{equation}\label{rtimesgrading}
R^\times=\bigcup_{\mu\in\Delta}\ (\mu+S_{\mu}),
\end{equation}
and $\{S_\mu\}_{\mu\in\Delta}$ is a reduced root system extended
by $G=\langle \cup_{\mu\in\Delta}\ S_{\mu}\rangle$, as defined in
Section 2. We note that
$$\rad W=(\rad V)\cap W,$$
by the Kac Conjecture \eqref{kac}.

We can give the graded structure of the core $\mathcal L_c$ from
\eqref{rtimesgrading}. For each $\mu\in\Delta$ and $g\in G$, if
$g\in S_{\mu}$, where we let
$$(\mathcal L_c)_\mu^g:
=\mathcal L_c\cap \mathcal L_{\mu+g},$$ and if $g\notin S_{\mu}$,
where we let $(\mathcal L_c)_\mu^g:=0$. 
Then, we can easily show
that
$$\mathcal L_c
=\bigoplus_{\mu\in\Delta\cup\{0\}}\ \bigoplus_{g\in G}\ (\mathcal
L_c)_\mu^g,$$ where $(\mathcal L_c)_0^g:=\sum_{\mu\in\Delta}\
\sum_{g=h+k}\ [(\mathcal L_c)_\mu^h,(\mathcal L_c)_{-\mu}^k]$, and
that
\begin{equation}\label{coreistorus}
\text{ $\mathcal L_c$ is a locally Lie $G$-torus of type
$\Delta$,}
\end{equation}
or more precisely, of type $\{S_\mu\}_{\mu\in\Delta}$.
Furthermore, if we let
$$\mathcal L_c^g
:=\bigoplus_{\mu\in\Delta\cup\{0\}}\ (\mathcal L_c)_\mu^g,$$ 
then
we obtain a $G$-graded Lie algebra
$$
\mathcal L_c=\bigoplus_{g\in G}\ \mathcal L_c^g.
$$

\enskip

Next, we note some properties related to $R^0$ for LEALAs. As
mentioned in Section 2, $S_\mu=S_\nu$ if $\mu$ and $\nu$ are the
same length for $\mu,\nu\in\Delta$. If we let
$$S:=S_\mu$$ for a short root $\mu$,
then $S$ contains all $S_\nu$, as in \eqref{shortbig}, and $S$
satisfies $0\in S$ and $2S-S\subset S$. In addition,
$$
\text{$S$ spans $\rad W$}
$$
(see [Thm 8, Y2]).

\begin{lemma}\label{tames+s}
Let $\mathcal L$ be a LEALA. Then, $S+S\subset R^0$, and $S+S=R^0$
if $\mathcal L$ is tame.

Moreover, we have $\rad V=\spa_\Bbb Q R^0$. In particular, if
$\mathcal L$ has nullity, then (nullity of $\mathcal L$) = (null
dimension of $\mathcal L$).
\end{lemma}
\proof The first statement follows from \eqref{S4} in Section 1,
but we present this for convenience with respect to the next
statement. Let $s,s'\in S$. Then, $\mathcal L_{-\mu+s}\neq 0$ and
$\mathcal L_{\mu+s'}\neq 0$ for $\mu\in\Delta_{sh}$, and
$[\mathcal L_{-\mu+s},\mathcal L_{\mu+s'}]\neq 0$, by
$\sll_2$-theory. (Consider the $\sll_2$-subalgebra generated by
$\mathcal L_{\mu-s}$ and $\mathcal L_{-\mu+s}$, and let it act on
$\mathcal L_{\mu+s'}$.) 
Therefore, $0\neq[\mathcal
L_{-\mu+s},\mathcal L_{\mu+s'}]\subset\mathcal L_{s+s'}$ and hence
$s+s'\in R^0$. Thus, $S+S\subset R^0$.

Suppose that $\mathcal L$ is tame. Let $\sigma\in R^0$. If
$\alpha+\sigma\notin R$ for all $\alpha\in R^\times$, then
$\mathcal L_\sigma$ centralizes the core, and thus $\mathcal
L_\sigma$ is in the core. Therefore, $\mathcal
L_\sigma=\sum_{\mu\in\Delta,\ s+s'=\sigma} [\mathcal L_{\mu+s},
\mathcal L_{-\mu+s'}]$, and thus $\sigma=s+s'$ for some $s,s'\in
S_\mu=S_{-\mu}\subset S$. However, $0\neq \mathcal
L_{\mu+s}=\mathcal L_{\mu-s'+\sigma}$ and $0\neq\mathcal
L_{\mu-s'}$ since $-s'\in S_{\mu}$. Therefore, $\mu-s'+\sigma\in
R$ with $\mu-s'\in R^\times$, which is a contradiction. Thus,
$\alpha\in R^\times$ exists such that $\alpha+\sigma\in R$. (This
property is called {\bf nonisolated}. Therefore, we have shown
that any isotropic root is nonisolated if $\mathcal L$ is tame.)
Note that $\alpha=\mu+s$ for some $\mu\in\Delta$ and $s\in S$.
Hence, $s+\sigma\in S$, so $\sigma\in S-S=S+S$. Thus, $S+S= R^0$.

For the last statement, it is sufficient to show that $\rad
V\subset V^0:=\spa_\Bbb Q R^0$ (the other inclusion is clear).
Since $V=W+V^0$ (where $W = \spa_\Bbb Q R^\times$), it is
sufficient to show that $(\rad V)\cap W=\rad W\subset V^0$.
However, this is clear since $\rad W=\spa_\Bbb Q S$, as above.
\qed

\medskip

Note that if we put
\begin{equation*}\label{}
R^0_c:=\{\delta\in R^0\mid \mathcal L_\delta\cap \mathcal L_c\neq
0\},
\end{equation*}
then \eqref{S4} in Section 2 means that we always have
\begin{equation}\label{coreanisoroots}
R^0_c=S+S.
\end{equation}

\begin{remark}\label{Remark 3.2}
(1) In fact, the rank of $\langle R^0\rangle$ as a torsion-free
abelian group is always of the null dimension since the null
dimension is now simply the $\Bbb Q$-dimension of $\spa_\Bbb Q
R^0$ by Lemma \ref{tames+s}.

(2) There are notions of null dimension and nullity for LEARS
$(W,R^\times)$, i.e., (null dimension of $R^\times$) := $\dim\rad
W$ and (nullity of $R^\times$) := $\rank\langle S\rangle$ if
$\langle S\rangle$ is free (see \cite{Y3}).

For example, if $S=\Bbb Q$, then the null dimension is $1$, and
the rank (= the largest cardinality of linearly independent
elements over $\mathbb{Z}$) of $\Bbb Q$ as a torsion-free group is
also $1$. However,  we do not say that the nullity is $1$ when
$S=\Bbb Q$.

In general, (null dimension of $\mathcal L$) $\geq$ (null
dimension of $R^\times$). If $\mathcal L$ has nullity, then so
does $R^\times$and (nullity of $\mathcal L$) $\geq$ (nullity of
$R^\times$) since any subgroup of a free abelian group is free
(e.g., see \cite{G}). If $\mathcal L$ is tame, then (null
dimension of $\mathcal L$) $=$ (null dimension of $R^\times$), and
if $\mathcal L$ has nullity, then
$$
\text{ (nullity of $\mathcal L$) = (null dimension of $\mathcal
L$) = (nullity of $R^\times$) = (null dimension of $R^\times$) }
$$
since $S+S=R^0$.

\end{remark}

Now, we present some basic properties of the center of a LEALA.
\begin{proposition}\label{centerofleala}
Let $(\mathcal L,\mathcal H, \mathcal B)$ be a LEALA over $F$ with
the center $Z(\mathcal L)$, and $R^0$ is the set of isotropic
roots of $\mathcal L$. Then:

(1) We have
$$\sum_{\delta\in R^0} Ft_{\delta}\subset Z(\mathcal L) \subset\mathcal H,$$
where $t_{\delta}$ is a unique element in $\mathcal H$ defined by
\eqref{fpat} in Lemma \ref{FPAT}.

(2) Let $\mathcal L_c$ be the core of $\mathcal L$ and
$R^0_c=\{\delta\in R^0\mid \mathcal L_\delta\cap\mathcal L_c\neq
0\}$. Then, for $\delta\in R^0_c$, we have $t_\delta\in\mathcal
L_c$ and
 $$\sum_{\delta\in R_c^0} Ft_{\delta}= Z(\mathcal L_c)\cap\mathcal H
 \subset  Z(\mathcal L).$$

(3) Let $R^\times$ be the set of anisotropic roots of $\mathcal L$
(which is a LEARS). Let $m=\dim_\Bbb Q ({\rm rad}\ W)$ be the null
dimension of $R^\times$, i.e., the dimension of the radical of the
induced form from $\mathcal B$ on $W=\spa_\Bbb Q R^\times$. Then,
$m\geq \dim_{F}  \big(Z(\mathcal L_c)\cap\mathcal H\big)$, and if
$m\geq 1$, then $\dim_{F}   \big(Z(\mathcal L_c)\cap\mathcal
H\big)\geq 1$. Hence,  $m=1$ implies that $\dim_{F}
\big(Z(\mathcal L_c)\cap\mathcal H\big)=1$ and $\dim_{F}
Z(\mathcal L)\geq 1$.

(4) If $\mathcal L$ is tame, then $\sum_{\delta\in R^0}
Ft_{\delta}= \sum_{\delta\in R_c^0} Ft_{\delta} =Z(\mathcal L_c)
\cap\mathcal H=  Z(\mathcal L)$.

Furthermore, let $n$ be the null dimension of $\mathcal L$, i.e.,
$n=\dim_\Bbb Q\spa_\Bbb Q R^0$. Then, $m=n\geq \dim_{F} Z(\mathcal
L)$. Moreover, if $n\geq 1$, then $\dim_{F}  Z(\mathcal L)\geq 1$.
Hence, $n=1$ implies that $\dim_{F}  Z(\mathcal L)=1$.

\end{proposition}
\proof (1):  Since each $\delta$ is an isotropic root, we have
$[t_{\delta},x]=0$ for any root vector $x\in\mathcal L_\xi$. In
fact, $[t_{\delta},x]=\xi(t_{\delta})x=(\xi,\delta)x=0$ since
$\delta$ is in the radical of the form (see \eqref{kac}). Hence,
$[t_{\delta},\mathcal L]=0$, i.e., $t_{\delta}\in Z(\mathcal L)$.
Thus, $\sum_{\delta\in R^0} Ft_{\delta}\subset Z(\mathcal L)$. The
second inclusion is clear due to the fact that $\mathcal H$ is
self-centralizing.

(2):
 For $\delta\in R^0_c$, let $0\neq x\in\mathcal L_\delta\cap\mathcal L_c$.
Then, $t_\delta=[x,y]$ for some $y\in\mathcal L_{-\delta}$, and
hence $t_\delta\in\mathcal L_c$ since $\mathcal L_c$ is an ideal.
Thus, $\sum_{\delta\in R_c^0} Ft_{\delta}\subset\mathcal
L_c\cap\mathcal H$, and by (1), 
we obtain $\sum_{\delta\in R_c^0}
Ft_{\delta}\subset\mathcal L_c\cap Z(\mathcal L) \subset
Z(\mathcal L_c)$. Therefore, we obtain $\sum_{\delta\in R_c^0}
Ft_{\delta}\subset Z(\mathcal L_c)\cap\mathcal H$.

For the other inclusion, let $x\in Z(\mathcal L_c)\cap\mathcal H$.
Since
$$\mathcal L_c\cap\mathcal H
=\sum_{\xi\in R^\times}[\mathcal L_\xi, \mathcal L_{-\xi}]
+\sum_{\delta\in R^0_c}[\mathcal L_\delta, \mathcal
L_{-\delta}],$$ 
we can write
$$x=\sum_{\xi\in R^\times}a_{\xi} t_{\xi}
+\sum_{\delta\in R^0_c}a_{\delta}t_{\delta},$$ where $a_{\xi},
a_{\delta}\in F$. Let $\Delta\subset R^\times$ be a locally finite
irreducible root system determined by a reflectable section of
$\bar R^\times$ and $S$ is a reflection space for a short root in
$\Delta$. Then, we know that $R^\times\subset\Delta+S$ and
$R^0_c=S+S$ (see \eqref{coreanisoroots}). Thus, we obtain
\begin{align*}
x&=\sum_{\alpha\in\Delta, \delta'\in S} a_{\alpha+ \delta'}
t_{\alpha+ \delta'}
+\sum_{\delta\in S+S}a_{\delta}t_{\delta}\\
&=\sum_{\alpha\in\Delta, \delta'\in S}  (a_{\alpha+ \delta'}
t_{\alpha}
+a_{\alpha+ \delta'} t_{ \delta'})+\sum_{\delta\in S+S}a_{\delta}t_{\delta}\\
&=\sum_{\alpha\in\Delta, \delta'\in S}   a_{\alpha+ \delta'}
t_{\alpha} +\sum_{\alpha\in\Delta, \delta'\in S} a_{\alpha+
\delta'} t_{ \delta'} +\sum_{\delta\in S+S}a_{\delta}t_{\delta},
\end{align*}
and hence,
$$
y:=\sum_{\alpha\in\Delta, \delta'\in S}   a_{\alpha+ \delta'}
t_{\alpha}\in Z(\mathcal L_c).
$$
However, $y\in\frak h\subset \frak g$, and since $\frak g$ is a
locally finite split {\bf simple} Lie algebra, then $y$ has to be
$0$. Therefore,
$$
x=\sum_{\alpha\in\Delta, \delta'\in S} a_{\alpha+ \delta'} t_{
\delta'} +\sum_{\delta\in S+S}a_{\delta}t_{\delta}\in
\sum_{\delta\in R_c^0} Ft_{\delta},
$$
and we obtain $Z(\mathcal L_c) \cap \mathcal{H} \subset
\sum_{\delta\in R_c^0} Ft_{\delta}$. Hence, $\sum_{\delta\in
R_c^0} Ft_{\delta} = Z(\mathcal L_c) \cap \mathcal{H}$. The second
inclusion follows from (1).

(3): We know that $R^\times\subset\Delta+S$ and $m=\dim_\Bbb Q
({\rm rad}\ W) = \dim_\Bbb Q\spa S$. 
However, since $R_c^0=S+S$,
we have $m=\dim_\Bbb Q\spa R_c^0$. Using $\mathcal B$, we define
an injective linear map
$$\varphi_\mathcal B\ :\ \mathcal H \longrightarrow \mathcal H^*,$$
where $\varphi_\mathcal B(h) \in \mathcal H^*$ for $h \in \mathcal
H$ is given by $\varphi_\mathcal B(h)(h') = \mathcal B(h,h')$ for
all $h' \in \mathcal H$. Note that $\varphi_\mathcal B(t_\mu) =
\mu$ for $\mu \in R$, where $t_\mu \in \mathcal H$ satisfies
$[x,y] = \mathcal B(x,y) t_\mu$ for $x \in \mathcal L_\mu$ and $y
\in \mathcal L_{-\mu}$. Set $\mathcal H^\circ = {\rm im}
\varphi_\mathcal B \subset \mathcal H^*$. 
If we put
$$t = \varphi_\mathcal B^{-1}\ : \ \mathcal H^\circ \longrightarrow \mathcal H$$
and $t_\nu = t(\nu) = \varphi_\mathcal B^{-1}(\nu) \in \mathcal H$
for $\nu \in \mathcal H^\circ$, 
then we find that $t_{\nu + \nu'}
= t_\nu + t_{\nu'}$ for all $\nu, \nu' \in \mathcal H^\circ$, and
$t_{a \nu} = at_\nu$ for $\nu \in \mathcal H^\circ$ and $a \in F$.
Since $R \subset \mathcal H^\circ$, there is a one to one
correspondence
$$\{\delta\in R_c^0\}
\leftrightarrow \{t_{\delta}\}_{\delta\in R_c^0},$$ and, in
particular, we can see that
$t_{\delta+\delta'}=t_{\delta}+t_{\delta'}$ for $\delta,\delta'\in
R_c^0$ and $t_{a\delta}=at_{\delta}$ for $\delta \in R_c^0$ and
$a\in F$. Thus, $m=\dim_\Bbb Q\sum_{\delta\in R_c^0} \Bbb Q
t_{\delta} \geq \dim_{F}\sum_{\delta\in R_c^0} F t_{\delta}
=\dim_{F}  \big(Z(\mathcal L_c)\cap\mathcal H\big)$. Finally, if
$m\geq 1$, then $0\neq\delta\in R_c^0$ exists and thus
$t_{\delta}\neq 0$. 
Thus, $Ft_{\delta}\neq 0$, and hence we obtain
the last statement.

(4): We have $R^0=S+S=R^0_c$ since $\mathcal L$ is tame (see Lemma
\ref{tames+s}). Hence, $\sum_{\delta\in R^0} Ft_{\delta}=
\sum_{\delta\in R_c^0} Ft_{\delta}$. 
Furthermore, by (2), we
already have $\sum_{\delta\in R_c^0} Ft_{\delta} =Z(\mathcal L_c)
\cap\mathcal H\subset  Z(\mathcal L)$. Moreover, for $x\in
Z(\mathcal L)$, we have $x\in Z(\mathcal L_c)$ since $\mathcal L$
is tame. Hence, $Z(\mathcal L_c) \cap\mathcal H=  Z(\mathcal L)$.
The remaining assertions follow from the fact that $R^0=R^0_c$
using (3) and Lemma \ref{tames+s}. \qed

\remark There are examples of a tame LEALA or EALA where the
nullity is $\infty$ but the center is simply $1$-dimensional. For
example, $\mathcal L=\sll_2(\Bbb C[t_i^{\pm 1}]_{i\in\Bbb
N})\oplus \mathbb{C}c\oplus \mathbb{C}d$ is a tame EALA over $\Bbb
C$ of type $\text A_1$, where $d=\sum_{i=1}^\infty a_id_i$ with
degree derivation $d_i=t_i\frac{\partial}{\partial t_i}$, and
$\{a_i\}_{i\in\Bbb N}\subset\Bbb C$ is linearly independent over
$\Bbb Q$. This $\mathcal L$ has nullity of $\infty$ but the center
is equal to $Fc$. Note that the Cartan subalgebra $\mathcal H$ of
$\mathcal L$ is simply $3$-dimensional (for details, see
\cite[Rem.5.2(2)]{MY}).

\endremark

\begin{lemma}\label{derofcenterlesscore}
Let $(\mathcal L,\mathcal H, \mathcal B)$ be a tame LEALA. Then,
we have the natural embedding
$$\mathcal L/Z(\mathcal L_c)
\hookrightarrow \der_F\mathcal L_c \quad\text{and}\quad \mathcal
L/Z(\mathcal L_c) \hookrightarrow \der_F\big(\mathcal
L_c/Z(\mathcal L)\big).$$ (Note that $Z(\mathcal L)=Z(\mathcal
L_c)\cap\mathcal H$ by Proposition \ref{centerofleala}.)

In particular, if $\mathcal N$ is a complement of the core
$\mathcal L_c$, i.e., $\mathcal L = \mathcal L_c \oplus \mathcal
N$, then $\mathcal N$ can be identified with a subspace of
$\oder_F\big(\mathcal L_c/Z(\mathcal L)\big)$, i.e., the outer
derivations.
\end{lemma}
\proof Since $\mathcal L_c$ is an ideal of $\mathcal L$, we find
that the restriction $\ad x\mid_{\mathcal L_c}$ for $x \in
\mathcal L$ is in $\der_F \mathcal L_c$ and that $Z(\mathcal L_c)$
is an ideal of $\mathcal L$. Let
$$f\ :\ \mathcal L \longrightarrow \der_F \mathcal L_c$$
be the induced map obtained by this restriction. Since $\mathcal
L$ is tame, we have $\ker f  = C_\mathcal L(\mathcal L_c) =
Z(\mathcal L_c)$. Hence, we obtain the first embedding. In
addition, $\ad x\mid_{\mathcal L_c}$ induces a derivation of
$\mathcal L_c/Z(\mathcal L)$ since $Z(\mathcal L) = Z(\mathcal
L_c) \cap \mathcal H \subset Z(\mathcal L_c)$. Let
$$f' : \mathcal L \longrightarrow \der_F
\big(\mathcal L_c/Z(\mathcal L)\big)$$ be the induced map. Let $x
\in \ker  f'$. Then, we have $[ x , \mathcal L_c ] \subset
Z(\mathcal L)=Z(\mathcal L_c) \cap \mathcal H$. For any
$w\in\mathcal L_c$, since $\mathcal L_c$ is perfect, we can write
$w = \sum_{ i } [u_i,v_i]$ for some $u_i,v_i \in \mathcal L_c$.
Then, $[x,y] =\sum_i [[x,u_i],v_i] + \sum_i[u_i,[x,v_i]] = 0$, and
thus $[ x , \mathcal L_c ] = 0$. Hence, $\ker  f' \subset
Z(\mathcal L_c)$. It is clear that $Z(\mathcal L_c) \subset \ker
f'$. Thus, $\ker  f' = Z(\mathcal L_c)$, and hence we obtain the
second embedding.

For the second assertion, suppose that $\ad x$ for $x \in \mathcal
N$ is inner in $\der_F\big(\mathcal L_c/Z(\mathcal L)\big)$, i.e.,
$\ad x=\ad y$ on $\mathcal L_c/Z(\mathcal L)$ for some
$y\in\mathcal L_c$. Then, we have $[x-y,\mathcal L_c]\subset
Z(\mathcal L)$. However, since $\mathcal L_c$ is perfect, for $w=
\sum_{ i } [u_i,v_i]$ ($u_i,v_i \in \mathcal L_c$), we have
$[x-y,w]= \sum_{ i }[[x-y,u_i],v_i] + \sum_{ i }[u_i,
[x-y,v_i]]=0$. Hence, $x-y\in C_\mathcal L(\mathcal
L_c)=Z(\mathcal L_c)$ by tameness. In particular, $x-y\in\mathcal
L_c$, but 
$x\in\mathcal L_c$, which forces $x$ to be $0$.
Therefore, $\ad x$ is an outer derivation of $\mathcal
L_c/Z(\mathcal L)$. \qed

\enskip

Finally, we give some definitions for later use.

\definition
Let $V$ be a vector space over $\Bbb Q$, and $G$ is an additive
subgroup of $V$. Let
$$\mathcal A=\bigoplus_{g\in G}\ \mathcal A^g$$ be a $G$-graded algebra.
Define a linear transformation $d_i$ on $\mathcal A$ by
$$
d_i(a_g)=g_ia_g
$$
for $a_g\in\mathcal A^g$, where $g_i$ is the $i$-coordinate of $g$
obtained by a fixed basis of $V$. Note that $d_i$ depends on a
basis of $V$. Then, $d_i$ is a derivation of $\mathcal A$ where we
have
$$d_i(a_ga_h)=(g_i+h_i)a_ga_h
=g_ia_ga_h+h_ia_ga_h =d_i(a_g)g_h+a_gd_i(a_h)$$ for
$a_h\in\mathcal A^h$ and $h\in G$. We refer to each $d_i$ as an
{\bf $i$-th coordinate-degree derivation}.

If $\dim_FV=1$, then $d_1$ is simply called a {\bf degree
derivation}.
\enddefinition

We define a standard LEALA.
\begin{definition}\label{defstand}
If a LEALA $\mathcal L$ contains all coordinate-degree derivations
that act on the $G$-graded core, i.e., a locally Lie $G$-torus,
then $\mathcal L$ is called {\bf standard}. This concept depends
on the $G$-graded structure of the core, which is not unique.
Thus, when we use this term more precisely, we say that $\mathcal
L$ is standard (or {\bf non-standard}) relative to the locally Lie
$G$-torus.
\end{definition}

We define the minimality of a LEALA (see \cite{N2} and Remark
\ref{neebsmin}).

\definition\label{ourmin}
A LEALA $\mathcal L$ is called {\bf minimal} if $\mathcal L$ is
the only LEALA that contains $\mathcal L_c$ and which is contained
in $\mathcal L$ (equivalently, if there is no LEALA $\mathcal L'$
that satisfies $\mathcal L_c \subset \mathcal L' \subsetneq
\mathcal L$).
Note that if the nullity is positive, then $\mathcal L_c$ is never
a LEALA. Thus, if $\mathcal L$ has positive nullity and $\mathcal
L_c$ is a hyperplane in $\mathcal L$ (i.e., $\dim \mathcal
L/\mathcal L_c = 1$), then $\mathcal L$ is minimal.
\enddefinition

\example Let $\mathcal L^{ms} =\sll_\Bbb N(F[t^{\pm 1}])\oplus Fc
\oplus Fd^{(0)}$ (as explained in the Introduction), $\mathcal L_1
=\sll_\Bbb N(F[t^{\pm 1}])\oplus Fc \oplus F(e_{11}+d^{(0)})$,
where $e_{11}$ is the matrix unit of size $\Bbb N$ (only
$(1,1)$-entry is $1$ and all the other entries are $0$), and
$\mathcal L_2 =\sll_\Bbb N(F[t^{\pm 1}])\oplus Fc \oplus
F(p+d^{(0)})$, where
$$\displaystyle{
p= \dia  (1, \frac{1}{2}, \frac{1}{3},  \ldots,  \frac{1}{n},
\ldots )} \quad\text{(a diagonal matrix of size $\Bbb N$)}.
$$
Then, these three Lie algebras are all minimal LEALAs. (See
Definition \ref{defLALA}. In fact, these are minimal LALAs.) In
addition, $\mathcal L^{ms}$ is standard, but $\mathcal L_1$ and
$\mathcal L_2$ are not standard.

In Example \ref{nonstaiso}, we show that $\mathcal L_1$ is
isomorphic to $\mathcal L^{ms}$. We note that the concept of
standard is not an isomorphic invariant because it depends on the
grading of the core. In Example \ref{counterexofst}, we also show
that $\mathcal L_2$ is not isomorphic to $\mathcal L^{ms}$.
\endexample

\enskip

\section{LEALAs of nullity $0$}
\label{}

We classified LEALAs of nullity 0 in \cite[Thm 8.7]{MY}. Now, we
describe the tame LEALAs of nullity 0 in a slightly different
manner compared with the description in \cite{MY}.

Let $M:=M_{\frak I}(F)$, $M_{2\frak I+1}(F)$ or $M_{2\frak I}(F)$
 be the space of matrices of an infinite size $\frak I$,
 $2\frak I+1$, or $2\frak I$, respectively,
 and $T_{\frak I}$,
 $T_{2\frak I+1}$, 
 or $T_{2\frak I}$ is the subspace of $M$ that
comprises diagonal matrices. Let $T'$ be a complement of
$F\iota_\frak I$ in $T_{\frak I}$, where $\iota_{\frak I}$ is the
identity matrix such that
$$T_{\frak I}=T' \oplus F\iota_\frak I.$$
Then, the following list comprises infinite-dimensional {\bf
maximal tame LEALAs of nullity $0$}. (The term ``maximal'' is used
in the usual sense, i.e., no tame LEALA contains each listed LEALA
of each type.)

\begin{itemize}
\item Type $\text A_{\mathfrak I}$:
\begin{equation}\label{slt'}
\text{ $\sll_{\mathfrak I}(F)+ T'$ with a Cartan subalgebra $T'$}
\end{equation}
(Note that $T'$ is the unique modulo $F\iota_\frak I$. In
addition, see Remark \ref{modLEALA0} and Lemma \ref{centginA}),

\item Type $\text B_{\mathfrak I}$: $\text o_{2{\mathfrak
I}+1}(F)+T^+$ with a Cartan subalgebra $T^+$, where
$$T^+:=\{x\in T_{2\frak I+1}\mid sx=-xs\},$$

\item Type $\text C_{\mathfrak I}$: $\text {sp}_{2{\mathfrak
I}}(F)+T^+$ with a Cartan subalgebra $T^+$, where
$$T^+:=\{x\in T_{2\frak I}\mid sx=-xs\},$$

\item Type $\text D_{\mathfrak I}$: $\text o_{2{\mathfrak
I}}(F)+T^+$ with a Cartan subalgebra $T^+$, where
$$T^+:=\{x\in T_{2\frak I}\mid sx=-xs\},$$
and each matrix $s$ is the same as $s$ defined in
\eqref{definings}.
\end{itemize}

We note that $F\iota_\frak I$ is the center of $\sll_{\mathfrak
I}(F)+ T_\frak I$, and that
$$\sll_{\mathfrak I}(F)+ T'
\cong\big(\sll_{\mathfrak I}(F)+ T_\frak I \big)/F\iota_\frak I$$
for any $T'$. It is sometimes better to embed $T'$ into $T_\frak
I/F\iota_\frak I$.

As with locally finite split simple Lie algebras, each of type
$\text B_\frak I$, $\text C_\frak I$, or $\text D_\frak I$
 is the fixed algebra of $\sll_{2\frak I+1}(F)+T_{2{\mathfrak I}+1}$
or $\sll_{2\frak I}(F)+T_{2{\mathfrak I}}$ by the automorphism
$\sigma$ defined in \eqref{originalinv}. This is why we write
$T^+$ because this is the eigenspace of eigenvalue $1$ of
$\sigma$. We write the eigenspace of eigenvalue $-1$ of $\sigma$
as $T^-$.

Any subalgebra of a maximal tame LEALA of nullity $0$ that
contains each locally finite split simple Lie algebra is a tame
LEALA of nullity $0$. Thus, let $\mathcal L$ be a tame LEALA of
nullity $0$. Then,
\begin{itemize}
\item Type $\text A_{\mathfrak I}$: $\sll_{\mathfrak I}(F)\subset
\mathcal L\subset \sll_{\mathfrak I}(F)+ T'$
 with a Cartan subalgebra $\mathcal L\cap T'$,

\item Type $\text B_{\mathfrak I}$: $\text o_{2{\mathfrak
I}+1}(F)\subset \mathcal L\subset \text o_{2{\mathfrak
I}+1}(F)+T^+$
 with a Cartan subalgebra $\mathcal L\cap T^+$,

\item Type $\text C_{\mathfrak I}$: $ \text {sp}_{2{\mathfrak
I}}(F)\subset \mathcal L\subset \text {sp}_{2{\mathfrak
I}}(F)+T^+$
 with a Cartan subalgebra $\mathcal L\cap T^+$,

\item Type $\text D_{\mathfrak I}$: $ \text o_{2{\mathfrak
I}}(F)\subset \mathcal L\subset \text o_{2{\mathfrak I}}(F)+T^+$
 with a Cartan subalgebra $\mathcal L\cap T^+$.
\end{itemize}
(We describe the defining bilinear form $\mathcal B$ shortly.)

\enskip

We consider $T_{{\mathfrak I}}$ more carefully. Set
$$
T_{{\mathfrak I}}^{as}=\{d\in T_{{\mathfrak I}}\mid \ \text{$d$ is
{\bf almost scalar}\}} = \{ d \in T_\frak I \mid d -a\iota_\frak I
\in {\rm gl}_\frak I(F)\ \mbox{for some}\ a \in F \},
$$
i.e., $d$ has {\bf only finitely many different diagonal entries
from the identity} $\iota_\frak I$. Clearly, $T_{{\mathfrak
I}}^{as}$ is a subspace of $M_{{\mathfrak I}}(F)$.

\lemma\label{twodec} Let $\frak h$ be the diagonal subalgebra of
$\sll_\frak I(F)$. Then, we have
\begin{equation}\label{asdetail}
T_{{\mathfrak I}}^{as}=\frak h\oplus F\iota_\frak I\oplus Fe_{jj},
\end{equation}
where $e_{jj}$ is the matrix in $M_{{\mathfrak I}}(F)$ such that
the $(j,j)$-entry is $1$ and all the other entries are $0$ for any
fixed index $j\in\frak I$. In particular, we have
$$\gl_\frak I(F)=\sll_\frak I(F)\oplus Fe_{jj}$$
for any $j\in\frak I$.

Furthermore, let $I$ be any finite subset of $\frak I$, and
$\iota_I:=\sum_{i\in I}e_{ii}$. Then, we have
\begin{equation}\label{asaface}
T_{{\mathfrak I}}^{as}=\frak h_I\oplus F\iota_I\oplus T_{\frak
I\setminus I}^{as},
\end{equation}
where $\frak h_I$ is the subspace of $\frak h$ such that all
$(k,k)$-components of $k\in\frak I\setminus I$ are $0$, and
$T_{\frak I\setminus I}^{as}$ is the subspace of $T_{{\mathfrak
I}}^{as}$ such that all $(i,i)$-components of $i\in I$ are $0$.

Moreover, we have
\begin{equation}\label{asafacemore}
T_{{\mathfrak I}}=\frak h_I\oplus F\iota_I\oplus T_{\frak
I\setminus I},
\end{equation}
where $T_{\frak I\setminus I}$ is the subspace of $T_{{\mathfrak
I}}$ such that all $(i,i)$-components of $i\in I$ are $0$.

\endlemma
\proof It is clear that $T_{{\mathfrak I}}^{as}\supset\frak
h\oplus F\iota_\frak I\oplus Fe_{jj}$. For the other inclusion,
let $x\in T_{{\mathfrak I}}^{as}$. Then, $a\in F$ exists such that
$y:=x-a\iota_\frak I\in T_{\frak I}\cap \gl_{\frak I}(F)$. Hence,
$y=y-\tr(y)e_{jj}+\tr(y)e_{jj}$ and note that
$h:=y-\tr(y)e_{jj}\in\frak h$. Thus, $x=h+a\iota_\frak
I+\tr(y)e_{jj}\in\frak h\oplus F\iota_\frak I\oplus Fe_{jj}$. This
completes the description of (29).

For the second decomposition (30), we have $T_{{\mathfrak
I}}^{as}= T_I\oplus T_{\frak I\setminus I}^{as}$, where $T_I$ is
the subset of $T_{{\mathfrak I}}^{as}$ such that all
$(k,k)$-components of $k\in\frak I\setminus I$ are $0$. However,
it is then easy to see that $T_I=\frak h_I\oplus F\iota_I$. The
last decomposition (31) is now clear. \qed

\enskip

We have not mentioned the defining bilinear form $\mathcal B$ of a
tame LEALA $\mathcal L$ of nullity $0$. Thus, as described in
\cite{MY}, let $\frak g$ be one of the the locally finite split
simple Lie algebra
$$\sll_{\mathfrak I}(F),
\ \text o_{2{\mathfrak I}+1}(F), \ \text {sp}_{2{\mathfrak I}}(F)
\quad\text{or}\quad \text o_{2{\mathfrak I}}(F),$$ contained in
$\mathcal L$, as defined above. The restriction $\mathcal
B_{\mathcal L\times \frak g}$ of $\mathcal B$ to the space
$\mathcal L\times \frak g$ is a nonzero scalar multiple of the
trace form, and the remaining part, i.e., the restriction to
$\frak C\times \frak C$, where $\frak C$ is a complement of $\frak
g$, can be any symmetric bilinear form.

In fact, in 
\cite{MY}, we did not state clearly why the
restriction $\mathcal B_{\mathcal L\times \frak g}$ of $\mathcal
B$ is a nonzero scalar multiple of the trace form. However, this
follows from the perfectness of $\frak g$ and the invariance of
$\mathcal B$. We summarize this phenomenon in a slightly more
general setup. Let us refer to a symmetric invariant bilinear form
simply as a {\bf form} for convenience.

\lemma\label{perfectidealpr} Let $L$ be a Lie algebra with a form
$B$ and let $\frak g$ be a perfect ideal of $L$. If any form of
$\frak g$ is equal to $B':=B\mid_{\frak g\times\frak g}$ up to a
scalar, then any invariant bilinear form on $L\times\frak g$ or on
$\frak g\times L$ is equal to $B\mid_{L\times\frak g}$ or
$B\mid_{\frak g\times L}$ up to a scalar. In this case,
``invariant on $L\times\frak g$'' means that
$B([x,y],z)=B(x,[y,z])$ for $x,y\in L$ and $z\in\frak g$.
\endlemma
\proof Let $E$ be an invariant bilinear form on $L\times\frak g$.
For $x\in L$ and $y\in\frak g$, since $y=\sum_{i}[u_i, v_i]$ for
some $u_i,v_i\in\frak g$, then we have
$$
E(x,y)=E(x,\sum_{i}[u_i, v_i])= c\sum_{i}B'([x,u_i], v_i)=
c\sum_{i}B([x,u_i], v_i)= cB(x,\sum_{i}[u_i, v_i])=cB(x,y)$$ for
some $c\in F$. We can prove the result for $\frak g\times L$ in a
similar manne. \qed

\enskip

Recall that the associative algebra
\begin{equation*}\label{}
M_\frak I^{\fin}(F)= \{x\in M_\frak I(F)\mid \text{{\small each
row and column of} $x$ {\small have only finitely many nonzero
entries}}\}
\end{equation*}
is a Lie algebra under the commutator. Using the matrix $s\in
M_\frak I^{\fin}(F)$ defined in \eqref{definings}, we can define
an automorphism of $M_\frak K^{\fin}(F)$, where $\frak K=2\frak I$
or $2\frak I+1$, by the same definition of $\sigma$ in
\eqref{originalinv}. We also denote the automorphism by $\sigma$.
Thus, each fixed Lie algebra $M^{\fin}_\frak K(F)^\sigma$ contains
a locally finite split simple Lie algebra $\frak g:=\sll_\frak
K(F)^\sigma$.

\lemma\label{partiallyunique} Let $L$ be any subalgebra of
$M^{\fin}_\frak I(F)$, and let $M$ be any subalgebra of $\gl_\frak
I(F)$. Then, the trace form $\tr$ on $L\times M$ and $M\times L$
is well defined and it is invariant.

Hence, if $L$ contains $\sll_\frak I(F)$, then any invariant
bilinear form on $L\times\sll_\frak I(F)$ or on $\sll_\frak
I(F)\times L$ is equal to $c\tr$ for some $c\in F$. In particular,
$\sll_\frak I(F)$ is a perfect ideal of $L$.

Moreover, if $L$ is a subalgebra of $M^{\fin}_\frak K(F)^\sigma$
that contains $\frak g=\sll_\frak K(F)^\sigma$, then $\frak g$ is
a perfect ideal of $L$, and any invariant bilinear form on
$L\times\frak g$ or on $\frak g\times L$ is equal to $c\tr$ for
some $c\in F$.

\endlemma
\proof Since $xy\in\gl_\frak I(F)$ for $x\in L$ and $y\in M$, then
the trace form $\tr(xy)$ is well defined. To show the invariance,
i.e., $\tr([A,B]y)=\tr(A[B,y])$ for $A,B\in L$ and $y\in M$, it is
sufficient to show this for $y=e_{ij}$ (the matrix unit of
$(i,j)$-component).

Let $A=(a_{mn})$, $B=(b_{mn})$, and $C=(c_{mn}) = [ A , B ]$.
Then, $c_{mn} = \sum_{k}(a_{mk}b_{kn} - b_{mk}a_{kn})$ and
$\tr([A, B]y)= \tr ( (c_{mn}) e_{ij} ) = c_{ji}
=\sum_{k}(a_{jk}b_{ki}-b_{jk}a_{ki})$ and
$$\tr(A[B,y])=\tr ((a_{mn})\big(\sum_{m}b_{mi}e_{mj}-\sum_{n}b_{jn}e_{in} ) )
=\sum_{k}(a_{jk}b_{ki}-a_{ki}b_{jk}).$$ Therefore, the trace form
is invariant. We can prove this for the case where $M\times L$ in
a similar manner. We note that $\sll_\frak I(F)$ or $\frak g$ is a
perfect ideal of $L$. By \cite[Lem. II.11]{NS}, any form on
$\sll_\frak I(F)$ is equal to $c\tr$ for some $c\in F^\times$.
Therefore, the second and last statements follow from Lemma
\ref{perfectidealpr}. \qed

\begin{remark}
We employ the notation given in Lemma \ref{partiallyunique}. We
can identify $M^{\fin}_\frak I(F)$ with the derivation algebra
$\der\big(\gl_\frak I(F)\big)$, and $M^{\fin}_\frak K(F)^\sigma$
with the derivation algebra $\der\frak g$ (see \cite{N1}).
\end{remark}

\enskip

Suppose that $\mathcal B$ is a symmetric invariant bilinear form
on
$$\mathcal M_\frak I:=\sll_\frak I(F)+T_\frak I.$$
Then, by Lemma \ref{partiallyunique},
 the restriction of $\mathcal B$ to
 $\mathcal M_\frak I\times\sll_\frak I(F)$ or $\sll_\frak I(F)\times \mathcal M_\frak I$ is
 equal to
$c\tr$ for some $c\in F$. We claim that such a form $\mathcal B$
does exist. Therefore, we select any complement $\frak h^c$ of
$\frak h$ in $T_\frak I$, i.e.,
$$T_\frak I=\frak h^c\oplus \frak h.$$
Let
$$\psi:\frak h^c\times \frak h^c\longrightarrow F
$$
be an arbitrary symmetric bilinear form. Now, we define a
symmetric bilinear form $\mathcal B$ on $\mathcal M_\frak I$ as
$$
\mathcal B(x, y) =\psi(x, y)
$$
on $\frak h^c$,
 and $c\tr$
 on $\mathcal M_\frak I\times\sll_\frak I(F)$
 and $\sll_\frak I(F)\times \mathcal M_\frak I$.
To show that $\mathcal B$ is invariant, we prove the following.

\begin{claim}\label{tonull0case}
Let $x\in T_\frak I\setminus F\iota_\frak I$ and $y_k\in
\sll_\frak I(F)$ for $k=1,2,\ldots, r$. Then,
 a finite subset $I$ of $\mathfrak I$,\ \
$0\neq h\in\frak h$
 and $g\in T_\frak I$
 exist such that
$y_k\in\sll_I(F)$ for all $k$, $h\in\frak h_I$,
 $$x=h+g,\quad
[x, y_k]=[h, y_k] \quad\text{and}\quad \mathcal B(x, y_k)=\mathcal
B(h, y_k)$$
 for all $k$.
 Moreover,
 if $\mathcal{B}$ is nontrivial, then
 $y\in \sll_\frak I(F)$ and $h'\in\frak h$
exist such that $[x, y]\neq 0$ and
\begin{equation}\label{nonzeropairnull0}
\mathcal B(x, h')\neq 0.
\end{equation}

\end{claim}
\proof Let $I$ be a finite subset of $\frak I$ such that
$y_k\in\sll_I(F)$ for all $k$. Moreover, if the $I\times I$-block
submatrix of $x$ is a scalar matrix, then we enlarge $I$ until the
$I\times I$-block submatrix of $x$ is not a scalar matrix. For
$I$, by \eqref{asafacemore} in Lemma \ref{twodec},
$0\neq h\in\frak h_I$
exists such that $x=h+b\iota_I+x'$ for some $b\in F$ and $x'\in
T_{\frak I\setminus I}$. Put $g:=b\iota_I+x'$. Then, clearly $[g,
y_k]=0$. In addition, we have $\mathcal B(g, y_k)=c\tr (gy_k)
=cb\tr (y_k)=0$ since $\tr(y_k)=0$.

To show the second statement, it is sufficient to select $y \in
\sll_I(F)$ and $h' \in \frak h_I$ such that $[h,y]\neq 0$ and
$\tr(hh')\neq 0$. \qed

\begin{claim}\label{}
$\mathcal B$ is invariant.
\end{claim}
\proof It is sufficient to consider the case that involves some
elements in $\frak h^c$. Since $\frak h^c$ is an abelian
subalgebra, the case that involves three elements in $\frak h^c$
is clear.

For the case that involves one element in $\frak h^c$, let $x\in
\frak h^c$ and $y,z\in\sll_\frak I(F)$. Then, it is sufficient to
show that
$$
\mathcal B([x, y],z)=\mathcal B(x, [y, z]).
$$
If $x\in F\iota$, then both sides are clearly $0$. Thus, by Claim
\ref{tonull0case}, we can change $x$ into $h$ for $y$ and $[y,z]$
such that $\mathcal B([x, y],z)=\mathcal B([h, y],z)$ and
$\mathcal B(x, [y, z])=\mathcal B(h, [y, z])$. This follows from
the invariance on $\sll_\frak I(F)$.

The case that involves two elements in $\frak h^c$ can be shown in
a similar manner. Let $x,y\in \frak h^c$ and $z\in\sll_\frak
I(F)$. Then, it is sufficient to show that
$$
\mathcal B(x, [y, z])=0 \quad\text{and}\quad \mathcal B([x,
z],y)=\mathcal B(x, [z, y]).
$$
Again, if $x$ or $y\in F\iota$, then both sides of both equations
are clearly $0$. Thus, by Claim \ref{tonull0case}, the left-hand
side of the first equation is equal to $\mathcal B(h, [h', z])$
for some $h,h'\in\frak h_I$, and this is equal to $0$ by the
invariance on $\sll_\frak I(F)$. For the second equation, change
$x$ into $h$ for $z$ and $[z,y]$ such that (LHS) $=\mathcal B([h,
z],y)$ and (RHS) $=\mathcal B(h, [z, y])$. However, these are
equal according to the case involving one element, as described
above. Thus, we have proved that the symmetric bilinear form
$\mathcal B$ is invariant. \qed

\medskip

The radical of $\mathcal B$ is contained in $F\iota_\frak I$
whenever the restriction to $\sll_\frak I(F)$ is not zero. In
fact, this follows from \cite[Lem. 8.5]{MY} since the center of
$\mathcal M_\frak I = {\rm sl}_\frak I(F) + T_\frak I$ is equal to
$F\iota_\frak I$. However, for convenience, we show this directly.
First, let us mention the graded structure of $\mathcal M_\frak
I$.

Let $\frak g:=\sll_\frak I(F)$ and let $\frak g=\frak
h\oplus\big(\bigoplus_{\mu\in\text A_\frak I\subset \frak h^*}\
\frak g_\mu\big)$ be the root-space decomposition of $\frak g$
relative to $\frak h$. We extend each root $\mu\in\frak h^*$ to an
element in $T_\frak I^*$ as follows.

Let $\text A_\frak I=\{\pm(\epsilon_i-\epsilon_j)\mid i,j\in\frak
I\}$, where $\epsilon_i$ is the linear form of $\gl_\frak I(F)$
determined by $e_{kl}\mapsto \delta_{lk}\delta_{ki}$. Since an
element $p\in T_\frak I$ can be written as $p= \dia
(a_{ii})_{i\in\frak I}$, we can define $\epsilon_i(p)=a_{ii}$. In
this manner, we can embed $\text A_\frak I$ into $T_\frak I^*$.
Thus, $\mathcal M:=\mathcal M_\frak I$ has the root-space
decomposition
$$\mathcal M
=\bigoplus_{\mu\in T_\frak I^*}\ \mathcal M_{\mu}$$ relative to
$T_\frak I$, where $\mathcal M_{\mu}=\frak g_\mu$ for $\mu\neq 0$
and $\mathcal M_0=T_\frak I$, and $\mathcal M_{\mu}=0$ if
$\mu\notin\text A_\frak I$. This is an $\langle\text A_\frak
I\rangle$-graded Lie algebra, and $\mathcal B$ is {\bf graded} in
the sense that $\mathcal B(\mathcal M_\xi,\mathcal M_\eta)=0$,
unless $\xi+\eta=0$ for all $\xi,\eta\in\text A_\frak I$. In
general, a symmetric invariant bilinear form on a Lie algebra with
a root-space decomposition relative to a subalgebra is graded.

In particular, the radical of $\mathcal B$ is graded. Thus, we can
check the nondegeneracy for each homogeneous element. The elements
of degree $\mu\in\text A_\frak I$ cannot be in the radical by
Lemma \ref{partiallyunique}. For the elements of degree $0$, the
only candidate is an element in $F\iota_\frak I$ by
\eqref{nonzeropairnull0}, which implies that the radical of
$\mathcal{B})$ is contained in $F\iota_\frak I$.

\medskip

Therefore, we have the following.

\begin{lemma}
Let $\mathcal B$ be nontrivial. Then, the radical of $\mathcal B$
is equal to $F\iota_\frak I$ if $\mathcal B(\iota_\frak
I,\iota_\frak I)= 0$, and $\mathcal B$ is nondegenerate if
$\mathcal B(\iota_\frak I,\iota_\frak I)\neq 0$. \qed
\end{lemma}

Thus, for any symmetric bilinear form $\psi$ on $\frak h^c$ with
the radical $F\iota_\frak I$, the quotient Lie algebra $\mathcal
M_\frak I/F\iota_\frak I$ with the induced form $\bar B$ is a
LEALA of type $\text A_\frak I$ of nullity $0$. Note that
$\mathcal M_\frak I/F\iota_\frak I$ is isomorphic to $\mathcal
M_\frak I':=\sll_\frak I(F)\oplus \frak t$, where $\frak t$ is a
complement of $\frak h\oplus F\iota_\frak I$ in $T_\frak I$.
Conversely, if $\psi'$ is any symmetric bilinear form on $\frak
t$, we can define a symmetric nondegenerate invariant form
$\mathcal B'$ on $\mathcal M_\frak I'$ as described above, and
$\mathcal M_\frak I'$ is isomorphic to $\mathcal M_\frak
I/F\iota_\frak I$. By a similar argument, we can say that a LEALA
of type $\text A_\frak I$ of nullity $0$ is isomorphic to a
subalgebra of $\mathcal M_\frak I/F\iota_\frak I$ that contains
$\sll_\frak I(F)=(\sll_\frak I(F)+F\iota_\frak I)/F\iota_\frak I$
with the induced form $\bar B$.

\example\label{nullity0ex} The centerless Lie algebra $\gl_\frak
I(F)=\sll_\frak I(F)\oplus Fe_{jj}$ is an example of a LEALA of
type $\text A_\frak I$ of nullity $0$, where $e_{jj}$ is the
matrix unit for $j\in\frak I$. However, $\gl_n(F)=\sll_n(F)\oplus
Fe_{jj}$ has the center $F\iota_n$ if $j \in \{ 1,2, \cdots , n
\}$, where $\iota_n$ is the identity matrix on $\gl_n(F)$, and
this is a non-tame EALA of nullity $0$.

Suppose that $\mathcal B$ is a nondegenerate form on $\gl_\frak
I(F)$. Then, $\mathcal B$ is a nonzero scalar multiple $c\in F$ of
the trace form,  except on $Fe_{jj}\times Fe_{jj}$, by Lemma
\ref{partiallyunique}. Conversely, we can take any value to
$\mathcal B(e_{jj},e_{jj})$ and extend a nondegenerate form to
$\gl_\frak I(F)$.

For the finite case where $\gl_n(F)=\sll_n(F)\oplus Fe_{jj}$,
suppose that $\frak B$ is a nondegenerate form on $\gl_n(F)$.
Since $\rad\frak B$ is in the center of $\gl_n(F)$, we find that
$\frak B$ is nondegenerate $\Longleftrightarrow$ $\frak
B(\iota_n,\iota_n)\neq 0$. Moreover, this is equivalent to
\begin{equation}
\label{nondegcond} \frak B(e_{jj},e_{jj})\neq \frac{n-1}{n}c.
\end{equation}
In fact, consider the expression
$\iota_n=\iota_n-ne_{jj}+ne_{jj}$, where we note that
$\tr(\iota_n-ne_{jj})=0$. Since $x:=\iota_n-ne_{jj}\in\sll_\frak
I(F)$, we have $\frak B(\iota_n,\iota_n)=$
\begin{align*}
\frak B(x+ne_{jj},x+ne_{jj})
&=\frak B(x,x)+2n\frak B(x,e_{jj})+n^2\frak B(e_{jj},e_{jj})\\
&=c\tr(x^2)+2nc\tr(xe_{jj})+n^2\frak B(e_{jj},e_{jj})\\
&=
c\tr(\iota_n-2ne_{jj}+n^2e_{jj})+2nc\tr(e_{jj}-ne_{jj})+n^2\frak B(e_{jj},e_{jj})\\
&=c(n-2n+n^2)+2nc(1-n)+n^2\frak B(e_{jj},e_{jj})\\
&=cn-cn^2+n^2\frak B(e_{jj},e_{jj}).
\end{align*}
Hence, $\frak B(\iota_n,\iota_n) = 0$ if and only if $n^2\frak
B(e_{jj},e_{jj})=c(n^2-n)$, and thus \eqref{nondegcond} holds.
\endexample

\remark\label{modLEALA0} In the classification of tame LEALAs of
nullity $0$ of type $\text A_{\frak I}$ in \cite{MY}, we select
$\sll_{\frak I}(F)\oplus \frak t'$ for a complement $\frak t'$ of
$T_{\mathfrak I}^{as}$ in $T_\frak I$ as the maximal one. However,
a subalgebra of the bigger Lie algebra $\sll_{\frak I}(F)\oplus
\frak t$ defined above is actually a maximal tame LEALA of nullity
$0$, which is shown essentially by the following lemma.

\lemma\label{centginA} Let $p\in T_\frak I$. Suppose that
$[p,\sll_\frak I(F)]=0$. Then, $p\in F\iota_\frak I$. In
particular, $\sll_\frak I(F)+T'$ is a tame LEALA of nullity $0$
for any complement $T'$ of $F\iota_\frak I$ in $T_\frak I$.
\endlemma
\proof Let $I$ be any finite subset of $\frak I$. Decompose
$p=h\oplus s\iota_I\oplus q$ in $T_{{\mathfrak I}}=\frak h_I\oplus
F\iota_I\oplus T_{\frak I\setminus I}$ for some $s\in F$ (see
\eqref{asafacemore}). Since $\sll_I(F)\subset\frak g$ and
$[\iota_I,\sll_I(F)]=0$, we have $0=[p,\sll_I(F)]=[h,\sll_I(F)]$.
Hence, $h=0$, and thus $p=s\iota_I\oplus q$. For a different
subset $I'$, we have $p=s'\iota_{I'}\oplus q'$. However, for
$I''=I\cup I'$, we have $p=s''\iota_{I''}\oplus q''$. Since $I,
I'\subset I''$, we have $s=s'=s''$. Therefore, $p=s\iota_\frak I$.
\qed

\endremark

\enskip

Now, we consider the forms on the other types $\text B_\frak I$,
$\text D_\frak I$ and $\text C_\frak I$. Let $\mathcal B$ be a
symmetric invariant form on
$$\mathcal M_\frak K=\sll_\frak K(F)+T_\frak K$$
such that the restriction to $\sll_\frak K(F)$ is not zero, where
$\frak K=2\frak I$ or $2\frak I+1$. Let $\mathcal M_\frak
K^\sigma$ be the fixed algebra by the automorphism $\sigma$
defined above with the restricted form $\mathcal B^\sigma$. Then,
$\mathcal B^\sigma$ is still invariant, and by Lemma
\ref{partiallyunique}, the restriction to $\sll_\frak K(F)^\sigma$
is equal to $c\tr$ for some $c\in F^\times$.

Moreover, $\mathcal B^\sigma$ is nondegenerate. This follows from
\cite[Lem. 8.5]{MY} since $\mathcal M_\frak K^\sigma$ has a
trivial center. We can also show this using the following lemma,
which is similar to Lemma \ref{twodec}. Recall that $T^+$ denotes
the eigenspace of eigenvalue $+1$ of $\sigma$, and $T^-$ is the
eigenspace of eigenvalue $-1$ of $\sigma$.
\begin{lemma}\label{otherdecnull0}
Let $I$ be any finite subset of $\frak I$ and fix some index
$i_0\in I$. Then, we have
$$
T_{2\frak I}^+=\frak h_{2I}^+\oplus T_{2\frak I\setminus 2I}^+
\quad\text{and}\quad T_{2\frak I}^- =\frak h_{2I}^-\oplus
F(e_{i_0i_0}+e_{\frak I+i_0,\frak I+i_0})\oplus T_{2\frak
I\setminus 2I}^-,
$$
where $\frak h_{2I}^+$ or $\frak h_{2I}^-$ is a subset of $\frak
h^+$ or $\frak h^-$ such that all $(k,k)$ and $(\frak I+k,\frak
I+k)$ components for $k\in\frak I\setminus I$ are $0$, and
$T_{2\frak I\setminus 2I}^+$ or $T_{2\frak I\setminus 2I}^-$ is a
subset of $T_{2\frak I}^+$ or $T_{2\frak I}^-$ such that all
$(i,i)$ and $(\frak I+i,\frak I+i)$ components for $i\in I$ are
$0$.

Furthermore, we have
$$
T_{2\frak I+1}^+ =\frak h_{2I+1}^+\oplus T_{(2\frak I+1)\setminus
(2I+1)}^+ \quad\text{and}\quad T_{2\frak I+1}^- =\frak
h_{2I+1}^-\oplus Fe_{2\frak I+1,2\frak I+1}\oplus T_{(2\frak
I+1)\setminus (2I+1)}^-,
$$
where $\frak h_{2I+1}^+$ or $\frak h_{2I+1}^-$ is a subset of
$\frak h^+$ or $\frak h^-$ such that the $(k,k)$ and $(\frak
I+k,\frak I+k)$ components of all $k\in\frak I\setminus I$ are
$0$, and $T_{(2\frak I+1)\setminus (2I+1)}^+$ or $T_{(2\frak
I+1)\setminus (2I+1)}^-$ is a subset of $T_{2\frak I+1}^+$ or
$T_{2\frak I+1}^-$ such that the $(2\frak I+1,2\frak I+1)$
component and the $(i,i)$ and $(\frak I+i,\frak I+i)$ components
of all $i\in I$ are $0$.

Moreover, we have
\begin{equation}\label{asaface2null0}
T_{2\frak I}^- =\frak h_{2I}^-\oplus F\iota_{2I}\oplus T_{2\frak
I\setminus 2I}^- \quad\text{and}\quad T_{2\frak I+1}^- =\frak
h_{2I+1}^-\oplus F\iota_{2I+1}\oplus T_{(2\frak I+1)\setminus
(2I+1)}^-.
\end{equation}

\end{lemma}
\proof All of these statements are clear except
\eqref{asaface2null0}. To show this, we consider the two
equalities
$$T_{2I}^-=\frak h_{2I}^-\oplus
F(e_{i_0i_0}+e_{\frak I+i_0,\frak I+i_0}) \quad\text{and}\quad
T_{2I+1}^-=\frak h_{2I+1}^-\oplus Fe_{2\frak I+1,2\frak I+1},$$
where $T_{2I}^-$ is a subset of $T_{2\frak I}^-$ such that $(i,i)$
and the $(\frak I+i,\frak I+i)$ components of all $i\in \frak
I\setminus I$ are $0$, and $T_{2I+1}^-$ is a subset of $T_{2\frak
I+1}^-$ such that $(i,i)$ and the $(\frak I+i,\frak I+i)$
components of all $i\in \frak I\setminus I$ are $0$. However, as
in the proof of Lemma \ref{twodec}, for $y\in T_{2I}^-$ or
$T_{2I+1}^-$, it follows from the equation that
$$y=y-\frac{1}{2}\tr(y)(e_{i_0i_0}+e_{\frak I+i_0,\frak I+i_0})
+\frac{1}{2}\tr(y)(e_{i_0i_0}+e_{\frak I+i_0,\frak I+i_0})$$ or
$$y=y-\tr(y)e_{2\frak I+1,2\frak I+1}
+\tr(y)e_{2\frak I+1,2\frak I+1}.$$ Hence, \eqref{asaface2null0}
holds. \qed

\medskip

\begin{corollary}\label{toaffinecase2null0}
Let $x\in T^+$ or $x\in T^-\setminus F\iota$. Then,
some
$0\neq h\in\frak h^{\pm }$
exists such that
 $ \mathcal B(x, h)\neq 0$.

 \end{corollary}
\proof By \eqref{asaface2null0} in Lemma \ref{otherdecnull0},
a finite subset $I\subset\mathfrak I$
and $0\neq h'\in\frak h_{2I}^{\pm}$ or $\frak h_{2I+1}^{\pm}$
exist such that $x=h'+b\iota_{2I}+x'$ or $x=h'+b\iota_{2I+1}+x'$
for some $b\in F$ and $x'\in T_{2\frak I\setminus 2I}$ or
$T_{2\frak I+1\setminus 2I+1}$ (since $x\notin F\iota$). Since the
trace form is nondegenerate on $\frak h^{\pm}_{2I}$ or $\frak
h^{\pm}_{2I+1}$, we can select $h\in\frak h^{\pm}_{2I}$ or $\frak
h^{\pm}_{2I+1}$ such that $\tr(h'h)\neq 0$. Then, we have
$\mathcal B(x, h) =\tr(h'h)+b\tr (h)+\tr(x'h)=\tr(h'h)\neq 0$
(since $x'h=0$). \qed

\enskip

By Corollary \ref{toaffinecase2null0} related to $T^+$, we can
also see that $\mathcal B^\sigma$ is nondegenerate. (We use the
result related to $T^-$ later.) Moreover, the restriction of
$\mathcal B^\sigma$ to any subalgebra $\mathcal L$ of $\mathcal
M_\frak K^\sigma$ that contains $\sll_\frak K(F)^\sigma$ is still
a nondegenerate form.

Conversely, let $U$ be a complement of $\frak h^\sigma$ in
$\mathcal L\cap T^\sigma$, and $\varphi$ is an arbitrary symmetric
bilinear form on $U$. Then, we can extend $\varphi$ to a
nondegenerate form on $\mathcal L$, using Lemma
\ref{otherdecnull0} (or embedding $\mathcal L$ into $\mathcal
M_\frak K$) and Corollary \ref{toaffinecase2null0} again.
Consequently, we can say that a LEALA of type $\text X_\frak
I\neq\text A_\frak I$ of nullity $0$ is isomorphic to a subalgebra
of $\mathcal M_\frak K^\sigma$ that contains $\sll_\frak
K(F)^\sigma$.

\enskip

\section{LALAs}
\label{sec:examplesL}

The next interesting objects are LEALAs of null dimension 1. In
fact, our aim in this study is to classify tame LEALAs of nullity
1.

\begin{definition}\label{defLALA}
A tame LEALA of nullity 1 is called a {\bf LALA}.
\end{definition}

We know that the core of a LALA is a locally Lie $G$-torus (see
\eqref{coreistorus}), and since $R^\times$ is a LEARS of nullity
$1$, the core is a locally Lie $1$-torus. Moreover, using the
notations in Section 4, we have the following.
\begin{lemma}\label{dercentless}
Let $\mathcal L$ be a LALA. Then:

(1) The core $\mathcal L_c$ is a universal covering of a locally
loop algebra.

(2) $Z(\mathcal L)=Z(\mathcal L_c)$, and a natural embedding
$\ad\mathcal L \hookrightarrow \der_F\big(\mathcal L_c/Z(\mathcal
L_c)\big)$ exists.

In particular, if $\mathcal N$ is a complement of the core
$\mathcal L_c$, i.e., $\mathcal L = \mathcal L_c \oplus \mathcal
N$, then $\mathcal N$ can be identified with a subspace of
$\oder_F\big(\mathcal L_c/Z(\mathcal L_c)\big)$, i.e., the outer
derivations of a locally loop algebra.
\end{lemma}
\proof By Proposition \ref{centerofleala}, $\mathcal L_c$ has a
nontrivial center. Hence, by Theorem \ref{Theorem 2.4},
 (1) is true.
 For (2),
we have $Z(\mathcal L_c)= Z(\mathcal L)$ by Proposition
\ref{centerofleala}. Since $\ad\mathcal L\cong \mathcal
L/Z(\mathcal L)$, we obtain the embedding using Lemma
\ref{derofcenterlesscore}. \qed

\enskip

To complete the classification of LALAs, we need to classify a
complement of the core. First, we give some examples of LALAs.

\example\label{exampleLALA} Let $\mathfrak I$ be an arbitrary
index set. We can construct 14 minimal standard LALAs (see
Defintion \ref{defstand} and \ref{ourmin}) from the 14 locally
loop algebras $L(\text X_{\mathfrak I}^{(i)})$ given in Section 3.
Thus,
$$\mathcal L^{ms}=\mathcal L^{ms}(\text X_{\mathfrak I}^{(i)}) :=L(\text X_{\mathfrak I}^{(i)})\oplus Fc\oplus Fd^{(0)}$$
is a LALA of type $\text X_{\mathfrak I}^{(i)}$, where $c$ is
central and $d^{(0)}$ is the degree derivation, i.e.,
$$d^{(0)}(t^m)=mt^m$$
with a Cartan subalgebra
$$\mathcal H =\mathfrak h\oplus Fc\oplus Fd^{(0)},$$
where $\mathfrak h$ is the subalgebra of $\mathfrak g(\text
X_{\mathfrak I})$, which comprises diagonal matrices if $\mathfrak
I$ is infinite or any Cartan subalgebra if $\mathfrak I$ is
finite. In addition, a nondegenerate invariant symmetric bilinear
form $\mathcal B$ on $\mathcal L ^{ms}$ is an extension of the
form defined in Section 3 for loop algebras using the trace form
or the Killing form if $\mathfrak I$ is finite, and a
nondegenerate symmetric associative bilinear form on $F[t^{\pm
1}]$, and by defining $\mathcal B(c,d^{(0)})=1$. In particular, we
define $\mathcal B(d^{(0)},d^{(0)})=0$ as usual, although
 $\mathcal B(d^{(0)},d^{(0)})$ can be any number in $F$.
These LALAs are minimal LALAs. Note that any standard LALA
contains a minimal standard LALA. In addition, we note that if
$\mathfrak I$ is finite, then LALAs are automatically minimal
standard LALAs, which are the affine (Kac-Moody) Lie algebras.
Note that a minimal standard LALA $\mathcal{L}^{ms}$ is also
denoted by $\mathcal{L}(0)$.

\enskip

Now, we give examples of bigger (and the biggest) LALAs when
$\mathfrak I$ is infinite. Note that
$$\sll_\frak I(F)+T=\gl_\frak I(F)+T,$$
where $T=T_\frak I$ is the subspace of all the diagonal matrices
in the matrix space $M_\frak I(F)$ of size $\frak I$, which is a
Lie algebra with the split center $F\iota$, where $\iota$ is the
diagonal matrix and its diagonal entries are all $1$. Thus, its
loop algebra
\begin{equation}\label{defofuloop}
\mathcal U=\mathcal U_\frak I:=\big(\sll_\frak I(F)+T\big)\otimes
F[t^{\pm 1}]
\end{equation}
is a Lie algebra with the split center $\iota\otimes F[t^{\pm
1}]$.

Assume that $\mathcal B$ is a symmetric invariant bilinear form on
$\mathcal U$, which is not a zero on $\sll_\frak I(F)$. Then, by
Lemma \ref{uniquenessform} and Lemma \ref{perfectidealpr},
$\mathcal B$ is unique up to a scalar to $\tr\otimes\epsilon$ on
\begin{equation}\label{halfandhalfdef}
\big(\sll_\frak I(F)\otimes F[t^{\pm 1}]\big)\times\mathcal U
\quad\text{and}\quad \mathcal U\times\big(\sll_\frak I(F)\otimes
F[t^{\pm 1}]\big),
\end{equation}
i.e., for $x,y\in\mathcal U$, and if $x$ or $y\in \sll_\frak
I(F)$, then
\begin{equation}\label{defBdown}
\mathcal B(x\otimes t^m, y\otimes t^n) =a\tr(xy)\delta_{n,-m}
\end{equation}
for some $a\in F^\times$. We claim that such a form $\mathcal B$
does exist. As in the case of nullity $0$, we select a complement
$\frak h^c$ of $\frak h$ in $T$, i.e., $T=\frak h^c\oplus \frak
h$. For each $m\in\Bbb Z$, let
$$\psi_m:\frak h^c\times \frak h^c\longrightarrow F
$$
be an arbitrary bilinear form. We define a symmetric bilinear form
$\mathcal B$ on $\mathcal U$ as
$$
\mathcal B(x\otimes t^m, y\otimes t^n) =\psi_m(x, y)\delta_{n,-m}
$$
on $\frak h^c\otimes F[t^{\pm 1}]$,
 and \eqref{defBdown}
 on \eqref{halfandhalfdef}.
We can prove that $\mathcal B$ is invariant in a similar manner to
the case of nullity $0$ using the following claim (which can also
be proved in a similar manner to Claim \ref{tonull0case}).

\begin{claim}\label{toaffinecase}
Let $x\in T\setminus F\iota$ and $y_k\in \sll_\frak I(F)$ for
$k=1,2,\ldots, r$. Then, a finite subset $I$ of $\mathfrak I$,\ \
$0\neq h\in\frak h$
exist and $g\in T$
 such that
$y_k\in\sll_I(F)$ for all $k$, and $h\in\frak h_I$,
 $$x=h+g,\quad
[x\otimes t^m, y_k\otimes t^n]=[h\otimes t^m, y_k\otimes t^n]
\quad\text{and}\quad \mathcal B(x\otimes t^m, y_k\otimes
t^n)=\mathcal B(h\otimes t^m, y_k\otimes t^n)$$
 for all $m, n\in\Bbb Z$ and all $k$.
 Moreover,
$y\in \sll_\frak I(F)$ and $h'\in\frak h$
exist such that
\begin{equation}\label{nonzeropair}
[x\otimes t^m, y\otimes t^n]\neq 0 \quad\text{and}\quad \mathcal
B(x\otimes t^m, h'\otimes t^{-m})\neq 0.
\end{equation}
\qed
 \end{claim}

Now, we can use a general construction, i.e., a one-dimensional
central extension by the $2$-cocycle
$$\varphi(u,v):=\mathcal B\big(d^{(0)}(u),v\big)$$
for $u,v\in\mathcal U$, where $d^{(0)}$ is the degree derivation
on $\mathcal U$. This is well known (e.g., see \cite{AABGP}), but
for convenience, we show that $\varphi$ is a $2$-cocycle in a
slightly more general setup. Note that $d^{(0)}$ is a skew
derivation relative to $\mathcal B$, i.e.,
$$\mathcal B\big(d^{(0)}(u),v\big)=-\mathcal B\big(u,d^{(0)}(v)\big).$$
More generally, for a $\Bbb Z$-graded algebra
$A=\bigoplus_{m\in\Bbb Z}\ A_m$ with a symmetric {\bf graded
bilinear form} $\psi$, the degree derivation $d^{(0)}$ is skew
relative to $\psi$. In fact, for $x=\sum_m x_m$ and $y=\sum_m
y_m\in A$, we have $\psi\big(d^{(0)}(x),y\big) =\sum_mm\psi(x_m,y)
=\sum_mm\psi(x_m,y_{-m}) =\sum_mm\psi(x,y_{-m})
=-\sum_mm\psi(x,y_{m}) =-\psi\big(x,d^{(0)}(y)\big)$. Hence,
$d^{(0)}$ is skew.

In general, on a Lie algebra $L$ with a symmetric invariant
bilinear form $B$, we can define $\varphi(u,v):=B\big(d(u),v\big)$
for any skew derivation $d$ and $u,v\in L$. Then, $\varphi(u,v)$
is a $2$-cocycle (which is also well known). In fact, the first
condition of the cocycle clearly holds, i.e., $\varphi(u,u)=0$ for
all $u\in L$, since $\varphi(u,u) = B(d(u),u) = -B(u,d(u)) =
-B(d(u),u) = -\varphi(u,u)$. For the second condition, we have
\begin{align*}
&\varphi([u,v],w)+\varphi([v,w],u)+\varphi([w,u],v)\\
&=B\big(d([u,v]),w)-B\big([v,w],d(u)\big)-B\big([w,u],d(v)\big)\\
&=B\big([d((u),v]),w)+B\big([u,d(v)]),w)-B\big([v,w],d(u)\big)-B\big([w,u],d(v)\big)\\
&=B\big(d(u),[v,w])-B\big(d(v),[u,w])-B\big([v,w],d(u)\big)-B\big([w,u],d(v)\big)=0.
\end{align*}
Thus, we obtain a $1$-dimensional central extension
$$
\tilde{\mathcal U}:=\mathcal U\oplus Fc
$$
using the $2$-cocyle $\varphi(u,v)=\mathcal
B\big(d^{(0)}(u),v\big)$ given above. Then,
$$
\hat{\mathcal U}=\hat{\mathcal U}_\frak I:=\tilde{\mathcal
U}\oplus Fd^{(0)}
$$
is naturally a Lie algebra that defines $$[c,d^{(0)}]=0,$$
anti-symmetrically. Thus, the center of $\hat{\mathcal U}$ is
equal to $Fc\oplus F\iota$. We also extend the form $\mathcal B$
by
$$\mathcal B(c,d^{(0)})=1
\quad\text{and}\quad \mathcal B(\mathcal U,d^{(0)})=0,
$$
symmetrically (where the value of $\mathcal B(d^{(0)},d^{(0)})$
can be any). Then, we can also check that this extended form is
invariant.

Let $\frak g:=\sll_\frak I(F)$ and let
$$\frak g=\frak h\oplus\bigoplus_{\mu\in\text A_\frak I\subset \frak h^*}\
\frak g_\mu$$ be the root-space decomposition of $\frak g$
relative to $\frak h$. Let
$${ \mathcal H} :=T\oplus Fc\oplus Fd^{(0)}.$$
We extend each root $\mu\in\frak h^*$ to an element in $\mathcal
H^*$ as follows. First, we can extend $\mu$ to $T$ as in the case
of nullity $0$. Then, we define $\mu(Fc\oplus Fd^{(0)})=0$. We
also define $\delta\in \mathcal H^*$ as $\delta(T\oplus Fc)=0$ and
$\delta(d^{(0)})=1$. Then, $\hat{\mathcal U}$ has the root space
decomposition
$$\hat{\mathcal U}
=\bigoplus_{\xi\in\mathcal H^*}\ \hat{\mathcal U}_{\xi}$$ relative
to $\mathcal H$, where $\hat{\mathcal U}_{\mu+m\delta}=\frak
g_\mu\otimes t^m$ for $\mu\in\text A_\frak I$, $\hat{\mathcal
U}_{m\delta}=T\otimes t^m$ for $m\neq 0$ and $\hat{\mathcal
U}_0=\mathcal H$, and $\hat{\mathcal U}_{\xi}=0$ if
$\xi\notin\text A_\frak I^{(1)} = (\text{A}_\frak I \cup \{ 0 \})
+ \mathbb{Z}\delta$. For convenience, we assume that $0 \not\in
\text{A}_\frak I$ but $0 \in \text{A}_\frak I^{(1)}$. In the
above, $\hat{\mathcal{U}}$ is an $\langle\text A_\frak
I^{(1)}\rangle$-graded Lie algebra, and $\mathcal B$ is graded in
the sense that $\mathcal B(\hat{\mathcal U}_\xi,\hat{\mathcal
U}_\eta)=0$, unless $\xi+\eta=0$ for all $\xi,\eta\in\text A_\frak
I^{(1)}$. In particular, the radical of $\mathcal B$ is graded.

\claim The radical of $\mathcal B$ is contained in $\iota\otimes
F[t^{\pm 1}]$.

\endclaim

\proof Since the radical of $\mathcal B$ is graded, we can check
the nondegeneracy for each homogeneous element. It is clear that
the elements of degree $\mu+m\delta$ for $\mu\in\text A_\frak I$
cannot be in the radical. The elements of degree $m\delta$ are
also outside of the radical by \eqref{nonzeropair}. Hence, the
radical should be in $\iota\otimes F[t^{\pm 1}]$. \qed

\enskip

Now, it is easy to check that $(\hat{\mathcal U},\mathcal
H,\mathcal B)$ is a LEALA of nullity $1$ by defining
$\psi_0(\iota,\iota)\neq 0$. Since the center of $\hat{\mathcal
U}$ is equal to $Fc\oplus F\iota$, this is not tame. However,
since $\iota\otimes F[t^{\pm 1}]$ is an ideal of $\hat{\mathcal
U}$, the quotient LEALA
$$
\mathcal L^{max}:=\hat{\mathcal U}/\big(\iota\otimes F[t^{\pm
1}]\big)
$$
is tame, by defining $\psi_0(\iota,\iota)= 0$. Thus, $\mathcal
L^{max}$ is a LALA, which is isomorphic to the Lie algebra
\eqref{maxreali} described in the Introduction. The core $\mathcal
L^{max}_c$ is equal to $\sll_\frak I(F)\otimes F[t^{\pm 1}]\oplus
Fc$. (As stated in the Introduction, we omit bars for the quotient
Lie algebra.) Moreover, it is easy to check that a $1$-dimensional
extension of the core, such as
$$\mathcal L(p)=\mathcal L^{max}_c\oplus F(d^{(0)}+p)$$
for some $p\in T$, is a minimal LALA of type $\text A_\frak
I^{(1)}$ (which is a subalgebra of $\mathcal L^{max}$). In
addition, we can show that any homogeneous subalgebra of $\mathcal
L^{max}$ that contains some $\mathcal L(p)$ is a LALA. In Section
6, we show that any LALA of type $\text A_\frak I^{(1)}$ is a
homogeneous subalgebra of $\mathcal L^{max}$ that contains some
$\mathcal L(p)$.

\medskip

We describe the other untwisted LALAs using $\hat{\mathcal
U}_{2\frak I}$ and $\hat{\mathcal U}_{2\frak I+1}$, and the
automorphism $\sigma$ is again defined in \eqref{originalinv}. Let
$$
T=T^\sigma\oplus T^-
$$
be the decomposition of $T=T_{2\frak I}$ or $T_{2\frak I+1}$,
where $T^\sigma$ is the eigenspace of eigenvalue $1$ (the fixed
algebra of $T$ by $\sigma$) and $T^-$ is the eigenspace of
eigenvalue $-1$. Instead of $T^\sigma$, we use $T^+$ because we
consider the fixed algebra by another automorphism $\tau$ later.
Thus, we have
\begin{align}\label{descriptionoftsigma}
&T^\sigma=\{(a_{kk})\in T_{2{\mathfrak I}+1}\mid
a_{ii}=-a_{{\mathfrak I}+i,{\mathfrak I}+i}\ (\forall
i\in\mathfrak I),\ a_{2{\mathfrak I}+1,2{\mathfrak I}+1}=0\}\quad
\text{and}
\nonumber \\
&T^-=\{(a_{kk})\in T_{2{\mathfrak I}+1}\mid a_{ii}=a_{{\mathfrak
I}+i,{\mathfrak I}+i}\ (\forall i\in\mathfrak I)\}\quad \text{for
$\text B_{\mathfrak I}^{(1)}$,}
\nonumber \\
&T^\sigma=\{(a_{kk})\in T_{2{\mathfrak I}}\mid
a_{ii}=-a_{{\mathfrak I}+i,{\mathfrak I}+i}\ (\forall
i\in\mathfrak I)\}\quad \text{and}
\nonumber \\
&T^-=\{(a_{kk})\in T_{2{\mathfrak I}}\mid a_{ii}=a_{{\mathfrak
I}+i,{\mathfrak I}+i}\ (\forall i\in\mathfrak I)\} \quad\text{for
$\text C_{\mathfrak I}^{(1)}$ or $\text D_{\mathfrak I}^{(1)}$.}
\end{align}
Note that $F\iota$ is $\sigma$-invariant and $F\iota\subset T^-$.
Since $\frak h$ is $\sigma$-invariant, we have
$$
\frak h=\frak h^\sigma\oplus \frak h^-, \quad \frak h^\sigma =
\mathfrak h \cap T^\sigma \subset T^\sigma \quad\text{and}\quad
\frak h^- = \mathfrak h \cap T^- \subset T^-.
$$
Let us extend the automorphism on $\hat{\mathcal U}=\hat{\mathcal
U}_{2\frak I}$ or $\hat{\mathcal U}_{2\frak I+1}$ as
$$
\hat\sigma(x\otimes t^k):=\sigma(x)\otimes t^k,\quad
\hat\sigma(c):=c\quad\text{and}\quad \hat\sigma(d^{(0)}):=d^{(0)}.
$$
Then, the fixed algebra $\hat{\mathcal U}^{\hat\sigma}$ with the
restriction of the form $\mathcal B$ is a LALA of type $\text
B_\frak I^{(1)}$, $\text C_\frak I^{(1)}$ or $\text D_\frak
I^{(1)}$, depending on the type of $\sigma$. In particular,
$$
\hat{\mathcal U}^{\hat\sigma} =\big( (\frak g+T^\sigma)\otimes
F[t^{\pm 1}] \big) \oplus Fc\oplus Fd^{(0)},
$$
where $\frak g=\sll_{2\frak I+1}(F)^\sigma$ or $\sll_{2\frak
I}(F)^\sigma$ is a locally finite split simple Lie algebra of each
type.

The nondegeneracy of the restricted form $\mathcal B$ follows from
the next lemma, where the proof is similar to the case in nullity
$0$.

\begin{lemma}\label{toaffinecase2}
Let $0\neq x\in T^\sigma$ or $x\in T^-\setminus F\iota$. Then,
some
$h\in\frak h^{\sigma }$ or $h\in\frak h^-$
exist such that
 $$
 \mathcal B(x\otimes t^m, h\otimes t^{-m})\neq 0$$
 for all $m\in\Bbb Z$.
 \qed
 \end{lemma}

As in the case of type $\text A_\frak I^{(1)}$, a $1$-dimensional
extension of the core $\hat{\mathcal U}^{\hat\sigma}_c$, such as
$$\mathcal L(p)=\hat{\mathcal U}^{\hat\sigma}_c\oplus F(d^{(0)}+p)$$
for some $p\in T^{ \sigma}$, is a minimal LALA of each type. In
addition, we can check that any homogeneous subalgebra of
$\mathcal{L}^{max} = \hat{\mathcal U}^{\hat\sigma}$ that contains
some $\mathcal L(p)$ is a LALA of each type. In Section 6, we show
that any LALA of each type is a homogeneous subalgebra of
$\mathcal{L}^{max} = \hat{\mathcal U}^{\hat\sigma}$ that contains
some $\mathcal L(p)$.

\enskip

We now give examples of twisted LALAs. 
Again, we use the
automorphism $\sigma$ defined in \eqref{originalinv} to obtain the
type $\text C_\frak I$ or $\text B_\frak I$, and we extend the
automorphism on $\hat{\mathcal U}=\hat{\mathcal U}_{2\frak I}$ or
$\hat{\mathcal U}_{2\frak I+1}$ as
\begin{equation}\label{deftwisting}
\hat\sigma(x\otimes t^k):=(-1)^k\sigma(x)\otimes t^k,\quad
\hat\sigma(c)=c\quad\text{and}\quad \hat\sigma(d^{(0)}):=d^{(0)}.
\end{equation}
Then, the fixed algebra $\hat{\mathcal U}^{\hat\sigma}$ with the
restriction of the form $\mathcal B$ is a LALA of type $\text
C_\frak I^{(2)}$ or $\text {BC}_\frak I^{(2)}$, depending on the
type of $\sigma$. In particular,
\begin{equation}\label{twbc-c}
\hat{\mathcal U}^{\hat\sigma}=\big( (\frak g^\sigma\oplus
T^\sigma)\otimes F[t^{\pm 2}] \big) \oplus \big( (\frak
g^-+T^-)\otimes tF[t^{\pm 2}] \big) \oplus Fc\oplus Fd^{(0)},
\end{equation}
where $\frak g^\sigma=\sll_{2\frak I}(F)^\sigma=\text{sp}_{2\frak
I}(F)$ or $\sll_{2\frak I+1}(F)^\sigma=\text o_{2\frak I+1}(F)$,
and $\frak g^-$ is the minus space of $\sll_{2\frak I}(F)$ or
$\sll_{2\frak I+1}(F)$ by $\sigma$. Since $\iota\otimes tF[t^{\pm
2}]$ is an ideal of $\hat{\mathcal U}$, the quotient LEALA
$$
\overline{\hat{\mathcal U}^{\hat\sigma}} :=\hat{\mathcal
U}^{\hat\sigma}/\big(\iota\otimes tF[t^{\pm 2}]\big)
$$
is tame, by defining $\psi_0(\iota,\iota)= 0$. Thus,
$\overline{\hat{\mathcal U}^{\hat\sigma}}$ is a LALA. (For the
type $\text C_\frak I^{(2)}$, this is isomorphic to the Lie
algebra \eqref{maxtwisted} described in the Introduction.)

The nondegeneracy of the restricted form $\mathcal B$ follows from
Lemma \ref{toaffinecase2}. As in the untwisted case, a
$1$-dimensional extension of the core $(\overline{\hat{\mathcal
U}^{\hat\sigma}})_c$, such as
$$\mathcal L(p):=(\overline{\hat{\mathcal U}^{\hat\sigma}})_c\oplus F(d^{(0)}+p)$$
for some $p\in T^\sigma$, is a minimal LALA of each type. We can
also show that any homogeneous subalgebra of $\mathcal{L}^{max} :=
\overline{\hat{\mathcal U}^{\hat\sigma}}$ that contains some
$\mathcal L(p)$ is a LALA of each type. In Section 7, we show that
any LALA of each type is a homogeneous subalgebra of
$\mathcal{L}^{max}$ that contains some $\mathcal L(p)$.

For the type $\text B_\frak I^{(2)}$, as described by Neeb in
\cite[App.1]{N2}, we introduce an automorphism $\tau$ on the
untwisted LALA $\mathcal M:=\hat{\mathcal U}^{\hat\sigma}_{2\frak
I+2}$ of type $\text D_{\frak I+1}^{(1)}$, which is defined by $s=
\begin{pmatrix}
   0   &\iota_{\frak I+1}       \\
  \iota_{\frak I+1}     & 0   \\
\end{pmatrix}
$. For convenience, let $\frak I+1=\{j\mid j\in\frak
I\}\cup\{j_0\}$ and
$$
2\frak I+2= (\frak I+1)+(\frak I+1) =\big(\{j\mid j\in\frak
I\}\cup\{j_0\}\big)\cup\big(\{-j\mid j\in\frak
I\}\cup\{-j_0\}\big).$$ Let
$$g=\iota_{2\frak I}+e_{j_0,-j_0}+e_{-j_0,j_0}$$
be the matrix of exchanging rows or columns, and let $\tau$ be an
involutive automorphism of $\text o_{2\frak I+2}(F)$ defined by
$$
\tau(x)=gxg.
$$
Then, we can see that the fixed algebra $\text o_{2\frak
I+2}(F)^\tau=\text o_{2\frak I+1}(F)$ (which is of type $\text
B_\frak I$) and the minus space
\begin{equation}\label{definefraks}
\frak s:=\{x\in\text o_{2\frak I+2}(F)\mid \tau(x)=-x\}
\end{equation}
by $\tau$ is isomorphic to $F^{2\frak I+1}$ as a natural $\text
o_{2\frak I+1}(F)$-module with
$$\frak s_0=\frak s\cap\frak h^\sigma
=F(e_{j_0j_0}-e_{-j_0,-j_0}).$$ We can extend $\tau$ on $\text
o_{2\frak I+2}(F)+T_{2\frak I+2}^\sigma$. Then, we have
\begin{equation}\label{plusandminusofTforD}
\big({T}^\sigma_{2\frak I+2}\big)^\tau = {T}^\sigma_{2\frak
I+1}(\cong {T}^\sigma_{2\frak I}) \quad\text{and}\quad \{x\in
T^\sigma_{2\frak I+2}\mid \tau(x)=-x\}=\frak s_0.
\end{equation}
We can also extend $\tau$ on $\mathcal M$ in the same 
manner as
\eqref{deftwisting}, i.e.,
\begin{equation*}\label{deftwistingtau}
\hat\tau(x\otimes t^k):=(-1)^k\tau(x)\otimes t^k,\quad
\hat\tau(c):=c\quad\text{and}\quad \hat\tau(d^{(0)}):=d^{(0)},
\end{equation*}
and we obtain a LALA $\mathcal M^{\hat\tau}$ of type $\text
B_\frak I^{(2)}$. In particular, we have
$$
\mathcal M^{\hat\tau}=\big( (\text o_{2\frak
I+1}(F)+{T}^\sigma_{2\frak I+1}) \otimes F[t^{\pm 2}] \big) \oplus
\big( \frak s\otimes tF[t^{\pm 2}] \big) \oplus Fc\oplus Fd^{(0)}.
$$
(The odd degree part of $t$ is the same as that in an affine Lie
algebra of type $\text B_\ell^{(2)}=\text D_{\ell+1}^{(2)}$.)

The nondegeneracy of the restricted form $\frak B$ follows from
Lemma \ref{toaffinecase2}. As in the above, a $1$-dimensional
extension of the core $\mathcal M^{\hat\tau}_c$, such as
$$\mathcal L(p)=\mathcal M^{\hat\tau}_c\oplus F(d^{(0)}+p)$$
for some $p\in {T}^\sigma_{2\frak I+1}$, is a minimal LALA of type
$\text B_\frak I^{(2)}$. In addition, we can show that any
homogeneous subalgebra of $\mathcal{L}^{max} := \mathcal
M^{\hat\tau}$ that contains some $\mathcal L(p)$ is a LALA of type
$\text B_\frak I^{(2)}$. In Section 7, we show that any LALA of
type $\text B_\frak I^{(2)}$ is a homogeneous subalgebra of
$\mathcal{L}^{max}$ that contains some $\mathcal L(p)$.

\endexample

\enskip

\section{Classification of the untwisted LALAs}
\label{sec:classificationunt}

Let $\mathcal L$ be an untwisted LALA of infinite rank, i.e., the
core $\mathcal L_c$ is a universal covering of an untwisted
locally loop algebra of type $\text A_\frak I^{(1)}$, $\text
B_\frak I^{(1)}$, $\text C_\frak I^{(1)}$ or $\text D_\frak
I^{(1)}$ for an infinite index $\frak I$. By selecting a
homogeneous complement of the $\Bbb Z$-graded core, we can write
$$\mathcal L=\mathcal L_c\oplus \bigoplus_{m\in\Bbb Z} D^m.$$
Note that the complement is assumed to be included in the null
space:
$$\bigoplus_{m\in\Bbb Z} D^m\subset
\bigoplus_{\delta\in R^0}\mathcal L_\delta =\bigoplus_{m\in\Bbb
Z}\mathcal L_{m\delta_1} \quad\text{and}\quad D^m\subset\mathcal
L_{m\delta_1},$$ where $\delta_1$ is a generator of $\langle
R^0\rangle_\Bbb Z=\langle S+S\rangle_\Bbb Z\cong\Bbb Z$ (see
\eqref{shortZ} and Lemma \ref{tames+s}).

Let
$$\mathcal L_c':=\mathcal L_c/Z(\mathcal L_c)$$
be the centerless core. Moreover, let $(\frak g,\frak h)$ be the
grading pair of the Lie $1$-torus $\mathcal L_c$ such that $\frak
h$ is the set of diagonal matrices of a locally finite split
simple Lie algebra $\frak g$:
$$\frak g
=\frak h\oplus\bigoplus_{\alpha\in\Delta}\frak g_\alpha =[\mathcal
L_c^0,\mathcal L_c^0] \subset \mathcal L_c =\bigoplus_{m\in\Bbb
Z}\ \mathcal L_c^m,$$ where
$$\mathcal L_c^m=\bigoplus_{\alpha\in\Delta\cup\{0\}}
\ (\mathcal L_c)_\alpha^m.$$ We identify the grading pair $(\frak
g,\frak h)$ of the Lie $1$-torus $\mathcal L_c'$ and $\mathcal
L_c$. Moreover, we identify $\mathcal L_c'$ with
$$L:=\frak g\otimes F[t^{\pm 1}].$$

Now, we classify the {\bf diagonal derivations} of an untwisted
locally loop algebra $L$ in general. Let
$$(\der_FL)^0_0
=\{d\in\der_FL\mid d(\frak g_\alpha\otimes t^m)\subset \frak
g_\alpha\otimes t^m \ \text{for all $\alpha\in\Delta$ and
$m\in\Bbb Z$}\}.$$ We refer to such an element as a {\bf diagonal
derivation of degree $0$}. Let $d \in (\det_F L)_0^0$. Note that
since $\frak g_0=\frak h=\sum_{\alpha\in\Delta} [\frak
g_\alpha,\frak g_{-\alpha}]$, then
we have
\begin{align*}
d(\frak h\otimes t^m) &=\sum_{\alpha\in\Delta} d([\frak
g_\alpha,\frak g_{-\alpha}]\otimes t^m) =\sum_{\alpha\in\Delta}
d([\frak g_\alpha\otimes t^m,\frak g_{-\alpha}\otimes 1])\\
&= \sum_{\alpha\in\Delta} \big([d(\frak g_\alpha\otimes t^m),\frak
g_{-\alpha}\otimes 1]
+[\frak g_\alpha\otimes t^m,d(\frak g_{-\alpha}\otimes 1)]\big)\\
&\subset \sum_{\alpha\in\Delta} [\frak g_\alpha\otimes t^m,\frak
g_{-\alpha}\otimes 1] =\frak h\otimes t^m.
\end{align*}
In addition, we note that $d\mid_\frak g$ is a diagonal derivation
of $\frak g$. Hence, by Neeb \cite{N1}, we obtain $d\mid_\frak
g=\ad p$ for a certain diagonal matrix $p$ of an infinite size. In
particular, we have $p\in P$, where
\begin{equation}\label{defP}
\text{ $P=T_{\frak I}$ for $\text A_{\frak I}$, and $T^+_{2\frak
I}$ or $T_{2\frak I+1}^+$ for the other types}
\end{equation}
as defined in Example \ref{exampleLALA}. Put
$$
d':=d-\ad p\in(\der_FL)^0_0.
$$
Then, we have
$$
d'(\frak g\otimes 1)=0.
$$
In particular, we have $d'(\frak h\otimes 1)=0$. Thus, for $0\neq
x\otimes t\in \frak g_\alpha\otimes t$, if
\begin{equation}\label{dxt}
d'(x\otimes t)=a x\otimes t
\end{equation}
for $a\in F$, then
\begin{equation}\label{dytinv}
d'(y\otimes  t^{-1})=-a y\otimes t^{-1}
\end{equation}
for all $y\in\frak g_{-\alpha}$. In fact, for $y \not= 0$, since
$0\neq [x,y]=h\in\frak h$ and $d'(y\otimes  t^{-1})=b y\otimes
t^{-1}$ for some $b\in F$, then we have
\begin{align*}
0&=d'(h\otimes 1)=d'([x\otimes t,y\otimes t^{-1}])
=[d'(x\otimes t),y\otimes t^{-1}]+[x\otimes 1,d'(y\otimes t^{-1})]\\
&=(a+b)[x\otimes t,y\otimes t^{-1}] =(a+b)[x,y]\otimes 1 = (a+b)h
\otimes 1.
\end{align*}
Hence, $b=-a$.

\begin{lemma}\label{simg}
Let $\frak g =\frak h\oplus\bigoplus_{\alpha\in\Delta}\frak
g_\alpha$ be a locally finite split simple Lie algebra. Then,
$\mathcal U(\frak g).\frak g_{\beta}=\frak g$ for any
$\beta\in\Delta$, where $\mathcal U(\frak g)$ is the universal
enveloping algebra of $\frak g$.
\end{lemma}
\proof Since $\mathcal U(\frak g).\frak g_{\beta}$ is a nonzero
ideal of $\frak g$, it must be equal to $\frak g$ by simplicity.
\qed

\enskip

By Lemma \ref{simg}, for a fixed $\alpha \in \Delta$, three
subspaces
$$
\text{$\frak g\otimes 1$,\quad  $\frak g_\alpha\otimes t$,\quad
and \quad $\frak g_{-\alpha}\otimes t^{-1}$}
$$
generate $L$ as a Lie algebra.

Let
$$
d'':=d'-a d^{(0)},
$$
where $d^{(0)}=\displaystyle{t\frac{d}{dt}}$. Then, we have
$d''(\frak g\otimes 1)=d'(\frak g\otimes 1)=0$ and using
\eqref{dxt},
$$d''(x\otimes t)=
d'(x\otimes t)-ax\otimes t=0$$ for $x \in \frak g_\alpha$.
Similarly, using \eqref{dytinv},
$$d''(y\otimes t^{-1})=
d'(y\otimes t^{-1})+ay\otimes t^{-1}=0$$ for $y \in \frak
g_{-\alpha}$. Thus, we have $d''(L)=0$ and $d''=0$. Hence, we
obtain
\begin{equation}\label{00deriv}
d=\ad p+a d^{(0)}, \quad\text{and} \quad (\der_FL)^0_0=\ad P\oplus
Fd^{(0)}.
\end{equation}

We define the {\bf shift map} $s_m$ for $m\in\Bbb Z$ on $L=\frak
g\otimes F[t^{\pm 1}]$ by
$$
s_m(x\otimes t^k):=x\otimes t^{k+m}
$$
for all $k\in\Bbb Z$. (Shift maps were discussed in the
classification of affine Lie algebras by Moody in \cite{Mo}.)
Clearly, the shift maps have the property
$$
s_m([x,y])=[s_m(x),y]=[x,s_m(y)]
$$
for $x, y\in L$. (In other words, the shift maps are in the
centroid of $L$.) Thus, $s_m\circ d$ is a derivation for any
derivation $d$ of $L$. In fact, for $x, y\in L$,
$$s_m\circ d([x,y])=s_m([d(x),y]+[x,d(y)])
=[s_m\circ d(x),y]+[x,s_m\circ d(y)].$$

Now, let
$$d\in(\der_FL)^m_0
=\{d\in\der_FL\mid d(\frak g_\alpha\otimes t^k)\subset \frak
g_\alpha\otimes t^{k+m} \ \text{for all $\alpha\in\Delta$ and
$k\in\Bbb Z$}\}.$$ Then, we have
$$s_{-m}\circ  d\in(\der_F L)_0^0.$$
Hence, by \eqref{00deriv}, $p=p_{d}\in P$ and some $a=a_{d}\in F$
exist such that
$$s_{-m}\circ d=\ad p+a  d^{(0)},$$
and thus
$$ d=s_m\circ (\ad p+a  d^{(0)}).$$
Therefore, we have classified diagonal derivations of the
untwisted locally loop algebra. Thus:
\begin{theorem} \label{diagonalderunt}
For all $m\in\Bbb Z$, we have
$$(\der_FL)^m_0
=s_m\circ(\der_FL)^0_0 =s_m\circ(\ad P\oplus Fd^{(0)}),$$ where
$P$ is defined in \eqref{defP}. \qed
\end{theorem}

The following property of diagonal derivations is useful later.
\begin{lemma} \label{classifydegreeder}
For all $m\in\Bbb Z$, let
$$(\der'_FL)^m_0
:=\{d\in (\der_FL)^m_0 \mid s_n\circ d=d\circ s_n \ \text{for some
$0\neq n\in\Bbb Z$}\}$$ and
$$(\der''_FL)^m_0
:=\{d\in (\der_FL)^m_0 \mid s_n\circ d=d\circ s_n \ \text{for all
$n\in\Bbb Z$}\}.$$ Then, we have
$$(\der'_FL)^m_0
=s_m\circ\ad  P =(\der''_FL)^m_0.
$$
\end{lemma}
\proof First, it is clear that
$$(\der'_FL)^m_0
\supset (\der''_FL)^m_0 \supset s_m\circ\ad P$$ for all $m\in\Bbb
Z$. Thus, it is sufficient to show that
\begin{equation}\label{otherinc}
(\der'_FL)^m_0 \subset s_m\circ\ad P.
\end{equation}
Therefore, let $s_m\circ(\ad p+a d^{(0)})\in(\der'_FL)^m_0 \subset
(\der_F L)^m_0$. Then, for
$$
h\otimes t^k\in \frak h\otimes t^k\subset L,
$$
we have
$$
s_n\circ s_m([p+a d^{(0)},h\otimes t^k]) =s_n(akh\otimes t^{k+m})
=akh\otimes t^{k+m+n}
$$
and
$$
[s_m\circ (p+a d^{(0)}),h\otimes t^{k+n}] =a(k+n)h\otimes
t^{k+n+m}
$$
for some $n\neq 0$. Hence, $an=0$, and we obtain $a=0$. Therefore,
we obtain
$$s_m\circ(\ad p+a d^{(0)})
=s_m\circ \ad p\in s_m\circ\ad P.$$ Thus, we have shown
\eqref{otherinc}. \qed

\enskip

\remark\label{azamresult} We can use some results given by Azam
related to the derivations of tensor algebras (see \cite[Thm
2.8]{A2}). However, their direct application to our tensor algebra
$\frak g\otimes_F  F[t^{\pm 1}]$ yields an isomorphism such that
$$
\der_F(\frak g\otimes_F F[t^{\pm 1}])\cong \der_F\frak
g\overleftarrow\otimes_F F[t^{\pm 1}] \ \oplus\ C(\frak
g)\overrightarrow\otimes_F \der_FF[t^{\pm 1}],
$$
where $C(\frak g)$ is the centroid of $\frak g$ and
$\overleftarrow\otimes_F$ and $\overrightarrow\otimes_F$ are
special types of tensor products (since $\frak g$ is
infinite-dimensional). Thus, we need to perform some more work to
obtain our desired form as given above. We only need a special
type of subspace, i.e., $(\der_FL)^m_0$, so we can approach them
directly without using Azam's result. In addition, we investigate
derivations of twisted locally loop algebras later that are not
tensor algebras.

\endremark

Now, we return to classifying $D^m$. Let $d\in D^m$. Then, $\ad
d\in(\der_F L)^m_0$ by Lemma \ref{dercentless}. Hence, by Theorem
\ref{diagonalderunt},
$p=p_{d}\in P$ (see \eqref{defP}) and some
$a=a_{d}\in F$ exist such that
$$\ad d=s_m\circ (\ad p+a  d^{(0)}).$$
We claim that $a=0$ for all $m\neq 0$. First, we note that
$h,\ h'\in\frak h$ exist such that $\tr(hh')\neq 0$.
In addition,
we have
$$\mathcal B(h\otimes t,h'\otimes t^{-1})
=\mathcal B(h\otimes t^m,h'\otimes t^{-m}) =c\tr(hh')\neq 0
$$
for all $m\in\Bbb Z$ and some $0\neq c\in F$ since $\mathcal
B=c\tr\otimes\epsilon$ (see Lemma \ref{uniquenessform} and we note
that $\epsilon$ is defined in the previous paragraph). Now, using
a pair $h$ and $h'$, we have
\begin{align*}
\mathcal B([d, h\otimes t],h'\otimes t^{-m-1})
&=a\mathcal B(h\otimes t^{m+1},h'\otimes t^{-m-1})\\
&=a\mathcal B(h\otimes t,h'\otimes t^{-1})
\end{align*}
while
\begin{align*}
\mathcal B([d, h\otimes t],h'\otimes t^{-m-1})
&=-\mathcal B(h\otimes t,[d,h'\otimes t^{-m-1}])\\
&=a(m+1)\mathcal B(h\otimes t,h'\otimes t^{-1}).
\end{align*}
Hence, $a=a(m+1)$, i.e., $am=0$. Thus, $m\neq 0$ implies that
$a=0$.

Moreover, suppose that $a=a_{d}=0$ for all $d\in D^0$. Then, $\ad
D^0\subset\ad P$
 (see \eqref{defP})
 and for the Cartan subalgebra $\mathcal H$ of the original LALA of
$\mathcal L$, we have $\mathcal H=\frak h\oplus Fc\oplus D^0$.
However, this contradicts the axiom $\mathcal L_0=\mathcal H$
since $\big[\frak h\otimes F[t^{\pm 1}],\mathcal H\big]=0$. Hence,
$p\in P$ exists such that $\ad p+d^{(0)}\in \ad D^0$.

Consequently, we obtain
$$\ad D^m\subset s_m\circ\ad P$$
for $m\neq 0$, and
$$\ad p+d^{(0)}\in\ad D^0\subset \ad P+Fd^{(0)}$$
for some $p\in P$. \remark In some cases, $d^{(0)}\notin \ad D^0$.
Thus, a LALA is not always standard. We can easily construct a
non-standard LALA even if $\dim_FD^0\geq 2$.
\endremark

Finally, we investigate the bracket on $D:=\bigoplus_{m\in\Bbb
Z}D^m$. Let $D':=\bigoplus_{m\neq 0}D^m$. First, note that
$[D',D']$ acts trivially on $L$ since $[\ad (p\otimes t^m),\ad
(p'\otimes t^n)] =\ad [p\otimes t^m,p'\otimes t^n]=0$ in $\der_F
L$. Hence,
$$[D',D']\subset Fc=Ft_{\delta_1}\subset \mathcal H,$$
by tameness. In addition, for $d_m\in D^m$ ($m\neq 0$) and $d_n\in
D^n$ ($n\neq 0$),
by the fundamental property \eqref{fpat} of a LEALA (see Lemma
\ref{FPAT}), we have,
$$
[d^m,d^n]=\delta_{m,-n}\mathcal B(d_m,d_n)t_{m\delta_1}
=m\delta_{m,-n}\mathcal B(d_m,d_n)t_{\delta_1}.
$$
Note that $\mathcal B(d_m,d_n)$ can be zero since $h\in\frak h$
exists such that $\tr(d_mh)\neq 0$ (and thus $\mathcal
B(d_m,h)\neq 0$).

Next, since $D^0\subset \mathcal H$, we have $[D^0,D^0]=0$.
Moreover, for $d\in D^0$ such that $\ad_L d=\ad_L p\in D^0$, we
have $[d,D^m]=0$. For the last case, i.e., for $d\in D^0$ such
that $\ad_L d=\ad_L p+a d^{(0)}\in \ad D^0$ and $d_m\in D^m$, we
have
$$[d,d_m]=[a  d^{(0)},d_m]=amd_m.$$

Now, $\iota\otimes t^m$ centralizes the $\mathcal L_c$, and hence
we obtain the following identifiication:
$$\mathcal L\cong
\big(\mathcal L+\iota\otimes F[t^{\pm 1}]\big) /\iota\otimes
F[t^{\pm 1}].$$ Thus:

\begin{theorem}\label{classiuntw}
Let $\mathcal L$ be an untwisted LALA. Then, $\mathcal L$ is
isomorphic to that in Example \ref{exampleLALA}. \qed
\end{theorem}

\enskip

\section{Classification of the twisted LALAs}
\label{sec:classificationt}

As mentioned earlier, each twisted loop algebra $M$ is a
subalgebra of an untwisted loop algebra $\tilde M$. In particular,
we have
\begin{align*}
\text{$M$ has type $\text B_\frak I^{(2)}$} &\ \Longrightarrow\
\text{$\tilde M$ has type $\text D_{\frak I+1}^{(1)}$}\\
\text{$M$ has type $\text C_\frak I^{(2)}$} &\ \Longrightarrow\
\text{$\tilde M$ has type $\text A_{2\frak I}^{(1)}$}\\
\text{$M$ has type $\text {BC}_\frak I^{(2)}$} &\ \Longrightarrow\
\text{$\tilde M$ has type $\text A_{2\frak I+1}^{(1)}$}.
\end{align*}

\smallskip

\remark\label{atypelabel} In the case where $\frak I$ is finite,
such as $\frak I = \{ 1,2, ... , n \}$, the type $\text A_\frak I$
usually means the Lie algebra $\sll_{\frak I+1}(F)$. Therefore, it
may be better to write
\begin{align*}
\text{$L$ has type $\text C_\frak I^{(2)} = \text C_n^{(2)}$} &\
\Longrightarrow\
\text{$\tilde L$ has type $\text A_{2\frak I-1}^{(1)} = \text A_{2n-1}^{(1)}$}\\
\text{$L$ has type $\text {BC}_\frak I^{(2)} = \text{BC}_n^{(2)}$}
&\ \Longrightarrow\ \text{$\tilde L$ has type $\text A_{2\frak
I}^{(1)} = \text A_{2n}^{(1)}$}
\end{align*}
in order to follow the common notations. However, in this study,
we use the type of the Lie algebra $\sll_{\frak I}(F)$ as $\text
A_\frak I$, instead of $\text A_{\frak I+1}$ provided that $\frak
I$ is an infinite set, as mentioned in the Introduction.
\endremark

\smallskip

First, we provide some basic lemmas for twisted locally loop
algebras, as follows:
$$
\begin{array}{l}
\text{(1)}\quad \frak ( g^\sigma\otimes F[t^{\pm 2}]) \oplus (
\frak g^-\otimes tF[t^{\pm 2}] ) \quad\text{for type ${\text
C}_{\frak I}^{(2)}$ or ${\text BC}_{\frak I}^{(2)}$,\quad and}\\
[0.2cm] \text{(2)}\quad \text ( o_{2\frak I+2}(F)^\tau\otimes
F[t^{\pm 2}] ) \oplus ( \frak s\otimes tF[t^{\pm 2}] )
\quad\text{for type ${\text B}_{\frak I}^{(2)}$}
\end{array}
$$
where in (1) $\frak g=\sll_{2\frak I+1}(F)$ or $\frak
g=\sll_{2\frak I}(F)$ (see \eqref{twbc-c} for $\sigma$), and in
(2) $\frak s$ is the minus space of $\text o_{2\frak I+2}(F)$ by
$\tau$, as described in \eqref{definefraks}. Note that $\text
o_{2\frak I+2}(F)^\tau=\text o_{2\frak I+1}(F)$, which has type
${\text B}_{\frak I}$.

\medskip

\lemma\label{irredg-} {\em (1)}\enskip $\frak g^-$ is an
irreducible $\frak g^\sigma$-module.

{\em (2)\enskip $\frak s$ is an irreducible $\text o_{2\frak
I+2}(F)^\tau$-module.

\endlemma

\proof For (1), it is sufficient to show that $w\in\mathcal
U(\frak g^\sigma)v$ for any $v,w\in\frak g^-$, where $\mathcal
U(\frak g^\sigma)$ is the universal enveloping algebra of $\frak
g^\sigma$. However, this is a local property. Thus, a
finite-dimensional split simple subalgebra $\frak f$ of $\frak g$
exists of the same type such that
 $v,w\in\frak f^-\subset \frak g^-$ and
 $\frak f^\sigma\subset \frak g^\sigma$.
 It is well known that this property holds in the finite-dimensional case
 (e.g., see \cite{K}).
 Thus, we are finished.
 Similarly,
 (2) holds.
\qed

\medskip

\lemma\label{centg+} {\em (1)}\enskip Let $C$ be the centralizer
of $\frak g^\sigma$ in $\frak g+T$. If $0\neq x\in C$, then $x\in
T^-\setminus\frak g^-$.

{\em (2)}\enskip Let $C$ be the centralizer of $\text o_{2\frak
I+2}(F)^\tau$ in $\sll_{2\frak I+2}(F)^\sigma+T_{2\frak
I+2}^\sigma =\text o_{2\frak I+2}(F)+T_{2\frak I+2}^\sigma$. Then,
$C=0$.
\endlemma

\proof For (1), we can write each Lie algebra as
$$\frak g + T = (\frak g+T)^\sigma\oplus (\frak g+T)^- = (\frak g^\sigma + T^\sigma) \oplus (\frak g^- + T^-).$$
Let
$$x=x_+\oplus x_-\in
(\frak g+T)^\sigma\oplus (\frak g+T)^- = (\frak g^\sigma +
T^\sigma) \oplus (\frak g^- + T^-)$$ be in $C$. Then, for any
$y\in\frak g^\sigma$, we have
$$0=[x,y]=[x_+,y]+[x_-,y].$$
Hence, $[x_+,y]=0$ and $[x_-,y]=0$. However, the centralizer
$C_{\frak g^\sigma+T^\sigma}\ \!\!(\frak g^\sigma)=0$ since $\frak
g^\sigma+T^\sigma$ is tame, as well as $Z(\frak g^\sigma) = 0$
(given that $\mathcal{L} = \mathcal{L}_c + T^\sigma$ is tame and
$\mathcal{L}_c = \frak ( g^\sigma\otimes F[t^{\pm 2}]) \oplus (
\frak g^-\otimes tF[t^{\pm 2}] ) \oplus Fc$,\ cf.~Section 4).
Hence, $x_+=0$, and we obtain $x\in (\frak g+T)^-=\frak g^-+T^-$.
If $x\in\frak g^-$, then $\mathcal U(\frak g^\sigma)x$, where
$\mathcal U(\frak g^\sigma)$ is the universal enveloping algebra
of $\frak g^\sigma$, is a $\frak g^\sigma$-submodule of $\frak
g^-$. However, since $\dim_F \frak g^->1$ and $\frak g^-$ is an
irreducible $\frak g^\sigma$-module (by Lemma \ref{irredg-}), then
we have $\frak g^- \supset \mathcal U(\frak g^\sigma)x = Fx$,
which implies that $x$ has to be $0$ in this case. Similarly, (2)
holds by property \eqref{plusandminusofTforD}. \qed

\lemma\label{-eleofh} {\em (1)}\enskip Let  $h\in T_{2\frak I+1}$
 for type ${\text BC}_{\frak I}^{(2)}$
 or
   $h\in T_{2\frak I}$
    for type ${\text C}_{\frak I}^{(2)}$.
   Suppose that $[h,\frak g^\sigma]\subset\frak g^-$,
   then
$h\in T_{2\frak I+1}^-$ or $h\in T_{2\frak I}^-$, respectively.

{\em (2)}\enskip Let $h\in T_{2\frak I+2}^\sigma$ for type $\text
B_{\frak I}^{(2)}$.
   Suppose that $[h,\text o_{2\frak I+2}(F)]\subset\frak s$,
   then $h\in \frak s_0 = \mathfrak s \cap \frak h^\sigma$
(for $\frak s_0$, see (44) in the last paragraph of Section 5).

\endlemma
\proof For (1), let $x\in\frak g^\sigma$ and $y=[h,x]\in\frak
g^-$. Then, $-y=\sigma(y) = [\sigma(h),x]$. Hence,
$[h+\sigma(h),x]=0$ for all $x\in\frak g^\sigma$. Therefore,
$h+\sigma(h)\in C$ in Lemma \ref{centg+}, and thus $h+\sigma(h)\in
T^-$. However since $h+\sigma(h)\in T^\sigma$, we obtain
$h+\sigma(h)=0$. Thus, $\sigma(h)=-h$, i.e., $h\in T_{2\frak
I+1}^-$ or $T_{2\frak I}^-$, respectively.

Similarly for (2), we obtain $h+\tau(h)=0$ by Lemma \ref{centg+}.
Thus, $h\in \frak s_0$.
    \qed

\enskip

Let $\mathcal L$ be a twisted LALA of infinite rank, i.e., the
core $\mathcal L_c$ is a universal covering of a twisted locally
loop algebra of type $\text B_\frak I^{(2)}$, $\text C_\frak
I^{(2)}$ or $\text {BC}_\frak I^{(2)}$ for an infinite index
$\frak I$. As in the untwisted case, by selecting a homogeneous
complement of the $\Bbb Z$-graded core, we can write
$$\mathcal L=\mathcal L_c\oplus \bigoplus_{m\in\Bbb Z} D^m,
\quad \bigoplus_{m\in\Bbb Z} D^m\subset \bigoplus_{\delta\in
R^0}\mathcal L_\delta =\bigoplus_{m\in\Bbb Z}\mathcal
L_{m\delta_1} \quad\text{and}\quad D^m\subset \mathcal
L_{m\delta_1},
$$
where $\delta_1$ is a generator of $\langle R^0\rangle_\Bbb Z$.
Let
$$\mathcal L_c':=\mathcal L_c/Z(\mathcal L_c)$$
be the centerless core and let $(\frak g,\frak h)$ be the grading
pair of the Lie $1$-torus $\mathcal L_c$ such that $\frak h$ is
the set of diagonal matrices of a locally finite split simple Lie
algebra $\frak g$, as before. According to this terminology,
$\mathcal L_c'=\mathcal L_c/Z(\mathcal L_c)$ can be identified
with
$$L:=\big(\frak g\otimes  F[t^{\pm 2}]\big)
\oplus \big(\frak s\otimes  tF[t^{\pm 2}]\big).$$

\smallskip

We note that the subalgebras $\frak g^+ = \frak g^\sigma$ and
$\frak g^-$ in the previous terminology correspond to $\frak g$
and $\frak s$ in this new terminology.

\medskip

Let $L$ be a locally loop algebra of type $X_\frak I^{(2)}$. Then,
$L$ is $\Delta$-graded, where $\Delta$ is a locally finite
irreducible root system of type $X_\frak I$. In addition, we can
see that $\frak g$ is $\Delta^{\rm red} \cup \{ 0 \}$-graded and
$\frak s$ is $\Delta' \cup \{ 0 \}$-graded, where $\Delta^{\rm
red}$ and $\Delta'$ are given as follows.

\bigskip

$$\begin{array}{c||c|c}
\Delta & \Delta^{\red} & \Delta'\\ \hline \hline \text B_\frak I &
\Delta & \Delta_{\sh}\ = (\text A_1)^{\times \frak I}\\ \hline
\text C_\frak I & \Delta & \Delta_{\sh}\ = \text D_\frak I\\
\hline \text{BC}_\frak I & \Delta_{\sh} \cup \Delta_{\lg}\ = \text
B_\frak I & \Delta\\ \hline
\end{array}$$

\bigskip

\noindent In this case, our new notation $(\text A_1)^{\times
\frak I}$ denotes the (orthogonal disjoint) union
$$(\text A_1)^{\times \frak I} = \dot{\sqcup}_{i \in \frak I} \Delta_i = \{ \pm \epsilon_i \mid i \in \frak I \}$$
of root systems $\Delta_i = \{ \pm \epsilon_i \}$ of type $\text
A_1$, which satisfy $\Delta_i \perp \Delta_j$ for distinct $i,j$.
In particular, we have
$$
\frak g = \bigoplus_{\alpha\in\Delta^{\rm red}\cup\{0\}}\ \frak
g_\alpha \quad \text{and}\quad \frak
s=\bigoplus_{\beta\in\Delta'\cup\{0\}}\ \frak s_\beta .
$$
As in the untwisted case, we can classify diagonal derivations of
a twisted locally loop algebra $L$ in general.

\medskip

Let
$$\begin{array}{l}
(\der_FL)^0_0 := \{d\in\der_FL\mid d(\frak g_\alpha\otimes
t^{2m})\subset \frak g_\alpha\otimes t^{2m}\\ [0.2cm] \qquad
\qquad \qquad \text{and}\ d(\frak s_\beta\otimes t^{2m+1})\subset
\frak s_\beta\otimes t^{2m+1} \ \text{for all
$\alpha\in\Delta^{\rm red}$, $\beta\in\Delta'$
 and $m\in\Bbb Z$}\}
 \end{array}
 $$
 and take $d\in (\der_FL)^0_0$.
Then, as before, $d\mid_\frak g$ is a diagonal derivation of
$\frak g$, and thus, by Neeb  \cite{N1}, $d\mid_\frak g=\ad p$ for
some $p\in P$, depending on the type of $\frak g$ (see
\eqref{defP}). Let
$$
d':=d-\ad p\in(\der_FL)^0_0.
$$
Then, we have $ d'(\frak g\otimes 1)=0 $. In particular, we have
$d'(\frak h\otimes 1)=0$. Thus, in the same manner as the
untwisted case, we can show that for $0\neq x\otimes t\in \frak
s_\beta\otimes t$, if
\begin{equation}\label{dxtt}
d'(x\otimes t)=a x\otimes t
\end{equation}
for $a\in F$, then
\begin{equation}\label{dytinvt}
d'(y\otimes  t^{-1})=-a y\otimes t^{-1}
\end{equation}
for all $y\in \frak s_{-\beta}$.

\begin{lemma}\label{simgt}
For the above $\frak s$, we have $\mathcal U(\frak g).\frak
s_{\beta}=\frak s$ for any $\beta\in\Delta'$, where $\mathcal
U(\frak g)$ is the universal enveloping algebra of $\frak g$.
\end{lemma}
\proof Since $\mathcal U(\frak g).\frak s_{\beta}$ is a nonzero
submodule of $\frak s$, then it must be $\frak s$ by the
irreducibility of $\frak s$. \qed

\enskip

By Lemma \ref{simgt}, for a fixed $\beta \in \Delta'$, the three
subspaces
$$
\text{$\frak g\otimes 1$,\quad  $\frak s_\beta\otimes t$,\quad and
\quad $\frak s_{-\beta}\otimes t^{-1}$}
$$
generate $L$ as a Lie algebra. As before, let $d'':=d'-a d^{(0)}$.
Then, we have $d''(\frak g\otimes 1)=d'(\frak g\otimes 1)=0$ and
using \eqref{dxtt}, we obtain $d''(x\otimes t)= d'(x\otimes
t)-ax\otimes t=0$ for $x \in \frak s_\beta$. Similarly, using
\eqref{dytinvt}, we have $d''(y\otimes t^{-1})= d'(y\otimes
t^{-1})+ay\otimes t^{-1}=0$ for $y \in \frak s_{-\beta}$. Thus, we
have $d''(L)=0$ and $d''=0$. Hence, we obtain
\begin{equation}\label{00derivt}
d=\ad p+a d^{(0)}.
\end{equation}

\medskip

Again, we define the shift map $s_{2m}$ for $m\in\Bbb Z$ on
$L=\big(\frak g\otimes  F[t^{\pm 2}]\big) \oplus \big(\frak
s\otimes  tF[t^{\pm 2}]\big)$ by
$$
s_{2m}(x\otimes t^{2k}):=x\otimes t^{2k+2m} \ \text{and}\ \
s_{2m}(v\otimes t^{2k+1}):=v\otimes t^{2k+2m+1}
$$
for $x\in\frak g$ and $v\in \frak s$. Let
$$\begin{array}{l}
(\der_FL)_0^{2m} :=\{d\in\der_FL\mid d(\frak g_\alpha\otimes
t^{2k})\subset \frak g_\alpha\otimes t^{2k+2m}\\ [0.2cm] \qquad
\qquad \qquad \text{and}\ d(\frak s_\beta\otimes t^{2k+1}) \subset
\frak s_\beta\otimes t^{2k+2m+1} \ \text{for all
$\alpha\in\Delta^{\rm red}$, $\beta\in\Delta'$
 and $k\in\Bbb Z$}\}
 \end{array}
 $$
 and we take $d_{2m} \in (\der_FL)_0^{2m}$.
Then, we have $s_{-2m}\circ d_{2m}\in(\der_FL)_0^0$. Hence, by
\eqref{00derivt}, some $p=p_{d_{2m}}\in P$ and $a=a_{d_{2m}}\in F$
exist such that $s_{-2m}\circ d_{2m}=\ad p+a d^{(0)}$, and thus
$$d_{2m}=s_{2m}\circ\ad p+at^{2m}d^{(0)}=s_{2m}\circ\ad p+at^{2m+1}\frac{d}{dt}.$$
Therefore, as in Theorem \ref{diagonalderunt}, we have the
following.

\begin{lemma} \label{diagonaldert}
For all $m\in\Bbb Z$, we have
$$(\der_F L)^{2m}_0
=s_{2m}\circ(\der_F L)^0_0 =s_{2m}\circ(\ad P\oplus Fd^{(0)}),$$
where $P$ is defined in \eqref{defP}. \qed
\end{lemma}

Moreover, as in Lemma \ref{classifydegreeder}, we have the
following.
\begin{lemma} \label{classifydegreedert}
For all $m\in\Bbb Z$, let
$$(\der'_F L)^{2m}_0
:=\{d\in (\der_F L)^{2m}_0 \mid s_{2n}\circ d=d\circ s_{2n} \
\text{for some $0\neq n\in\Bbb Z$}\}$$ and
$$(\der''_F L)^{2m}_0
:=\{d\in (\der_F L)^{2m}_0 \mid s_{2n}\circ d=d\circ s_{2n} \
\text{for all $n\in\Bbb Z$}\}.$$ Then, we have
$$(\der'_F L)^{2m}_0
=s_{2m}\circ\ad P =(\der''_F L)^{2m}_0.
$$
\qed
\end{lemma}

Now, we return to the classification of $D^m$. Let $d_{2m}\in
D^{2m}$. Then, $\ad d_{2m} =s_{2m}\circ\ad p+at^{2m}d^{(0)}$ for
some $p\in P$ and $a\in F$ by Lemma \ref{diagonaldert}. Then, as
in the untwisted case, we can show that $a=0$ for all $m\neq 0$,
using
$$\mathcal B([d, h\otimes t],h'\otimes t^{-m-1})
=-\mathcal B(h\otimes t^{2},[d_{2m},h'\otimes t^{-2m-2}])
$$
for some $h,h'\in\frak h$ such that $\tr(h,h')\neq 0$.
Furthermore, as in the untwisted case, some $p\in P$ exists such
that $\ad p+d^{(0)}\in\ad D^0$. Thus, the spaces $D^m$ for even
$m$ s coincide with those in Example \ref{exampleLALA}.

\enskip

Next, we determine $(\der_F L)_0^{2m+1}$, where
$$\begin{array}{l}
(\der_FL)_0^{2m+1} := \{d\in\der_FL\mid d(\frak g_\alpha\otimes
t^{2k})\subset
\frak s_\alpha\otimes t^{2k+2m+1}\\
\qquad \qquad \qquad \quad \text{and}\ d(\frak s_\beta\otimes
t^{2k+1}) \subset \frak g_\beta\otimes t^{2k+2m+2} \ \text{for all
$\alpha\in\Delta^{\rm red}$, $\beta\in\Delta'$
 and $k\in\Bbb Z$}\}.
 \end{array}
 $$

\lemma\label{oddcomshift} Let $q\in (\der_F L)_0^{2m+1}$. Then,
$q$ commutes with a shift map $s_{2i}$ for all $i\in\Bbb Z$.
\endlemma

\proof We note that
$$q(x_\alpha \otimes t^{2k}) = 0\quad (x_\alpha \in \frak g_\alpha,\ \alpha \in \Delta_{lg},\ k \in \mathbb{Z})$$
for $\text B_\frak I$ or $\text C_\frak I$, and that
$$q(x_\beta \otimes t^{2k+1}) = 0\quad (x_\beta \in \frak s_\beta,\ \beta \in \Delta_{ex},\ k \in \mathbb{Z})$$
for $\text{BC}_\frak I$, since $\frak s_\alpha = 0$ and $\frak
g_\beta = 0$. Therefore, in particular,
$$q \circ s_{2i} (z) = s_{2i} \circ q (z)$$
for $z = x_\alpha \otimes t^{2k}$ in the case of type $\text
B_\frak I$ or $\text C_\frak I$, and for $z = x_\beta \otimes
t^{2k+1}$ in the case of type $\text{BC}_\frak I$. For any other
given homogeneous element $x$ we can find suitable homogeneous
elements $y$ and $z$ such that $x = [ y , z ]$ in the following
sense.

$$\begin{array}{c||c|c|c}
& x & y & z\\ \hline \hline
\text B_\frak I\ \text{or}\ \text C_\frak I & x_\alpha \otimes t^{2k} & y_{\alpha'} \otimes t^{2k} & z_{\alpha''} \otimes 1\\
& x_\alpha \in \frak g_\alpha\ (\alpha \in \Delta_{sh}) &
y_{\alpha'} \in \frak g_{\alpha'}\ (\alpha' \in \Delta_{sh}) &
z_{\alpha''} \in \frak g_{\alpha''}\ (\alpha'' \in \Delta_{lg})\\
\cline{2-4}
& x_\beta \otimes t^{2k+1} & y_{\beta'} \otimes t^{2k+1} & z_{\alpha''} \otimes 1\\
& x_\beta \in \frak s_\beta\ (\beta \in \Delta_{sh}) & y_{\beta'}
\in \frak s_{\beta'}\ (\beta' \in \Delta_{sh}) & z_{\alpha''} \in
\frak g_{\alpha''}\ (\alpha'' \in \Delta_{lg})\\ \hline
\end{array}$$

\medskip

\noindent The table shown above is for $\text B_\frak I$ or $\text
C_\frak I$. For example, we can understand that for any $x =
x_\alpha \otimes t^{2k}\ (x_\alpha \in \frak g_\alpha,\ \alpha \in
\Delta_{sh})$, there are $y = y_{\alpha'} \otimes t^{2k}\
(y_{\alpha'} \in \frak g_{\alpha'}\, \alpha' \in \Delta_{sh})$ and
$z = z_{\alpha''} \otimes 1\ (z_{\alpha''} \in \frak
g_{\alpha''},\ \alpha'' \in \Delta_{lg})$ such that $x = [ y , z
]$.

Similarly, for $\text{BC}_\frak I$, we obtain the following table.

$$\begin{array}{c||c|c|c}
& x & y & z\\ \hline \hline
\ \ \ \text {BC}_\frak I\ \  \ & x_\alpha \otimes t^{2k} & y_{\beta'} \otimes t^{2k-1} & z_{\beta''} \otimes t\\
& x_\alpha \in \frak g_\alpha\ (\alpha \in \Delta_{sh}) &
y_{\beta'} \in \frak s_{\beta'}\ (\beta' \in \Delta_{sh}) &
z_{\beta''} \in \frak s_{\beta''}\ (\beta'' \in \Delta_{ex})\\
\cline{2-4}
& x_\alpha \otimes t^{2k} & y_{\alpha'} \otimes t^{2k-1} & z_{\beta''} \otimes t\\
& x_\alpha \in \frak g_\alpha\ (\alpha \in \Delta_{lg}) &
y_{\beta'} \in \frak s_{\beta'}\ (\alpha' \in \Delta_{lg}) &
z_{\alpha''} \in \frak s_{\beta''}\ (\beta'' \in \Delta_{ex})\\
\cline{2-4}
& x_\beta \otimes t^{2k+1} & y_{\beta'} \otimes t^{2k} & z_{\beta''} \otimes t\\
& x_\beta \in \frak s_\beta\ (\beta \in \Delta_{sh}) & y_{\alpha'}
\in \frak g_{\alpha'}\ (\alpha' \in \Delta_{sh}) & z_{\beta''} \in
\frak s_{\beta''}\ (\beta'' \in \Delta_{ex})\\ \cline{2-4}
& x_\beta \otimes t^{2k+1} & y_{\alpha'} \otimes t^{2k} & z_{\beta''} \otimes t\\
& x_\beta \in \frak s_\beta\ (\beta \in \Delta_{lg}) & y_{\alpha'}
\in \frak g_{\alpha'}\ (\alpha' \in \Delta_{lg}) & z_{\alpha''}
\in \frak s_{\beta''}\ (\beta'' \in \Delta_{ex})\\ \hline
\end{array}$$

\bigskip

\noindent In the expression $x = [ y , z ]$, we note that $q(z) =
0$ is always true for all $\text B_\frak I$, $\text C_\frak I$ and
$\text{BC}_\frak I$ as before, which is the most important fact in
this case. Hence, we obtain
$$\begin{array}{lll}
q \circ s_{2i}(x) & = & q \circ s_{2i} ([y,z]) = q([y,s_{2i}(z)])\\
& = & [q(y),s_{2i}(z)] + [y,q\circ s_{2i}(z)]\\
& = & [q(y),s_{2i}(z)]\\
& = & s_{2i}([q(y),z])\\
& = & s_{2i}([q(y),z] + [y,q(z)])\\
& = & s_{2i} \circ q([y,z])\\
& = & s_{2i} \circ q(x).
\end{array}$$
Therefore, $q \circ s_{2i} = s_{2i} \circ q$ on $L$. \qed

\enskip

\lemma\label{extder} Let $L=\big( \frak g\otimes F[t^{\pm 2}]
\big) \oplus \big( \frak s\otimes tF[t^{\pm 2}] \big)$ be a
twisted loop algebra, which is double graded by $\Delta\cup\{0\}$
and $\Bbb Z$ as above. Let $d$ be in $(\der_F L)_0^{2m+1}$ such
that $s_2\circ d=d\circ s_2$. Then, a unique derivation $\tilde d$
on $\tilde L$ exists such that
$$
\tilde d\mid_{L}=d,\quad \tilde d(x\otimes t^{2k+1})=s_1\circ
d(x\otimes t^{2k}) \quad\text{and}\quad \tilde d(v\otimes
t^{2k})=s_{-1}\circ d(v\otimes t^{2k+1})
$$
for all $x\in\frak g$,\  $v\in \frak s$, and $k\in\Bbb Z$.
Moreover,
$$\tilde d
\in(\der_F \tilde L)_0^{2m+1}\quad \text{such that}\quad s_k\circ
\tilde d=\tilde d\circ s_k \quad \text{for all $k\in\Bbb Z$}.$$

\endlemma
\proof The uniqueness is clear since the image of all the
homogeneous elements has been determined. Therefore, it is
sufficient to show that $\tilde d$ is a derivation. Thus, we need
to check the following:\quad For $x,y\in\frak g$ and $v,w\in \frak
s$,
\begin{itemize}
\item[(a)] $\tilde d([x\otimes t^{2k}, y\otimes t^{2\ell+1}])
=[\tilde d(x\otimes t^{2k}), y\otimes t^{2\ell+1}] +[x\otimes
t^{2k}, \tilde d(y\otimes t^{2\ell+1})]$

\item[(b)] $\tilde d([x\otimes t^{2k}, v\otimes t^{2\ell}])
=[\tilde d(x\otimes t^{2k}), v\otimes t^{2\ell}] +[x\otimes
t^{2k}, \tilde d(v\otimes t^{2\ell})]$

\item[(c)] $\tilde d([x\otimes t^{2k+1}, y\otimes t^{2\ell+1}])
=[\tilde d(x\otimes t^{2k+1}), y\otimes t^{2\ell+1}] +[x\otimes
t^{2k+1}, \tilde d(y\otimes t^{2\ell+1})]$

\item[(d)] $\tilde d([x\otimes t^{2k+1}, v\otimes t^{2\ell+1}])
=[\tilde d(x\otimes t^{2k+1}), v\otimes t^{2\ell+1}] +[x\otimes
t^{2k+1}, \tilde d(v\otimes t^{2\ell+1})]$

\item[(e)] $\tilde d([x\otimes t^{2k+1}, v\otimes t^{2\ell}])
=[\tilde d(x\otimes t^{2k+1}), v\otimes t^{2\ell}] +[x\otimes
t^{2k+1}, \tilde d(v\otimes t^{2\ell})]$

\item[(f)] $\tilde d([v\otimes t^{2k+1}, w\otimes t^{2\ell}])
=[\tilde d(v\otimes t^{2k+1}), w\otimes t^{2\ell}] +[v\otimes
t^{2k+1}, \tilde d(w\otimes t^{2\ell})]$

\item[(g)] $\tilde d([v\otimes t^{2k}, w\otimes t^{2\ell}])
=[\tilde d(v\otimes t^{2k}), w\otimes t^{2\ell}] +[v\otimes
t^{2k}, \tilde d(w\otimes t^{2\ell})]$.

\end{itemize}

All of these equations involve simple calculations, but we check
them to be sure.

For (a), we have
\begin{align*}
(LHS) &=\tilde d([x,y]\otimes t^{2k+2\ell+1}) =s_1\circ
d([x,y]\otimes t^{2k+2\ell})
=s_1\circ d([x\otimes t^{2k}, y\otimes t^{2\ell}])\\
&=s_1 \big([d(x\otimes t^{2k}), y\otimes t^{2\ell}]
+[x\otimes t^{2k}, d(y\otimes t^{2\ell})]\big)\\
&=[d(x\otimes t^{2k}), y\otimes t^{2\ell+1}] +[x\otimes t^{2k},
s_1\circ d(y\otimes t^{2\ell})]=(RHS).
\end{align*}

For (b), we have
\begin{align*}
(LHS) & =\tilde d([x,v]\otimes t^{2k+2\ell})
=s_{-1}\circ d([x,v]\otimes t^{2k+2\ell+1}])\\
&=s_{-1} \big([d(x\otimes t^{2k}), v\otimes t^{2\ell+1}]
+[x\otimes t^{2k}, d(v\otimes t^{2\ell+1})]\big)\\
&=[d(x\otimes t^{2k}), v\otimes t^{2\ell}] +[x\otimes t^{2k},
s_{-1} \circ d(v\otimes t^{2\ell+1})]=(RHS).
\end{align*}

For (c), we have
\begin{align*}
(LHS) & =\tilde d([x,y]\otimes t^{2k+2\ell+2})
=d([x,y]\otimes t^{2k+2\ell+2}])\\
&=d([x\otimes t^{2k},y\otimes t^{2\ell+2}])
=[d(x\otimes t^{2k}),y\otimes t^{2\ell+2}]+[x\otimes t^{2k},d(y\otimes t^{2\ell+2})]\\
&=s_1 \big( [d(x\otimes t^{2k}),y\otimes t^{2\ell+1}] \big)
+[x\otimes t^{2k},d\circ s_2(y\otimes t^{2\ell})]\\
&= [s_1\circ d(x\otimes t^{2k}),y\otimes t^{2\ell+1}]
+s_2 \big( [x\otimes t^{2k},d(y\otimes t^{2\ell})] \big)\\
&\text{(since $s_2$ and $d$ commute)}\\
&= [\tilde d(x\otimes t^{2k+1}),y\otimes t^{2\ell+1}] +[x\otimes
t^{2k+1},s_1\circ d(y\otimes t^{2\ell})]=(RHS).
\end{align*}

For (d), we have
\begin{align*}
(LHS) & =\tilde d([x,v]\otimes t^{2k+2\ell+2})
=s_{-1}\circ d([x,v]\otimes t^{2\ell+3}])\\
&=s_{-1}\circ d([x\otimes t^{2k}, v\otimes t^{2\ell+3}])\\
&=s_{-1} \big([d(x\otimes t^{2k}),v\otimes t^{2\ell+3}]
+[x\otimes t^{2k},d(v\otimes t^{2\ell+3})]\big)\\
&=s_{1}\big( [d(x\otimes t^{2k}),v\otimes t^{2\ell+1}]\big)
+s_{-2}\big([x\otimes t^{2k+1},d(v\otimes t^{2\ell+3})]\big)\\
&=[s_{1}\circ d(x\otimes t^{2k}),v\otimes t^{2\ell+1}]
+[x\otimes t^{2k+1},s_{-2}\circ d(v\otimes t^{2\ell+3})]\\
&=(RHS)\quad \text{(since $s_2$ and $d$ commute)}.
\end{align*}

For (e), we have
\begin{align*}
(LHS) & =\tilde d([x,v]\otimes t^{2k+2\ell+1}) =d([x,v]\otimes
t^{2k+2\ell+1}])
= d([x\otimes t^{2k}, v\otimes t^{2\ell+1}])\\
&=[d(x\otimes t^{2k}),v\otimes t^{2\ell+1}]
+[x\otimes t^{2k},d(v\otimes t^{2\ell+1})]\\
&=[s_{1}\circ d(x\otimes t^{2k}),v\otimes t^{2\ell}] +[x\otimes
t^{2k},\tilde d(v\otimes t^{2\ell+1})] =(RHS).
\end{align*}

For (f), we have
\begin{align*}
(LHS) & =\tilde d([v, w]\otimes t^{2k+2\ell+1}])
=s_{1}\circ d([v,w]\otimes t^{2k+2\ell}])\\
&=s_{1}\circ d([v\otimes t^{2k-1}, w\otimes t^{2\ell+1}])\\
&=s_{1} \big([d(v\otimes t^{2k-1}),w\otimes t^{2\ell+1}]
+[v\otimes t^{2k-1},d(w\otimes t^{2\ell+1})]\big)\\
&= s_2 [ d(v\otimes t^{2k-1}),w\otimes t^{2\ell}]\big)
+[v\otimes t^{2k},d(w\otimes t^{2\ell+1})]\\
&=[d(v\otimes t^{2k+1}),w\otimes t^{2\ell}] +s_{-1} \big(
[v\otimes t^{2k+1},d(w\otimes t^{2\ell+1})] \big) =(RHS).
\end{align*}

For (g), we have
\begin{align*}
(LHS) & =\tilde d([v, w]\otimes t^{2k+2\ell}]) =d([v,w]\otimes
t^{2k+2\ell}])
=d([v\otimes t^{2k-1}, w\otimes t^{2\ell+1}])\\
&=[d(v\otimes t^{2k-1}),w\otimes t^{2\ell+1}]
+[v\otimes t^{2k-1},d(w\otimes t^{2\ell+1})]\\
&=[s_{1}\circ d(v\otimes t^{2k-1}),w\otimes t^{2\ell}] +[v\otimes
t^{2k},s_{-1}\circ d(w\otimes t^{2\ell+1})] =(RHS).
\end{align*}

For the second assertion, it is clear that $\tilde d\in(\der_F
\tilde L)_0^{2m+1}$. In addition, since $d$ commutes with $s_{2}$,
then the same is true of $\tilde d$. Hence, by Lemma
\ref{classifydegreeder}, $\tilde d$ commutes with $s_k$ for all
$k\in\Bbb Z$. \qed

\enskip

Thus, together with Lemma \ref{diagonaldert}, we have classified
the diagonal derivations of twisted locally loop algebras.

\theorem\label{odddegreeth} Let $L$ be a twisted loop algebra.
Then, we have $ (\der_FL)^0_0 =\ad P\oplus Fd^{(0)}$, where $P$ is
defined in \eqref{defP}, and
$$(\der_F L)_0^{2m}=s_{2m}\circ(\der_F L)_0^0
\quad\text{and}\quad (\der_F L)_0^{2m+1}=s_{2m+1}\circ \ad T^-
$$
for all $m\in\Bbb Z$, where $T^-=\frak s_0$\ \ for\ \ $\text
{B}_\frak I^{(2)}$,\ \ $T^-=T^-_{2\frak I}$\ \ for\ \ $\text
C_\frak I^{(2)}$\ \ or\ \ $T^-=T^-_{2\frak I+1}$\ \ for\ \ $\text
{BC}_\frak I^{(2)}$, as defined in Example \ref{exampleLALA}.
\endtheorem
\proof By Lemma \ref{oddcomshift}, \ref{extder}, and the
classification of the untwisted case, if $d\in (\der_F
L)_0^{2m+1}$, then $\tilde d\in s_{2m+1}\circ(\der_F \tilde
L)_0^0$. In addition, by Lemma \ref{oddcomshift} and Lemma
\ref{classifydegreedert}, we obtain $\tilde d\in s_{2m+1}\circ\ad
P$. Thus, $\ad p:=s_{-2m-1}\circ \tilde d\in\ad P$, and we have
$[p,\frak g^+]\subset \frak g^-$ according to the terminology used
in Lemma \ref{-eleofh}. Hence, by Lemma \ref{-eleofh}, we obtain
$p\in T^-$. Therefore, $d\in s_{2m+1}\circ\ad T^-$. \qed

\remark If $L$ is a twisted loop algebra of type $\text {B}_\frak
I^{(2)}$, then $(\der_F L)_0^{2m+1}=s_{2m+1}\circ \ad \frak s_0
=\ad(\frak s_0\otimes t^{2m+1})$. Thus, there is no outer
derivation of odd degree.

\endremark

\medskip

We return to the classification of twisted LALAs. By Theorem
\ref{odddegreeth}, if $d\in D^{2m+1}$, then $\ad_L d\in
s_{2m+1}\circ\ad_LT^-$. The bracket on $D:=\bigoplus_{m\in\Bbb
Z}D^m$ can be investigated in the same manner as the untwisted
case. In particular, for type $\text {BC}_\frak I^{(2)}$ or $\text
{C}_\frak I^{(2)}$, we use the isomorphism
$$\mathcal L\cong
\big(\mathcal L+\iota\otimes F[t^{\pm 1}]\big) /\iota\otimes
tF[t^{\pm 2}].$$ 
Thus, $D^m$ for $m\in\Bbb Z$ is an exact example
for each type described in Example \ref{exampleLALA}. Thus, we
have completed the classification.

\theorem\label{classitw} Let $\mathcal L$ be a twisted LALA. Then,
$\mathcal L$ is isomorphic to that in Example \ref{exampleLALA}.
\qed
\endtheorem

\remark\label{twlalaofuntwlala} We can show that any twisted LALA
is the fixed algebra of some untwisted LALA. Moreover, for any
untwisted LALA $\mathcal L$ of type $\text A_\frak I^{(1)}$ or
$\text D_\frak I^{(1)}$, a twisted LALA $\mathcal L'$ exists,
which is a subalgebra of $\mathcal L$ such that $\mathcal L'$ is
the intersection of $\mathcal L$ and the fixed algebra of a
maximal untwisted LALA $\mathcal L^{max}$ that contains $\mathcal
L$. Note that a maximal twisted LALA is also unique up to
isomorphism, as in the case of a maximal untwisted LALA.
\endremark

\remark\label{onlytypeA} By Theorem \ref{classiuntw} and Theorem
\ref{classitw}, the LALAs in Example \ref{exampleLALA} comprise
all of the algebras. Given this fact, the following statement is
clear and it is a useful criterion.

If a diagonal matrix $p\in T$ with a trace that is a nonzero value
(e.g., $e_{ii}$ or $e_{ii}+e_{\frak I+i,\frak I+i}$, etc.) is used
in a LALA, then this LALA must be of type $\text A_\frak I^{(1)}$,
$\text C_\frak I^{(2)}$, or $\text {BC}_\frak I^{(2)}$. Moreover,
if the type is $\text C_\frak I^{(2)}$ or $\text {BC}_\frak
I^{(2)}$, then $p$ has to be used in odd degree.
\endremark

\enskip

\section{Standard LALAs}
\label{sec:}

We prove the following criterion whether a LALA is standard or not.
\begin{lemma}\label{cridegderi}
Let $(\mathcal L,\mathcal H, \mathcal B)$ be a LALA with center
$Fc$ and $\mathcal L_c$ is its core, which is a locally Lie
$1$-torus with grading pair $(\frak g,\frak h)$. If $0\neq
d\in\mathcal L$ exists such that $[d,\frak g]=0$ and $\mathcal
B(d,c)\neq 0$, then the action of $d$ on the $\Bbb Z$-graded core
coincides with a nonzero multiple of a degree derivation relative
to $\Bbb Z$, and thus $\mathcal L$ contains the degree derivation.
Hence, $\mathcal L$ is standard.
\end{lemma}
\proof Let $d=\sum_{\xi\in R}\ x_\xi$ for $x_\xi\in\mathcal
L_\xi$. If $\xi\in R^\times$, then $[\frak h,x_\xi] = Fx_\xi
\subset\mathcal L_\xi$, and thus $x_\xi=0$ since $[d,\frak g]=0$.
If $\xi\in R^0\setminus\{0\}$, then $x_\xi\in T\otimes t^{m}$ for
some $0\neq m\in\Bbb Z$, by Theorem \ref{classiuntw} and
\ref{classitw}. However, if $x_\xi\neq 0$, then a root vector
$y\in\frak g_\alpha$ ($\alpha\in\Delta$) exists such that
$[y,x_\xi]\neq0$, which is a contradiction. Hence, $x_\xi=0$.
Thus, $d=x_0\in\mathcal L_0=\mathcal H$. Then, by Theorems
\ref{classiuntw} and \ref{classitw},
$$d=p+a d^{(0)}+b c$$ for some
$p\in T = T \otimes t^0$ and $a,b\in F$, as well as $a\neq 0$,
since $\mathcal B(d,c)\neq 0$. Therefore, we have $0=[d,\frak
g]=[p,\frak g]$. However, unless $\mathcal L$ has type $\text
A_\frak I^{(1)}$, we have $p\in T^\sigma$, and thus $p$ must be
zero. If $\mathcal L$ has type $\text A_\frak I^{(1)}$, then $p\in
F\iota$, by Lemma \ref{centginA},
 and thus $p$ must again be zero (modulo $F\iota$).
 Thus, we obtain $d=a d^{(0)}+bc$.
\qed

\enskip

\remark\label{neebsmin} In 
\cite[Def.3.6]{N2}, Neeb defined a
minimal LALA $\mathcal L$, which is minimal in the sense described
above and that satisfies one more condition:
$$
\text{ $\exists\ d\in\mathcal H$ such that\ \ $W':=\spa_\Bbb
Q\{\alpha\in R^\times\mid \alpha(d)=0\}$\ \ is a reflectable
section}
$$
of $W=\spa_\Bbb Q R^\times$. Thus, $[\frak g,d]=0$. Moreover, if
$\delta(d)= 0$, where $\delta$ is a generator of $R^0\cong\Bbb Z$,
then $\alpha(d)=0$ for all $\alpha\in R^\times$.
Hence, $W'=W$, which is a contradiction. Thus, $\delta(d)\neq 0$.
However, $d$ is a nonzero multiple of a degree derivation modulo
of the center by Lemma \ref{cridegderi}, and thus a minimal LALA
in \cite{N2} is a minimal standard LALA in our sense.
\endremark

\example\label{nonstaiso} The minimal LALA $\mathcal L =\sll_\Bbb
N(F[t^{\pm 1}])\oplus Fc \oplus F(e_{11}+d^{(0)})$ is isomorphic
to a minimal standard LALA $\mathcal L^{ms}=\sll_\Bbb N(F[t^{\pm
1}])\oplus Fc \oplus Fd^{(0)}$. In fact, let $g= \dia  (t, 1, 1,
\ldots )$. Then, $g^{-1}Xg$ for $X\in \sll_\Bbb N(F[t^{\pm 1}])$
gives an automorphism $f$ of $\sll_\Bbb N(F[t^{\pm 1}])$.
Therefore, we can extend $f$ from $\mathcal L^{ms}$ onto $\mathcal
L$ such that $f(c)=c$ and $f(d^{(0)})=e_{11}+d^{(0)}$. Thus,
$\mathcal L$ is isomorphic to $\mathcal L^{ms}$, as in Lie
algebras.
\endexample

\example \label{counterexofst} Let $\displaystyle{ p= \dia  (1,
\frac{1}{2}, \frac{1}{3},  \frac{1}{4}, \ldots )} $ and put
$d=p+d^{(0)}$. Then, the minimal LALA $\mathcal L =\sll_\Bbb
N(F[t^{\pm 1}])\oplus Fc \oplus Fd$ is not isomorphic to a minimal
standard LALA $\mathcal L^{ms}$. In fact, if $\mathcal L$ is isomorphic to
$\mathcal L^{ms}$, then an isomorphism
$$\psi : \mathcal L^{ms} \longrightarrow \mathcal L$$
exists such that $\psi(d^{(0)}) = x + a d = x + a (d^{(0)} + p)$
for some $x \in \mathcal L_c = \sll_\Bbb N(F[t^{\pm 1}])\oplus Fc$
and some nonzero $a \in F$. Then, we have
$$\psi \circ \ad d^{(0)} \circ\psi^{-1} =  \ad(\psi(d^{(0)})) =
\ad (x + a  d^{(0)}+ a  p)$$ in $\der_F(\mathcal L)$. Now, we can
compare the eigenvalues of the same operators $\psi \circ \ad
d^{(0)} \circ\psi^{-1}$ and $\ad (x + a   d^{(0)}+ a   p)$. Note
that the eigenvalues of $\psi \circ \ad d^{(0)} \circ\psi^{-1}$
are all integers. We can select $h = e_{\ell \ell} - e_{\ell + 1,
\ell + 1} \in \sll_\Bbb N(F[t^{\pm 1}])$ such that
$$[x,h] = 0,$$
by taking $\ell > \!\! > 0$, where $e_{ij}$ is a matrix unit.
Then,
$$[x + a  d^{(0)} + a  p,h \otimes t] = a(h \otimes t),$$
which implies that $a$ is a nonzero integer since $a$ is an
eigenvalue of $\ad (x + a   d^{(0)}+ a   p)$. We can also choose
sufficiently large different integers $m\neq n >> 0$ that satisfy
$$[x,e_{mn}] = 0$$
and
\begin{equation}\label{noninteger}
\frac{a(n-m)}{mn} \not\in \mathbb{Z}.
\end{equation}
For these integers, $m$ and $n$, we can see that
$$[x + a  d^{(0)} + a  p,e_{mn}] = a \left( \frac{1}{m} - \frac{1}{n} \right) e_{mn}
= \frac{a(n-m)}{mn} e_{mn}.$$ Since
$$\frac{a(n-m)}{mn}$$
is an eigenvalue of $\ad (x + a   d^{(0)}+ a   p)$, it must
be an integer, which contradicts \eqref{noninteger}. Hence,
$\mathcal L$ is not isomorphic to $\mathcal L^{ms}$.

\enskip

\enskip

\end{document}